\documentclass[reqno,10pt]{amsart}
\usepackage{amsfonts}
\usepackage{graphicx}
\usepackage{amsfonts,amsmath,amsthm,amssymb}
\usepackage{hyperref}
\usepackage{marginnote}
\usepackage{color}
\usepackage{mathrsfs}

\numberwithin{equation}{section}

\newcommand{\R}{\mathbb{R}}

\newtheorem{theorem}{Theorem}[section]
\newtheorem{corollary}[theorem]{Corollary}
\newtheorem{lemma}[theorem]{Lemma}
\newtheorem{proposition}[theorem]{Proposition}
\newtheorem{remark}[theorem]{Remark}
\newtheorem{definition}[theorem]{Definition}

\def\v{\varepsilon}

\def\t{\theta}

\def\m{\mu}

\def\g{\gamma}
\def\d{\delta}

\def\r{\rho}
\def\s{\sigma}

\def\f{\frac}
\def\dd{{\rm d}}
\def\M{{\mathcal{M}}}
\def\hl{{\mathfrak{b}}}
\usepackage{tikz}

\begin{document}

\title[Global Solutions to the Compressible Euler Equations]{Global Solutions
of the Compressible Euler Equations with Large Initial Data of Spherical Symmetry and
Positive Far-Field Density}

\author[G.-Q. Chen]{Gui-Qiang G. Chen}
\address
{Gui-Qiang G. Chen, Mathematical Institute, University of Oxford,
Oxford OX2 6GG, UK; Academy of Mathematics and Systems Science,
Chinese Academy of Sciences, Beijing 100190, China.}
\email{chengq@maths.ox.ac.uk}

\author[Y. Wang]{Yong Wang}
\address
{Yong Wang,
Academy of Mathematics and Systems Science and University of Chinese Academy of Sciences,
Chinese Academy of Sciences, Beijing 100190, China.}
\email{yongwang@amss.ac.cn}

\begin{abstract}
We are concerned with the global existence theory for spherically symmetric solutions
of the multidimensional
compressible Euler equations with large initial data of positive far-field density
so that the total initial-energy is  unbounded.
The central feature of the solutions is the strengthening of waves as they move radially inward toward the origin.
For the large initial data of positive far-field density,
various examples have shown that the spherically symmetric solutions
of the Euler equations blow up near the origin at certain time.
A fundamental unsolved problem is whether the density of the global solution would form concentration
to become a measure near the origin for the case when the total initial-energy is unbounded and
the wave propagation is not at finite speed starting initially.
In this paper, we establish a global existence theory for spherically symmetric solutions
of the compressible Euler equations with large initial data of
positive far-field density and relative finite-energy.
This is achieved by developing a new approach via adapting a class of degenerate density-dependent viscosity
terms, so that a rigorous proof of the vanishing viscosity
limit of global weak solutions of the Navier-Stokes equations
with the density-dependent viscosity terms to the corresponding global solution
of the Euler equations with large initial
data of spherical symmetry and positive far-field density
can be obtained.
One of our main observations is that
the adapted class of degenerate density-dependent viscosity terms
not only includes the viscosity terms for the Navier-Stokes
equations for shallow water (Saint Venant) flows
but also, more importantly, is suitable to achieve our key objective of this paper.
These results indicate that concentration is not formed in the vanishing viscosity limit
for the Navier-Stokes approximations constructed in this paper
even when the total initial-energy is unbounded,
though the density may blow up near the origin at certain time
and the wave propagation is not at finite speed.
\end{abstract}

\keywords{Euler equations, Navier-Stokes equations, compressible flows, spherically symmetric solutions,
positive far-field density,
unbounded energy, relative finite-energy, large data, {\it a priori} estimate,
higher integrability, concentration, blow up, inviscid limit, compactness framework, approximate solutions}
\subjclass[2010]{\, 35Q35, 35Q31, 35B25, 35B44, 35L65, 35L67, 76N10}
\date{\today}
\maketitle

\setcounter{tocdepth}{1}
\thispagestyle{empty}

\section{Introduction}

We are concerned with the global existence theory for spherically symmetric solutions
of the multidimensional (M-D)
compressible Euler equations with large initial data of positive far-field density,
{\it i.e.}, given constant density $\bar{\rho}>0$ at infinity, so that the total initial-energy is unbounded.
The study of spherically symmetric
solutions dates back to the 1950s and is motivated by many important physical
problems such as flow in a jet engine inlet manifold and stellar dynamics including
gaseous stars and supernovae formation ({\it cf.} \cite{Courant-Friedrichs,Guderley,Rosseland,Slemrod,Whitham}).
The central feature of the solutions is the strengthening of waves as they move radially inward toward
the origin.
An existence theory has been established in Chen-Perepelitsa  \cite{Chen7}
and Chen-Schrecker \cite{Chen-Schrecker} via an approach of vanishing artificial viscosity
for the case when the initial data are of finite-energy, which requires that $\bar{\rho}=0$.
For the far-field density $\bar{\rho}>0$, various physical examples have shown that the spherically
symmetric solutions of the compressible Euler equations blow up more often near the origin at certain
time (see \cite{Courant-Friedrichs,Guderley,Whitham} and
the references cited therein). The fundamental unsolved problem is whether
the density would form concentration to become a measure near the origin for the case
when the total initial-energy is unbounded and the wave propagation is not at finite speed starting initially.
In this paper, we establish a global existence theory for spherically symmetric solutions
in $L^p_{\rm loc}$ of the compressible Euler equations with large initial data of
positive far-field density $\bar{\rho}>0$ and relative finite-energy in $\mathbb{R}^N$ for $N\geqq2$.
This is achieved by developing a new approach via adapting a class of degenerate density-dependent viscosity
terms, so that a rigorous proof of the vanishing viscosity
limit of global weak solutions of the compressible Navier-Stokes equations
with the density-dependent viscosity terms to the corresponding global solution
of the Euler equations with large initial
data of spherical symmetry and positive far-field density can be obtained.
One of our main observations is that the adapted class of
degenerate density-dependent viscosity terms
not only includes the viscosity terms for the Navier-Stokes
equations for shallow water (Saint Venant) flows,
among others ({\it cf.} Bresch-Dejardins \cite{BD-2003-CRMASP,BD-2004-CRMASP},
 Bresch-Dejardins-Lin \cite{BD-2003-CPDE}, Lions \cite{Lions-CNS}, and Mallet-Vasseur \cite{MV}),
but also, more importantly, is suitable to achieve our key objective
of this paper.
These results indicate that concentration is not formed in the vanishing viscosity limit
for the Navier-Stokes approximations constructed in this paper
even when the total initial-energy is unbounded,
though the density may blow up near the origin at certain time and the wave propagation is not at finite speed.

More precisely, the M-D Euler equations for compressible isentropic fluids take the form:
\begin{align}\label{1.1-1}
\begin{cases}
\partial_t \rho+\mbox{div}\M=0,\\[1mm]
\partial_t \M+\mbox{div}\big(\frac{\M\otimes\M}{\rho}\big)+\nabla p=0,
\end{cases}
\end{align}
for $(t, \mathbf{x})\in \R_+\times \mathbb{R}^N$ with $N\geqq 2$,
where  $\rho$ is the density, $p$ is the pressure,
and $\M \in\mathbb{R}^N$ represents the momentum; see
also Chen-Feldman \cite{Chen-Feldman2018} and Dafermos \cite{Dafermos}.
When $\rho>0$, $U=\frac{\M}{\rho}\in\mathbb{R}^N$ is the velocity.
The constitutive pressure-density relation for polytropic gases is
\begin{equation*}
p=p(\r)=\kappa\r^{\gamma},
\end{equation*}
where $\gamma>1$ is the adiabatic exponent; by scaling, constant $\kappa$ in the pressure-density relation may be
chosen as  $\kappa=\frac{(\gamma-1)^2}{4\gamma}$ without loss
of generality.
We are concerned with the Cauchy problem for \eqref{1.1-1} with the Cauchy
data:
\begin{equation}\label{1.1-2}
(\rho, \M)|_{t=0}=(\rho_0, \M_0)(\mathbf{x})\longrightarrow (\bar{\rho}, \mathbf{0}) \qquad\,\,\,\,\mbox{as $|\mathbf{x}|\to \infty$},
\end{equation}
where $(\bar{\rho}, \mathbf{0})$ is a constant far-field state, for which
the initial far-field velocity has been assumed to be zero in \eqref{1.1-2}
without loss of generality, owing to the Galilean invariance of system \eqref{1.1-1}.
Since a global solution of the Euler equations \eqref{1.1-1} normally contains the vacuum states
$\{(t,\mathbf{x})\,:\, \rho(t,\mathbf{x})=0\}$ where
the fluid velocity $U(t,\mathbf{x})$ is not well-defined (even though the far-field density is positive),
we will use the physical variables such as the momentum
$\M(t,\mathbf{x})$, or $\frac{\M(t,\mathbf{x})}{\sqrt{\rho(t,\mathbf{x})}}$,
which will be shown to be always well-defined, instead of $U(t,\mathbf{x})$, when the vacuum states are involved
throughout this paper.

In order to construct global spherically symmetric solutions in $L^p_{\rm loc}$
of the Euler equations \eqref{1.1-1} with large initial data of positive far-field density, $\bar{\rho}>0$,
the approach of vanishing artificial viscosity developed in \cite{Chen7,Chen-Schrecker} is no longer applied directly,
and the problem has been remained open.
To solve this problem, in this paper, we develop a different approach by adapting a class of degenerate density-dependent viscosity
terms so that
the required uniform estimates in terms of the viscosity coefficients can be achieved for the vanishing
viscosity limit. More precisely, we consider the M-D Navier-Stokes equations for compressible barotropic fluids
with the adapted class of degenerate density-dependent viscosity terms:
\begin{align}\label{1.1}
\begin{cases}
\partial_t \rho+\mbox{div} \M=0,\\[1mm]
\partial_t \M+\mbox{div}\big(\frac{\M\otimes\M}{\rho}\big)+\nabla p=\v \mbox{div}\big(\mu(\r)D(\frac{\M}{\rho})\big)+\v\nabla\big(\lambda(\r)\mbox{div}(\frac{\M}{\rho})\big),
\end{cases}
\end{align}
where $D(\frac{\M}{\rho})=\frac{1}{2}\big(\nabla (\frac{\M}{\rho})+(\nabla (\frac{\M}{\rho}))^\top\big)$ is the stress tensor,
and the shear and bulk viscosity coefficients $\mu(\rho)$ and $\lambda(\rho)$ depend on the density
and may vanish on the vacuum.
Indeed,
in the derivation of the Navier-Stokes equations from the Boltzmann equation by the Chapman-Enskog expansions,
the viscosity terms depend on the temperature, which are translated into the dependence on
the density for barotropic flows ({\it cf}. \cite{liu-xin-yang}).
Moreover,
for the shallow water (Saint Venant) models, $N=2, \gamma=2$,
and $(\mu(\rho),\lambda(\rho))=(\rho, 0)$ ({\it cf.} Lions \cite[\S 8.4]{Lions-CNS});
also see \cite{BD-2003-CRMASP,BD-2003-CPDE}
for such models in geophysical flows.
This indicates that it is of independent interest and importance to
analyze the Navier-Stokes equations \eqref{1.1} with the density-dependent viscosity terms.
In particular, we are also interested in the inviscid limit of the Navier-Stokes equations \eqref{1.1}.
Formally, as $\v\rightarrow0+$, the Navier-Stokes equations \eqref{1.1} converge to the Euler
equations \eqref{1.1-1}.
A fundamental problem in mathematical fluid dynamics is whether a rigorous proof of
the vanishing viscosity limit of the solutions of the Navier-Stokes equations \eqref{1.1}
to the Euler equations \eqref{1.1-1} could be provided.

There has been an extensive literature in the analysis of the vanishing artificial/numerical viscosity limit
to the isentropic Euler equations.
For the 1-D case with general $L^\infty$ initial data,
it has been analyzed by DiPerna \cite{R.J.DiPerna1}, Ding-Chen-Luo \cite{Chen5}, Ding \cite{Ding},
Chen \cite{Chen1,Chen3}, Lions-Perthame-Souganidis \cite{Lions P.-L.1}, Lions-Perthame-Tadmor \cite{Lions P.-L.2},
and Huang-Wang \cite{Huang-Wang} via the methods of entropy analysis and compensated compactness.
Also see DiPerna \cite{R.J.DiPerna2}, Morawetz \cite{C. Morawetz}, Perthame-Tzavaras \cite{B. Perthame},
and Serre \cite{D. Serre1} for
general $2\times2$ strictly
hyperbolic systems of conservation laws.
The vanishing artificial viscosity limit to general strictly hyperbolic systems of conservation laws
with general small BV initial data was first established by Bianchini-Bressan \cite{Bressan-2005}
via direct BV estimates with small oscillation; see also \cite{Bressan-2004,Wang-2012-SIAM} and the references cited therein
for the rate of convergence.

For the study of spherically symmetric weak solutions,
the local existence of such solutions outside a solid ball at the origin was
discussed in Makino-Takeno \cite{M-T} for the case $1<\gamma\leqq \frac{5}{3}$;
also see Yang \cite{Yang-1,Yang-2}.
A first global existence of spherically symmetric solutions in $L^\infty$ including the origin was
established in Chen \cite{Chen-1997} for a class of $L^\infty$ Cauchy data of arbitrarily large amplitude,
which model outgoing blast waves and large-time asymptotic solutions.
A compactness framework was established
in LeFloch-Westdickenberg \cite{Ph. LeFloch} to construct finite-energy solutions
to the isentropic Euler equations with spherical symmetry and finite-energy initial data
for the case $1<\gamma\leqq \frac{5}{3}$.
As indicated earlier,  the convergence of the vanishing artificial viscosity
approximate solutions
to the corresponding finite-energy entropy solution of the M-D Euler equations
with large initial data of spherical symmetry was
established in \cite{Chen7,Chen-Schrecker} for any $\gamma>1$ for the case $\bar{\rho}=0$.

For the compressible Navier-Stokes equations with constant viscosity
coefficients ({\it i.e.}, $\mu$ and $\lambda$ are constants),
the global existence of solutions has been studied extensively;
see \cite{D. Hoff1,Kanel} and the references cited therein
for the 1-D case.
For $\mathbf{x}\in \R^N, N\ge 2$, Lions \cite{Lions-CNS} first obtained
the global existence of renormalized solutions, provided that $\gamma$ is suitably large,
which was further extended by Feireisl-Novotny-Petzeltov\'{a} \cite{Feireisl} to $\gamma>\frac{N}{2}$
and by Plotnikov-Weigant \cite{Plotnikov} to $\gamma=\frac{N}{2}$,
and by Jiang-Zhang \cite{Song Jiang} to $\gamma>1$
under the spherical symmetry.
When $\mu$ and $\lambda$ depend on the density, the Navier-Stokes equations \eqref{1.1}
become degenerate when $\rho\to 0$.
Such cases were analyzed in Bresch-Desjardins-Lin \cite{BD-2003-CPDE}
based on the new mathematical entropy -- the BD entropy,
first discovered by Bresch-Desjardins \cite{BD-2003-CRMASP} for the particular case $(\mu,\lambda)=(\rho,0)$,
and later
generalized by Bresch-Desjardins  \cite{BD-2004-CRMASP} to include the case of any viscosity
coefficients $(\mu, \lambda)$ satisfying the BD relation:
$\lambda(\rho)=\rho\mu'(\rho)-\mu(\rho)$;
also see Bresch-Desjardins \cite{BD-2007-JMPA}.
When the initial data are of spherical symmetry, Guo-Jiu-Xin \cite{Guo-Jiu-Xin-2} obtained
the global existence of spherically symmetric weak solutions of the system for $\gamma\in(1,3)$
in a finite ball with Dirichlet boundary conditions.

The idea of regarding inviscid gases as viscous gases with vanishing physical viscosity
can date back the seminal paper by Stokes \cite{Stokes} and the important contributions
of Rankine \cite{ Rankine}, Hugoniot \cite{Hugoniot},
and Rayleigh \cite{Lord Rayleigh} ({\it cf.} Dafermos \cite{Dafermos}).
However, the first rigorous convergence analysis of the inviscid limit
from the barotropic Navier-Stokes to Euler equations
was made by Gilbarg \cite{D. Gilbarg} much later, in which the
existence and vanishing viscous limit of the Navier-Stokes shock layers were established.
For the convergence analysis confined in the framework of piecewise smooth solutions,
see \cite{GMWZ,D. Hoff-Liu,Xin-1993}
and the references cited therein.

The key objective of this paper is to establish the global existence of
spherically symmetric solutions of \eqref{1.1-1}:
\begin{align}\label{1.2}
\r(t,\mathbf{x})=\r(t,r),\quad \M(t,\mathbf{x})=m(t,r)\f{\mathbf{x}}{r}\quad\,\,\,\,\,\, \mbox{for $r=|\mathbf{x}|$},
\end{align}
subject to the
initial condition:
\begin{align}\label{initial}
(\rho, \M)(0,\mathbf{x})=
(\rho_0,\M_0)(\mathbf{x})
=(\rho_0(r),m_0(r)\frac{\mathbf{x}}{r}) \longrightarrow (\bar{\rho},\mathbf{0})
\quad\,\,\mbox{as $r\to \infty$}
\end{align}
with $\bar{\rho}>0$ and relative finite-energy.
To achieve this, we establish the vanishing viscosity limit of the corresponding spherically symmetric
solutions of the Navier-Stokes equations \eqref{1.1} with the
adapted class of degenerate density-dependent viscosity terms and
approximate initial data of similar form to
\eqref{initial}.
For spherically symmetric solutions of form \eqref{1.2},  systems \eqref{1.1-1} and \eqref{1.1}
become
\begin{align}\label{euler}
\begin{cases}
\r_t+ m_r+\f{N-1}r m=0,\\[1mm]
m_t+\big(\frac{m^2}{\rho}+p\big)_r+\f{N-1}r \frac{m^2}{\rho}=0,
\end{cases}
\end{align}
and
\begin{align}\label{1.3}
\begin{cases}
\r_t+ m_r+\f{N-1}r m=0,\\[1mm]
m_t+\big(\frac{m^2}{\rho}+p\big)_r+\f{N-1}r \frac{m^2}{\rho}
 =\v\big((\mu+\lambda)((\frac{m}{\rho})_r+\frac{N-1}{r}\frac{m}{\rho})\big)_r-\v\frac{N-1}{r}\frac{m}{\rho}\mu_r,
\end{cases}
\end{align}
respectively.

In Chen-Perepelitsa \cite{Chen6},
the vanishing viscosity limit of smooth solutions for the 1-D
Navier-Stokes equations to the corresponding relative finite-energy solution
of the Euler equations has been established
for $\bar{\rho}>0$; also see \cite{Chen7a} for the 1-D shallow water case.
In \cite{Chen7,Chen-Schrecker}, the convergence of artificial viscosity approximate smooth solutions
to the corresponding finite-energy entropy solution of
the Euler equations \eqref{euler}
with spherical symmetry and large initial data
has been
established for $\bar{\rho}=0$ (also see \cite{Schrecker}).
As indicated earlier, in this paper, we develop a different approach to investigate the vanishing physical
viscosity limit of the weak solutions of the M-D Navier-Stokes equations \eqref{1.1}
with spherical symmetry to the corresponding relative finite-energy solution of the Euler equations \eqref{1.1-1}
with large initial data of positive far-field density $\bar{\rho}>0$.
Owing to the non-zero initial density at infinity so that the total initial-energy is unbounded,
which may cause the possibility for additional nature of singularities
at origin $r=0$ and far-field $r=\infty$,
several key techniques for the previous uniform
estimates as in \cite{Chen6,Chen7,Chen-Schrecker}
no longer apply.
In particular, for the weak solutions of the Navier-Stokes equations,
it is essential to ensure enough decay of solutions {\it a priori} as $r\rightarrow \infty$
so that integration by parts on unbounded regions
can be performed for the key estimates in the proof.

\smallskip
We now describe some of our approach and techniques involved to solve the problem in this paper.
Owing to the singularity at $r=0$, it has not been clear yet whether there always exists a global smooth solution
of the Cauchy problem of the Navier-Stokes equations with smooth large initial data of spherical symmetry.
To achieve our key objective, the main point of this paper is first to obtain
global weak solutions of the compressible Navier-Stokes equations
with some uniform estimates and the $H^{-1}_{\rm loc}$--compactness, so that
the compactness framework in \cite{Chen6} can be applied.
For this purpose, we first construct smooth approximate solutions
$(\rho^{\v,\d,b},m^{\v,\d,b})$, depending on the three parameters $(\v,\d,b)$,
through the Navier-Stokes equations \eqref{1.3};
see \eqref{6.1}--\eqref{6.3}.
Noting that the spherically symmetric Navier-Stokes equations \eqref{1.3} become
singular at the origin,
we first remove the origin in the approximate problem.
For the smooth approximate solutions as designed,
it is direct to obtain the basic energy estimate, Lemma \ref{lem6.1}.
Under relation \eqref{6.1-3}, we also obtain the
BD entropy estimate, Lemma \ref{lem6.2}.
Similar to that in \cite{Chen6}, we can obtain the uniform higher integrability
of the density; see Lemma \ref{lem6.3}.

To employ the compactness framework in \cite{Chen6},
we still need the uniform higher integrability of the velocity, as described in Proposition \ref{lem6.4},
for all $\gamma>1$.
To prove this, we apply the relative entropy pair $(\tilde{\eta},\tilde{q})$ of the spherically symmetric
Euler equations \eqref{euler}
to obtain \eqref{6.82} in \S 4.
The most difficult terms are the second and third terms on the right-hand side of \eqref{6.82},
which are essential for the M-D case (these two terms do not appear for the 1-D case).
By a careful analysis on the relative entropy pair, we see that
\begin{align}\label{1.80}
m\partial_\rho \tilde{\eta}(\rho,m)
 +\frac{m^2}{\rho}\partial_m \tilde{\eta}(\rho,m)-\tilde{q}(\rho,m)
\leqq C_\gamma(\bar{\rho})\big(\frac{m^2}{\rho}+ e(\rho,\bar{\rho})\big)
\end{align}
for some constant $C_\gamma(\bar{\rho})>0$,
which implies that the third term on the right-hand side of \eqref{6.82}
can be bounded by using the basic energy at least locally;
see Lemma \ref{lem6.14} for the details.
In fact, estimate \eqref{1.80} is quite subtle.
Since the left-hand side of \eqref{1.80} contains the terms on $\frac{|m|^3}{\rho^2}$ and $\rho^{\gamma+\theta}$,
we have to deal with such terms; otherwise, the higher integrability of the velocity may not be obtained.
This is achieved by our observation of underlying cancellation by dividing it into several cases;
see \eqref{6.63}--\eqref{6.80} for the details of its proof.

From the expression of $\tilde{q}$ in \eqref{6.91},
in order to control the second term  $r^{N-1}\tilde{q}$ on the right-hand side of \eqref{6.82},
we need to obtain some decay rate estimate of $(\rho^{\v,\d,b}-\bar{\rho}, m^{\v,\d,b})(t,r)$
as $r\rightarrow\infty$.
To achieve this, we first obtain the upper and lower bounds of density $\rho^{\v,\d,b}$ so that they are independent of $b$.
With these bounds of the density and property \eqref{4.1a} satisfied by the approximate initial data,
we can prove a better decay estimate for $(\rho^{\v,\d,b}-\bar{\rho}, m^{\v,\d,b})(t,r)$, uniformly in $b$;
see Lemmas \ref{lem6.10}--\ref{lem6.11} in more detail.
Then the decay estimate allows us to control $r^{N-1}\tilde{q}$.
Since the boundary values of $(\rho^{\v,\d,b}, u_r^{\v,\d,b})(t,b)$ are determined by the equations
and may depend on $\v$,
we integrate \eqref{6.82} over $[0,T]\times[b-1,b]\times[d,D]$ to avoid the trace estimates,
so that Proposition \ref{lem6.4} is obtained.
Then we take the limit, $b\rightarrow\infty$,
to obtain the global existence of a strong solution
$(\rho^{\v,\d},\mathcal{M}^{\v,\d})=(\rho^{\v,\d},m^{\v,\d}\,\frac{\mathbf{x}}{r})$ for \eqref{1.1}
on  $[0,\infty)\times (\mathbb{R}^N\setminus B_\delta(\mathbf{0}))$ for each fixed $\delta>0$.
Noting that the second term on the right-hand side of \eqref{6.15-1} vanishes
when $b\rightarrow\infty$,
we obtain the desired estimates in Proposition \ref{prop5.1}.

By similar arguments as in \cite{MV,Guo-Jiu-Xin-2},
we can then take the limit, $\d\rightarrow 0+$, to obtain the global weak solution
$(\rho^\v, \mathcal{M}^\v)=(\rho^\v, m^\v\,\frac{\mathbf{x}}{r})$ of the Cauchy problem
for \eqref{1.1}.
To prove that
\begin{equation}\nonumber
\partial_t\eta(\rho^\v,m^\v)+\partial_rq(\rho^\v,m^\v) \qquad  \mbox{is compact in }  \    H^{-1}_{\rm loc}(\mathbb{R}^2_+),
\end{equation}
special care is required,  since $(\rho^\v, m^\v)$ is only a weak solution
and $\partial_t\eta(\rho^\v,m^\v)+\partial_rq(\rho^\v,m^\v) $
is only a local bounded Radon measure for each fixed $\v>0$.
Moreover, since the viscosity coefficients depend on the density,
we can not say that $(\frac{m^\v}{\rho^\v})_r$ is a function due to the possible
appearance of the vacuum in general so that it is not suitable to use the weak form of $(\rho^\v, m^\v)$
to prove the $H^{-1}_{\rm loc}$--compactness.
In fact, the $H^{-1}_{\rm loc}$--compactness is achieved through smooth approximate solutions and their limits.

Based on the uniform estimates and the $H^{-1}_{\rm loc}$--compactness,
we then employ the compactness framework in \cite{Chen6} to take the vanishing viscosity $\v\to 0$ for all $\gamma>1$.
On the other hand, we have to be careful to pass the limit, $\v\rightarrow0$, in the momentum equations (see \eqref{5.48}),
since it is quite delicate to vanish the right-hand side of \eqref{5.48} by using the uniform estimates
in Theorem \ref{thm5.10}.
To overcome this difficulty, we employ underlying cancellations and introduce a new function $V^\v $,
which is uniformly bounded in $L^{2}(0,T; L^2)$  so that the right-hand side of \eqref{5.48}
is expressed by \eqref{5.48-5}. Then we can vanish the viscosity terms by using the new expression.

\vspace{1.5mm}
The paper is organized as follows:
In \S 2, we first introduce the notion of relative finite-energy solutions of the Cauchy problem \eqref{1.1-1}--\eqref{1.1-2}
for the compressible
Euler equations and then state Main Theorem I: Theorem \ref{thm:2.1} for the global existence of such solutions.
To establish Theorem 2.2, we construct global weak solutions of the Cauchy problem \eqref{1.1} and \eqref{initial-data}
for the compressible Navier-Stokes equations and analyze their vanishing viscosity limit, as stated in Main Theorem II: Theorem \ref{thm1.1}.
We also give several related remarks.
In \S 3, we first construct global approximate smooth solutions
$(\rho^{\v,\delta,b}, m^{\v, \delta,b})$ and
make the basic energy estimate and the BD entropy estimate of $(\rho^{\v,\delta,b}, m^{\v, \delta,b})$,
uniformly bounded in $(\v, \delta, b)$, for the Navier-Stokes equations \eqref{6.1}.
In \S 4, we derive the higher integrability of the approximate smooth solutions
$(\rho^{\v,\delta,b}, m^{\v, \delta,b})$ uniformly in $b$.
In \S 5, we first take the limit, $b\to \infty$, of $(\rho^{\v,\delta,b}, m^{\v, \delta,b})$
to obtain global strong solutions $(\rho^{\v,\delta}, m^{\v, \delta})$ of system \eqref{6.1}
with some uniform bounds in $(\v, \delta)$, and then
we take the limit, $\delta\to 0+$, to obtain global, spherically symmetric weak solutions of
the Navier-Stokes equations \eqref{1.1} with some desired uniform bounds
and the $H_{\rm loc}^{-1}$--compactness,
which are essential for us to employ the compensated compactness framework in \S 6
to establish
Theorem 2.2.
In the appendix, we construct the approximate initial data with desired estimates, which
are used for the construction of the approximate solutions in \S 3.

Throughout this paper,
we denote $L^p(\Omega), W^{k,p}(\Omega)$, and $H^k(\Omega)$ as the standard Sobolev
spaces on domain $\Omega$ for $p\in [1,\infty]$.
We also use $L^p(\Omega;\, r^{N-1}\dd r)$ or $L^p([0, T)\times \Omega;\,r^{N-1}\dd r \dd t)$
for $\Omega\subset \R_+$ with measure $r^{N-1}\dd r\,$ or $r^{N-1}\dd r \dd t\,$ correspondingly, and
$L^p_{\rm loc}([0,\infty); r^{N-1}\dd r)$ to represent $L^p([0,R);\, r^{N-1}\dd r)$ for any fixed $R>0$.

\section{Mathematical Problems and Main Theorems}

In this section, we first introduce the notion of relative finite-energy solutions
of the Cauchy problem \eqref{1.1-1}--\eqref{1.1-2} for the compressible Euler equations.

\begin{definition}\label{definition-Euler}
A pair $(\rho, \M)$ is said to be a relative finite-energy solution of the Cauchy problem \eqref{1.1-1}--\eqref{1.1-2}
if the following conditions hold{\rm :}

\begin{enumerate}
\item[\rm (i)]	
$\rho(t,\mathbf{x})\geqq 0$ a.e., and $(\mathcal{M}, \frac{\mathcal{M}}{\sqrt{\rho}})(t,\mathbf{x})=\mathbf{0}$ {\it a.e.} on the vacuum states
$\{(t,\mathbf{x})\,:\, \rho(t,\mathbf{x})=0\}${\rm ;}

\item[\rm (ii)] For a.e. $t>0$, the total relative energy with respect to the far-field state $(\bar{\rho},\mathbf{0})$
is finite{\rm :}
\begin{align} \label{1.20}
\int_{\mathbb{R}^N}\Big(\frac{1}{2}\big|\frac{\M}{\sqrt{\rho}}\big|^2
     +e (\rho, \bar{\rho})\Big)(t,\mathbf{x})\, \dd \mathbf{x}
\leqq E_0,
\end{align}
where
\begin{equation}\label{1.20a}
E_0:=\int_{\mathbb{R}^N}\Big(\frac{1}{2}\big|\frac{\M_0}{\sqrt{\rho_0}}\big|^2
  +e (\rho_0, \bar{\rho})\Big)(\mathbf{x})\,\dd\mathbf{x}<\infty
\end{equation}
is the finite total relative initial-energy, and $e(\rho, \bar{\rho})$ is the relative internal energy respective to $\bar{\rho}>0$:
\begin{equation}\label{6.10}
e(\rho,\bar{\rho}):=\frac{\kappa}{\gamma-1} \big(\rho^\g-\bar{\r}^{\g}-\g \bar{\rho}^{\g-1} (\rho-\bar{\r})\big){\rm ;}
\end{equation}

\smallskip
\item[\rm (iii)] For any $\zeta(t,\mathbf{x})\in C^{1}_0([0,\infty)\times\mathbb{R}^N)$,
\begin{align}\label{1.18}
\int_{\mathbb{R}^{N+1}_+} \big(\rho \zeta_t + \M\cdot\nabla\zeta\big)\,\dd\mathbf{x}\dd t
+\int_{\mathbb{R}^N} (\rho_0 \zeta)(0,\mathbf{x})\,\dd\mathbf{x}=0{\rm ;}
\end{align}

\item[\rm (iv)] For all $\psi(t,\mathbf{x})=(\psi_1,\cdots,\psi_N)(t,\mathbf{x})\in \big(C^1_0([0,\infty)\times\mathbb{R}^N)\big)^N$,
\begin{align}\label{1.19}
&\int_{\mathbb{R}^{N+1}_+} \Big(\M\cdot\partial_t\psi +\frac{\M}{\sqrt{\rho}}\cdot \big(\frac{\M}{\sqrt{\rho}}\cdot \nabla\big)\psi
+p(\rho)\,\mbox{\rm div}\,\psi\Big)\, \dd\mathbf{x}\dd t\nonumber\\
&\quad+\int_{\mathbb{R}^N} \M_0(\mathbf{x})\cdot \psi(0,\mathbf{x})\,\dd\mathbf{x}
=0,
\end{align}
\end{enumerate}
where and whereafter we always use $\mathbb{R}_+^{N+1}:=\mathbb{R}_+\times \mathbb{R}^N=(0, \infty)\times \mathbb{R}^N$ for $N\ge 2$.
\end{definition}

Our first main theorem of this paper is

\begin{theorem}[Main Theorem I:  Existence of Spherically Symmetric Solutions of the Euler Equations]\label{thm:2.1}
Consider the Cauchy problem of the Euler equations \eqref{1.1-1} with
large initial data of spherical symmetry of form \eqref{initial}.
Let $(\rho_0, \M_0)(\mathbf{x})$ satisfy \eqref{1.20a} with the positive far-field density $\bar{\rho}>0$.
Then there exists a global
relative finite-energy solution $(\rho, \M)(t,\mathbf{x})$ of \eqref{1.1-1} and \eqref{initial} with spherical symmetry
of form \eqref{1.2} in the sense of  Definition {\rm \ref{definition-Euler}},
where $(\rho, m)(t,r)$ is determined by the corresponding
Cauchy problem of system \eqref{euler} with the initial data
$(\rho_0, m_0)(r)$ given in \eqref{initial}.
\end{theorem}

To establish Theorem \ref{thm:2.1},
we first construct global weak solutions of the Cauchy problem of the compressible Navier-Stokes equations \eqref{1.1}
with appropriately adapted degenerate density-dependent viscosity terms and approximate initial data:
\begin{equation}\label{initial-data}
(\rho, \M)|_{t=0}=(\rho_0^\v, \M_0^\v)(\mathbf{x})\longrightarrow (\rho_0, \M_0)(\mathbf{x}) \quad\,\,\,\mbox{as $\v\to 0$},
\end{equation}
constructed as in the appendix satisfying Lemmas \ref{lem8.2}--\ref{lem8.3} and Lemma \ref{lem8.4}(i).

For clarity, we adapt the viscosity terms with
$(\mu,\lambda)=(\rho,0)$
in \eqref{1.1}, as the case for the shallow water (Saint Venant) models,
and $\v\in (0,1]$
without loss of generality throughout this paper.
The arguments also work for a general class of degenerate density-dependent viscosity terms;
see Remark \ref{rem:2.7} below for more details.

\begin{definition}\label{definition-NS}
A pair $(\rho^\v, \M^\v)$ is said to be a weak solution of the Cauchy problem \eqref{1.1} and \eqref{initial-data}
with $(\mu, \lambda)=(\rho, 0)$
if the following conditions hold{\rm :}

\begin{enumerate}
\item[\rm (i)]	
$\rho^\v(t,\mathbf{x})\geqq 0$ a.e., and $(\mathcal{M}^\v,\frac{\mathcal{M}^\v}{\sqrt{\rho^\v}})(t,\mathbf{x})=\mathbf{0}$
a.e. on the vacuum states $\{(t,\mathbf{x})\,:\,\rho^\v(t,\mathbf{x})=0\}$,
\begin{align*}
&\rho^\v\in L^{\infty}(0,T; L^\gamma_{\rm loc}(\mathbb{R}^N)),
\quad   \nabla\sqrt{\rho^\v}\in  \big(L^{\infty}(0,T; L^2(\mathbb{R}^N))\big)^N, \\
& \frac{\M^\v}{\sqrt{\rho^\v}}\in \big(L^{\infty}(0,T; L^2(\mathbb{R}^N))\big)^N;
\end{align*}
	
\item[\rm (ii)] For any $t_2\geqq t_1\geqq 0$ and any $\zeta(t,\mathbf{x})\in C^1_0([0,\infty)\times\mathbb{R}^{N})$,
the mass equation $\eqref{1.1}_1$ holds in the sense{\rm :}
\begin{align*}
&\int_{\mathbb{R}^N}(\rho^\v \zeta )(t_2,\mathbf{x})\, \dd\mathbf{x}
-\int_{\mathbb{R}^N}(\rho^\v \zeta )(t_1,\mathbf{x})\, \dd\mathbf{x}\\
&=\int^{t_2}_{t_1}\int_{\mathbb{R}^N}\big(\rho^\v\zeta_t+\M^\v\cdot \nabla \zeta\big)(t,\mathbf{x})\, \dd\mathbf{x} \dd t;
\end{align*}
	
\item[\rm (iii)] For any $\psi=(\psi_1,\cdots,\psi_N)\in \big(C^2_0([0,\infty)\times \mathbb{R}^N)\big)^N$,
the momentum equations $\eqref{1.1}_2$ hold in the sense{\rm :}
\begin{align*}
&\int_{\mathbb{R}^{N+1}_+} \Big(\M^\v\cdot \psi_t
+\frac{\M^\v}{\sqrt{\rho^\v}}\cdot \big(\frac{\M^\v}{\sqrt{\rho^\v}}\cdot \nabla\big)\psi
+ p(\rho^\v)\, \mbox{\rm div}\,\psi\Big)\,\dd\mathbf{x}\dd t\nonumber\\
&\quad +\int_{\mathbb{R}^N} \M_0^\v(\mathbf{x})
\cdot\psi(0,\mathbf{x})\, \dd\mathbf{x}\nonumber\\
&=-\v\int_{\mathbb{R}^{N+1}_+}
\Big(\frac{1}{2}\M^\v\cdot \big(\Delta \psi
  +\nabla\mbox{\rm div}\,\psi\big)
+\frac{\M^\v}{\sqrt{\rho^\v}}\cdot \big(\nabla\sqrt{\rho^\v}\cdot \nabla\big)\psi
\nonumber\\
&\qquad\qquad\qquad+\nabla\sqrt{\rho^\v}\cdot \big(\frac{\M^\v}{\sqrt{\rho^\v}}\cdot \nabla\big)\psi\Big)\,  \dd\mathbf{x}\dd t.\nonumber
\end{align*}
\end{enumerate}
\end{definition}

\vspace{1.5mm}
Consider spherically symmetric solutions of form \eqref{1.2}.
Then systems \eqref{1.1-1} and \eqref{1.1} for such solutions
become \eqref{euler} and \eqref{1.3}, respectively.
A pair of functions
$(\eta(\rho,m),q(\rho,m))$
is called an entropy pair of the $1$-D Euler system ({\it i.e.}, system \eqref{euler} with $N=1$)
if they satisfy
$$
\partial_t\eta(\rho,m)+ \partial_r q(\rho,m)=0
$$
for any smooth solution $(\rho, m)$ of the 1-D Euler system;
see Lax \cite{P. D. Lax}.
Furthermore, $\eta(\rho,m)$ is called a weak entropy if
\begin{align}\nonumber
\eta|_{\rho=0}=0\qquad \text{for any fixed $u=\frac{m}{\rho}$}.
\end{align}

From now on, we also use $u=\frac{m}{\rho}$ and $m$ alternatively when $\rho>0$.

From \cite{Lions P.-L.2}, it is well-known that any weak entropy pair $(\eta,q)$ can be represented by
\begin{align}\label{weakentropy}
\begin{split}
\displaystyle\eta(\rho,m)&=\int_{\R}\chi(\rho;s-u)\psi(s)\,\dd s,\\
\displaystyle q(\rho,m)&=\int_{\R}(\theta s+(1-\theta)u)\chi(\rho;s-u)\psi(s)\,\dd s
\end{split}
\end{align}
when $\rho>0$, where the kernel is
$$
\chi(\rho;s-u)=[\rho^{2\theta}-(s-u)^2]_{+}^{\hl}
\quad\,\, \mbox{for $\hl:=\frac{3-\gamma}{2(\gamma-1)}>-\frac{1}{2}$ and $\theta:=\frac{\gamma-1}{2}$}.
$$
For instance, when $\psi(s)=\frac{1}{2}s^2$, the entropy pair consists of the mechanical energy and the associated energy flux:
\begin{align}\label{mechanical}
\eta^{*}(\rho,m)=\frac{1}{2}\frac{m^2}{\rho}+e(\rho),
\quad q^{*}(\rho,m)=\frac{1}{2}\frac{m^3}{\rho^2}+m e'(\rho),
\end{align}
where $e(\rho)=\frac{\kappa}{\gamma-1}\rho^{\gamma}$ represents the internal energy.
Since we expect that $(\rho, m)(t,r)\to (\bar{\rho}, 0)$ with $\bar{\rho}>0$ as $r\to \infty$,
we define the relative mechanical energy
\begin{align}\label{relative-entropy}
\bar{\eta}^{\ast}(\rho,m)=\frac{m^2}{2\rho}+e(\rho,\bar\rho)
\end{align}
with $e(\rho,\bar{\rho})$ defined by \eqref{6.10}
satisfying (see \cite{Chen6}):
\begin{align}\label{6.10-1}
e(\rho,\bar{\r})\geqq C_{\g} \rho (\rho^{\t}-\bar{\r}^{\t})^2
\end{align}
for some constant $C_\gamma>0$.

\begin{theorem}[Main Theorem II: Existence and Inviscid Limit for the Navier-Stokes Equations]\label{thm1.1}
Consider the compressible Navier-Stokes equations \eqref{1.1} with $N\ge 2$
and the spherically symmetric approximate initial data  \eqref{initial-data} satisfying that, as $\v\to 0$,
\begin{align}
&(\rho_0^\v,m_0^\v)(r)\to (\rho_0, m_0)(r) \qquad \mbox{in $L^1_{\rm loc}([0,\infty); r^{N-1}{\rm d}r)$}, \label{1.53o}\\
&E_0^\v:=\omega_N\int_0^\infty\bar{\eta}^{\ast}(\rho_0^\v,m_0^\v)\, r^{N-1}\dd r\to E_0,  \label{1.53}\\
&E_1^\v:=\v^2\int_0^\infty \big|\big(\sqrt{\rho_{0}^\v}\big)_r\big|^2\, r^{N-1}\dd r \to 0,  \label{1.53a}
\end{align}
and there exists a constant $C>0$ independent of $\v\in (0,1]$ such that
\begin{equation}\label{1.53b}
E_0^\v+E_1^\v\leqq C(E_0+1)
\end{equation}
for $E_0$ defined in \eqref{1.20a}
and $\omega_N=2\pi^{\frac{N}{2}}\Gamma(\frac{N}{2})^{-1}$
as the surface area of the unit ball in $\mathbb{R}^N$.
Then the following statements hold{\rm :}

\vspace{2mm}
\noindent{\bf Part I. Existence for the Navier-Stokes Equations \eqref{1.1}:}
For each $\v>0$, there exists a global spherically symmetric weak
solution
$$
(\rho^\v, \M^\v)(t,\mathbf{x})=(\rho^\v(t,r),m^\v(t,r) \frac{\mathbf{x}}{r})=(\rho^\v(t,r),\rho^\v(t,r)u^\v(t,r) \frac{\mathbf{x}}{r})
$$
of the Cauchy problem of \eqref{1.1} and \eqref{initial-data}
in the sense of Definition {\rm \ref{definition-NS}}, where
$u^\v(t,r)=\frac{m^\v(t,r)}{\rho^\v(t,r)}$ a.e. on $\{(t,r)\,:\, \rho^\v(t,r)\ne 0\}$ and $u^\v(t,r)=0$ a.e. on $\{(t,r)\,:\, \rho^\v(t,r)=0\}$.
Moreover, $(\rho^\v,m^\v)(t,r)$
satisfies the following uniform bounds{\rm :} For $t>0$,
\begin{align}
&\int_0^\infty  \bar{\eta}^{\ast}(\rho^\v,m^\v)(t,r)\,r^{N-1}\dd r
 +\v\int_{\mathbb{R}_+^2} \rho^\v(s,r)|u^\v(s,r)|^2\,r^{N-3}\dd r \dd s\nonumber\\
&\leqq \frac{E_0^\v}{\omega_N}\leqq C(E_0+1),\label{1.55}\\[1.5mm]
&\v^2 \int_0^\infty \big|\big(\sqrt{\rho^{\v}(t,r)}\big)_r\big|^2 \, r^{N-1}\dd r
  +\v\int_{\mathbb{R}_+^2} \big|\big((\rho^{\v}(s,r))^{\frac{\gamma}{2}}\big)_r\big|^2\,r^{N-1}\dd r \dd s \nonumber\\
&\leqq C (E_0+1),\label{1.56}
\end{align}
and
\begin{align}
&\int_0^T\int_{d}^D \big(\rho^{\v}(t,r)\big)^{\g+1}\,\dd r\dd t\leqq C(d,D,T, E_0),\label{1.57}\\[1mm]
&\int_0^T\int_0^D
\big(\rho^{\v}(t,r)|u^\v(t,r)|^3+\big(\rho^{\v}(t,r)\big)^{\gamma+\theta}\big)\,r^{N-1}\dd r\dd t
  \leqq C(D,T,E_0) \label{1.58}
\end{align}
for any fixed $T\in (0,\infty)$ and any compact subset $[d,D]\Subset (0,\infty)$,
where and whereafter we denote $\mathbb{R}^2_+:=\{(t,r)\ :\ t\in (0,\infty),\ r\in(0,\infty) \}$,
and $C>0$
and $C(d, D,T,E_0)>0$
as two universal constants independent of $\v$, but depending on $(\gamma,N)$
and $(d, D,T,E_0)$, respectively.

Let $(\eta, q)$ be an entropy pair defined in \eqref{weakentropy}
for a smooth compact supported function $\psi(s)$ on $\mathbb{R}$.
Then, for $\v\in (0,1]$,
\begin{align}\label{1.58-1}
\partial_t\eta(\rho^\v,m^\v)+\partial_rq(\rho^\v,m^\v) \qquad \mbox{is compact in }\, H^{-1}_{\rm loc}(\mathbb{R}^2_+),
\end{align}
where $H^{-1}_{\rm loc}(\mathbb{R}^2_+)$ represents $H^{-1}((0,T]\times\Omega)$ for any $T>0$ and bounded open subset $\Omega\Subset (0,\infty)$.

\vspace{2.5mm}
\noindent{\bf Part II. Inviscid Limit to the Euler Equations \eqref{1.1-1}:}
For the global weak solutions $(\rho^{\v}, \mathcal{M}^\v)$ of the compressible Navier-Stokes equations \eqref{1.1}
established in Part {\rm I},
there exist a subsequence {\rm (}still denoted{\rm )} $(\rho^{\v},m^{\v})$ and a vector function $(\rho, m)$
such that, as $\v\rightarrow0$,
\begin{align*}
&(\rho^{\v},m^{\v})\rightarrow (\rho,m)(t,r)\quad\,\, \mbox{in $(L^p_{\rm loc}\times L^q_{\rm loc})([0,\infty); r^{N-1}\dd r)$},\\
&\int_0^T\int_0^D \big|\big(\frac{m^\v}{\sqrt{\rho^\v}}\big)(t,r)-\big(\frac{m}{\sqrt{\rho}}\big)(t,r)\big|^2\,r^{N-1}\dd r
\dd t\to 0 \quad\, \mbox{for any fixed $T,D\in (0, \infty)$},
\end{align*}
where  $p\in[1,\gamma+1),\,\, q\in[1,\frac{3(\gamma+1)}{\gamma+3})$, and $(\rho, \mathcal{M})(t,\mathbf{x}):=(\rho(t,r), m(t,r)\frac{\mathbf{x}}{r})$ is a global relative finite-energy
solution of spherical
symmetry of the Euler equations \eqref{1.1-1} with initial data \eqref{initial} in the sense of
Definition {\rm \ref{definition-Euler}}.
\end{theorem}

\smallskip
\begin{remark}\label{remark:2.5a}
In Theorem {\rm \ref{thm1.1}}, the approximate initial data functions $(\rho_0^\v, m_0^\v)$
satisfying conditions {\rm \eqref{1.53o}}--{\rm \eqref{1.53a}}
are constructed in Lemmas {\rm \ref{lem8.2}--\ref{lem8.3}} and Lemma {\rm \ref{lem8.4}(i)}
in the appendix.
Then Theorem {\rm \ref{thm:2.1}} is a direct corollary of Theorem {\rm \ref{thm1.1}}.
\end{remark}

\begin{remark}
The main point of Theorem {\rm \ref{thm1.1}} is to construct suitable Navier-Stokes approximate solutions that converge strongly to
a global relative finite-energy solution of spherical symmetry
of the Euler equations \eqref{1.1-1} with initial data \eqref{initial} in the sense of Definition {\rm \ref{definition-Euler}}
under the relative finite-energy condition \eqref{1.20a} only.
We can follow the same arguments as in {\rm \S 3--\S 6} to obtain
a rigorous proof of the inviscid limit from the Navier-Stokes to Euler equations with fixed same initial data $(\rho_0, m_0)$
of appropriate regularity and decay at infinity.
\end{remark}

\begin{remark}\label{rem:2.7}
When both $\mu$ and $\lambda$ are constants, it is still an open problem for
the inviscid limit from \eqref{1.3} to \eqref{euler}, since the BD entropy estimate
is invalid for this case so that the required uniform estimate for the derivative of the density has not obtained yet.
On the other hand, our analysis in this paper applies to a class of more general viscosity coefficients $(\mu(\rho),\lambda(\rho))$.
For instance, our results hold for the class of $(\mu(\rho), \lambda(\rho))$ that satisfy the BD relation $($see \cite{BD-2004-CRMASP,MV}$)${\rm :}
\begin{equation}\label{6.1-3}
\lambda(\rho)=\rho \mu'(\rho)-\mu(\rho)
\end{equation}
with some additional conditions{\rm ;}
see also the approximate system \eqref{6.1}--\eqref{6.3}.
\end{remark}

\section{Approximate Solutions and Basic Uniform Estimates}

In this section, we first construct global approximate smooth solutions and make their basic energy estimate
and the BD entropy estimate, uniformly bounded with respect to the approximation parameters.

The main difficulty is to obtain some uniform estimates directly for the exact solutions of
the Navier-Stokes equations \eqref{1.1} with
approximate initial data \eqref{initial},
owing to the potential appearance of the vacuum and singularity of their limits at both the origin, $r=0$, and the far-field, $r=\infty$,
generically.
On the other hand, for our purpose, it suffices to obtain first uniform estimates for appropriately designed
approximate solutions of the Navier-Stokes equations \eqref{1.1}.
To achieve these, we construct the approximate solutions as the solutions of
the following approximate Navier-Stokes system with positive density ({\it i.e.}, $\rho>0$ so that
the velocity, $u=\frac{m}{\rho}$, is well-defined) in truncated domains:
\begin{align}\label{6.1}
\begin{cases}
\r_t+(\r u)_r+\f{N-1}r \r u=0,\\[1mm]
(\r u)_t+(\r u^2+p)_r+\f{N-1}r \r u^2=\v\big((\mu+\lambda)(u_r+\frac{N-1}{r}u)\big)_r-\v\frac{N-1}{r}u\mu_{r},
\end{cases}
\end{align}
where $t>0$ and $r\in [\d, b]$ with $\delta \in (0, 1]$ and $b\ge 1+\delta^{-1}$, and
\begin{equation}\label{6.1-2}
\mu(\rho)=\rho+\delta \rho^{\alpha},\qquad \lambda(\rho)=\delta (\alpha-1) \rho^{\alpha}
\end{equation}
with $\alpha\in(\frac{N-1}{N},1)$.
For concreteness, we take $\alpha=\frac{2N-1}{2N}$.
It is easy to check that $(\mu(\rho), \lambda(\rho))$ in \eqref{6.1-2}
satisfy relation \eqref{6.1-3}.

We impose \eqref{6.1} with the following approximate initial data:
\begin{align}\label{6.2}
(\rho, u)(0,r)=(\r_0^{\v,\d,b},u_0^{\v,\d,b})(r) \quad\,\, \mbox{for $r\in[\d,b]$},
\end{align}
and the boundary condition:
\begin{align}\label{6.3}
u(t,\d)=u(t,b)=0\quad\,\, \mbox{for $t>0$},
\end{align}
where $\r_0^{\v,\d,b}$ and $u_0^{\v,\d,b}$ are smooth functions
satisfying
\begin{equation}\label{6.3a}
0<(\beta \v)^{\frac{1}{4}}\leqq \r_0^{\v,\d,b}\leqq (\beta \v)^{-\frac{1}{2}}<\infty,
\end{equation}
for some small constant $\beta$ (determined in Lemma \ref{lem8.2}).

Such approximate initial data functions in \eqref{6.2}
have been constructed in the appendix, which satisfy all the properties
in Lemmas \ref{lem8.2}--\ref{lem8.4}.

For $N=2,3$, the existence of global smooth solutions $(\r^{\v,\d,b}, u^{\v,\d,b})$ of
\eqref{6.1}--\eqref{6.3} with $0<\r^{\v,\d,b}(t,r)<\infty$ can be established as
in Guo-Jiu-Xin \cite{Guo-Jiu-Xin-2}.
In fact, for any $N\geqq 2$, a similar global existence result for smooth solutions of the approximate system \eqref{6.1}--\eqref{6.3}
can be obtained by using analogous arguments as in \S 3 and \S 4.1 of \cite{Guo-Jiu-Xin-2}; see also \cite{D. Hoff1,Song Jiang}.
Since the upper and lower bounds of $\r^{\v,\d,b}$ in \cite{Guo-Jiu-Xin-2} depend on parameters $(\v,\d, b)$,
the key point of this section is to obtain some uniform estimates of $(\r^{\v,\d,b}, u^{\v,\d,b})$ independent of $(\d, b)$ so that
both limits $b\rightarrow\infty$ and $\d\rightarrow0+$ can be taken
to obtain the global weak solution of \eqref{1.1} and \eqref{initial-data}; see \S 5.

Throughout this section, for simplicity,  we always fix parameters $\v,\delta\in(0,1]$ and $b\ge 1+\delta^{-1}$,
use $u^{\v,\d,b}$ or $m^{\v,\d,b}$ alternatively since $\rho^{\v,\d,b}$ is positive,
and drop the superscripts of solution $(\rho^{\v,\d, b}, u^{\v,\d,b})(t,r)$
and the approximate initial data $(\rho^{\v,\d, b}_0, u^{\v,\d,b}_0)$, when no confusion arises.
We keep the superscripts when the initial data functions are involved.

\begin{lemma}[Basic Energy Estimate]\label{lem6.1}
The smooth solution $(\rho, u)$ of \eqref{6.1}--\eqref{6.3} satisfies that, for any $t>0$,
\begin{align}\label{6.11-1}
&\int_\d^b  \big(\frac12\rho u^2+e(\rho,\bar{\r}) \big)(t,r)\,r^{N-1}\dd r
  +\v\int_0^t\int_\d^b \big(\rho u_r^2+ \frac{N-1}{r^2}\rho u^2\big)(s,r)\,r^{N-1}\dd r \dd s\nonumber\\
&+\v\delta \int_0^t\int_\d^b \rho^{\alpha}\Big\{\alpha  u_r^2+2(\alpha-1)(N-1)\frac{u u_r}{r}\nonumber\\
&\qquad\qquad\qquad\quad  + \big(1+(N-1)(\alpha-1)\big)(N-1) \frac{u^2}{r^2}\Big\}(s,r)\, r^{N-1}\dd r \dd s\nonumber\\
&=\int_\d^b  \big(\frac12\rho_0 u_0^2+e(\rho_0,\bar{\r}) \big)(r)\,r^{N-1}\dd r=:\frac{E_0^{\v,\d,b}}{\omega_N},
\end{align}
where $E_0^{\v,\d,b}$ satisfies the properties stated in Lemma {\rm \ref{lem8.4}} in the appendix.
In particular,
there exists a positive constant $c_N>0$ {\rm (}depending only on $N${\rm )} such that
\begin{align}\label{6.11}
&\int_\d^b \big(\frac12\rho u^2+e(\rho,\bar{\r}) \big)(t,r)\, r^{N-1}\dd r
 +\v\int_0^t\int_\d^b \Big(\rho u_r^2+ \frac{\rho u^2}{r^2}\Big)(s,r)\,r^{N-1}\dd r \dd s\nonumber\\
&\quad+c_N\v\delta \int_0^t\int_\d^b \Big(\rho^{\alpha} u_r^2+\frac{\rho^{\alpha} u^2}{r^2}\Big)(s,r)\,r^{N-1}\dd r \dd s\nonumber\\
&\leqq \frac{E_0^{\v,\d,b}}{\omega_N}
\leqq C(E_0+1) \quad\mbox{for any $t>0$},
\end{align}
for some constant $C>0$ independent of $(\v,\d,b)$, where we have used \eqref{A.39f}.
\end{lemma}

\noindent{\bf Proof.}
Multiplying $\eqref{6.1}_2$ by $r^{N-1}u$ and performing integration by parts, we have
\begin{align}\label{6.11-2}
&\frac{\dd}{\dd t} \int_\d^b \frac12 \rho u^2\, r^{N-1}\dd r
+\int_\delta^b  p_r u\, r^{N-1}\dd r\nonumber\\
&=-\v\int_\d^b \Big((\mu+\lambda) \big(u_r+\frac{N-1}{r} u\big) (r^{N-1} u)_r-(N-1) \mu (r^{N-2} u^2)_r\Big) \dd r.
\end{align}

For the second term on the left-hand side of \eqref{6.11-2},
it follows from $\eqref{6.1}_1$ and integration by parts that
\begin{align}\label{6.11-3}
\int_\d^b  p_r u \, r^{N-1}\dd r
&=\frac{\kappa \gamma}{\gamma-1} \int_\d^b \rho u  (\rho^{\gamma-1})_r\, r^{N-1}\dd r\nonumber\\
&=-\frac{\kappa \gamma}{\gamma-1} \int_\d^b (\rho u r^{N-1})_r \rho^{\gamma-1}\, \dd r
\nonumber\\
&=\frac{\kappa}{\gamma-1} \int_\d^b (\rho^{\gamma} )_t\, r^{N-1} \dd r\nonumber\\
&=\frac{\kappa}{\gamma-1} \int_\d^b \big(\rho^{\gamma}
 -\bar\rho^{\gamma}-\gamma \bar{\rho}^{\gamma-1}(\rho-\bar{\rho})\big)_t\, r^{N-1}\dd r\nonumber\\
&=\frac{\dd}{\dd t}\int_\d^b e(\rho,\bar{\rho})(t,r)\, r^{N-1}\dd r.
\end{align}
For the viscous term, a direct calculation shows
\begin{align}\label{6.11-4}
&(\mu+\lambda) \big(u_r+\frac{N-1}{r} u\big) (u r^{N-1})_r-(N-1) \mu (u^2 r^{N-2})_r\nonumber\\[1mm]
&=\mu \big(r^{N-1} u_r^2+(N-1) r^{N-3} u^2\big)\nonumber\\
&\qquad\,\,+\lambda \big(r^{N-1} u_r^2+2(N-1) r^{N-2} u u_r+(N-1)^2 r^{N-3} u^2\big) \nonumber\\[1mm]
&=\delta \rho^{\alpha} \Big(\alpha r^{N-1} u_r^2+2(\alpha-1)(N-1) r^{N-2} u u_r \nonumber\\
&\qquad\quad\,\,+(N-1)(1+(\alpha-1)(N-1)) r^{N-3} u^2\Big)\nonumber\\[1mm]
&\quad\, +\rho \big(r^{N-1} u_r^2+(N-1) r^{N-3} u^2\big).
\end{align}
For the first term on the right-hand side of \eqref{6.11-4},
we calculate its discriminant:
\begin{align*}
&4(\alpha-1)(N-1)^2 -4\alpha (N-1)\big(1+(\alpha-1)(N-1)\big)\nonumber\\
&=4(N-1)^2 \big(1-\frac{N}{N-1} \alpha\big)<0,
\end{align*}
since $\alpha\in (\frac{N-1}{N}, 1)$.
Thus, there exists a positive constant $c_N>0$ such that
\begin{align}\label{6.11-5}
&(\mu+\lambda) \big(u_r+\frac{N-1}{r} u\big) (r^{N-1} u)_r-(N-1) \mu  (r^{N-2} u^2)_r
\nonumber\\
&\geqq \rho \big(u_r^2+ \frac{u^2}{r^2}\big)r^{N-1}+ c_N\delta \rho^{\alpha}\big(u_r^2+ \frac{u^2}{r^2}\big)r^{N-1}.
\end{align}
Integrating \eqref{6.11-2} over $[0,t]$ and using \eqref{6.11-3}--\eqref{6.11-5},
we obtain \eqref{6.11-1}--\eqref{6.11}.
$\hfill\Box$

For $(\mu,\lambda)$ determined by \eqref{6.1-2}, system \eqref{1.1} admits an additional {\it a priori} estimate
for the density (via the BD entropy), as observed by Bresch-Desjardins \cite{BD-2003-CRMASP,BD-2004-CRMASP}
(see also Bresch--Desjardins--Gerard-Varet \cite{BDG-2007})
with the Dirichlet boundary conditions in the $3$-D case.
For  the spherically symmetric problem, we have

\begin{lemma}[BD Entropy Estimate]\label{lem6.2}
The smooth solution of \eqref{6.1}--\eqref{6.3} satisfies
\begin{align}\label{6.12}
&\v^2 \int_\d^b \big((1 +\d \rho^{\alpha-1}+\d^2 \rho^{2(\alpha-1)})\frac{\rho_r^2}{\rho}\big)(t,r)\,r^{N-1} \dd r\nonumber\\[1mm]
&+\v\int_0^t\int_\d^b \big((1 +\d \rho^{\alpha-1})\rho^{\g-2}\rho_r^2\big)(s,r)\,r^{N-1}\dd r \dd s
 \leqq C(E_0+1),
\end{align}
where we have used
\begin{align}
\sup_{0<\v, \delta\leqq 1}\sup_{b\ge 1+\d^{-1}}\big(E_0^{\v,\d,b}+E_1^{\v,\d,b}\big)\leqq C(E_0+1),
\end{align}
which follows from {\rm (A.38)},
with
\begin{align}\label{6.13-1}
E_1^{\v,\d,b}:=\v^2\int_\d^b \big(1+2\alpha\d \rho_0^{\alpha-1}+\alpha^2\d^2 \rho_0^{2\alpha-2}\big)\big|(\sqrt{\rho_{0}})_r\big|^2\,r^{N-1}\dd r.
\end{align}
\end{lemma}

\noindent{\bf Proof.}
It is more convenient to deal with  \eqref{6.1} in the Lagrangian coordinates for this proof.
We divide the proof into four steps.

\medskip
1. For simplicity, denote
$L_b:=\int_\d^b\rho_0(r)r^{N-1}\,\dd r$.
Note that
$$
\frac{\dd}{\dd t}\int_\d^{b}\rho(t,r)\,r^{N-1} \dd r
=-\int_\d^{b}(\rho u r^{N-1})_r(t,r)\,\dd r=0.
$$
Then
\begin{equation}
\int_\d^{b}\rho(t,r)\,r^{N-1}\dd r=\int_\d^b\rho_0(r)\,r^{N-1}\dd r
=L_b\qquad \mbox{for all $t>0$}.\nonumber
\end{equation}
For $r\in[\d,b]$ and $t\in[0,T]$,
we define the Lagrangian
transformation:
\begin{equation}
x=\int_\d^r \rho(t,y)\,y^{N-1}\dd y,\qquad \tau=t,\nonumber
\end{equation}
which translates  domain $[0,T]\times[\d,b]$ into $[0,T]\times[0,L_b]$ and satisfies
\begin{equation}\label{6.8}
\begin{cases}
\frac{\partial x}{\partial r}=\rho r^{N-1}>0, \,\,\,\frac{\partial x}{\partial t}=-\rho ur^{N-1},
   \,\,\,\frac{\partial\tau}{\partial r}=0, \,\,\,\frac{\partial\tau}{\partial t}=1,\\[2mm]
\frac{\partial r}{\partial x}=\frac{1}{\rho r^{N-1}}>0, \,\,\,\frac{\partial r}{\partial\tau}=u,
   \,\,\,\frac{\partial t}{\partial\tau}=1, \,\,\,\frac{\partial t}{\partial x}=0.
\end{cases}
\end{equation}
Applying the Lagrange transformation, system \eqref{6.1} becomes
\begin{align}\label{6.9}
\begin{cases}
\rho_\tau+\rho^2(r^{N-1}u)_x=0,\\[2mm]
u_\tau+r^{N-1}p_x=\v r^{N-1}\big(\rho(\mu+\lambda)(r^{N-1}u)_x\big)_x-\v (N-1) r^{N-2} \mu_x u,
\end{cases}
\end{align}
and the boundary condition \eqref{6.3}
becomes
\begin{equation}\label{6.9a}
u(\tau,0)=u(\tau,L_b)=0\qquad \mbox{for $\tau>0$}.
\end{equation}

\smallskip
2. Multiplying $\eqref{6.9}_1$ by $\mu'(\rho)$ and using \eqref{6.1-3}, we have
\begin{align}\label{6.9-1}
\mu_{\tau}+\rho(\mu+\lambda)  (r^{N-1}u)_x=0.
\end{align}
Substituting \eqref{6.9-1} into the viscous term of $\eqref{6.9}_2$ leads to
\begin{align}\label{6.9-2}
u_{\tau}+r^{N-1}p_x=-\v r^{N-1} (\mu_{x})_\tau-\v (N-1) r^{N-2} \mu_x u.
\end{align}
Note from \eqref{6.8} that  $\frac{\partial r}{\partial\tau}=u$.
Then the last term of \eqref{6.9-2} is rewritten as
\begin{equation}\nonumber
\v (N-1) r^{N-2}u \mu_x=(N-1) r^{N-2} r_\tau \mu_x=(r^{N-1})_\tau \mu_x,
\end{equation}
which,  with \eqref{6.9-2}, yields
\begin{align}\label{6.9-4}
(u+\v r^{N-1} \mu_x)_{\tau}+r^{N-1}p_x=0.
\end{align}

\smallskip
3. Multiplying \eqref{6.9-4} by $u+\v r^{N-1} \mu_x$, we have
\begin{align}\label{6.9-5}
\frac{1}{2}\frac{\dd}{\dd\tau}\int_0^{L_b}(u+\v r^{N-1} \mu_x)^2\,\dd x +\v \int_0^{L_b}p_x \mu_x\,r^{2N-2} \dd x
+\int_0^{L_b}p_x u\,r^{N-1} \dd x=0.
\end{align}
For the last term on the left-hand side of \eqref{6.9-5},
it follows from integration by parts and $\eqref{6.9}_1$ that
\begin{align}\label{6.9-6}
\int_0^{L_b}p_x u\, r^{N-1}\, \dd x&=-\int_0^{L_b} p (u r^{N-1})_x\, \dd x
=\kappa \int_0^{L_b} \rho^{\gamma-2} \rho_\tau\, \dd x \nonumber\\
&=\frac{\kappa}{\gamma-1} \int_0^{L_b} (\rho^{\gamma-1})_\tau\, \dd x
 =\frac{\dd}{\dd\tau}\int_0^{L_b} \frac{e(\rho,\bar{\rho})}{\rho}\, \dd x.
\end{align}
Substituting \eqref{6.9-6} into \eqref{6.9-5} leads to
\begin{align}\label{6.9-7}
\frac{\dd}{\dd\tau}\int_0^{L_b}\Big(\frac{1}{2}(u+\v r^{N-1} \mu_x)^2+\frac{e(\rho,\bar{\rho})}{\rho}
\Big)\,\dd x
 +\v \int_0^{L_b}p_x \mu_x\,r^{2N-2} \dd x=0.
\end{align}
Integrating \eqref{6.9-7} over $[0,\tau]$ yields
\begin{align}\label{6.9-8}
&\int_0^{L_b}\Big(\frac{1}{2} (u+\v r^{N-1} \mu_x)^2+\frac{e(\rho,\bar{\rho})}{\rho}
\Big)\,\dd x
  +\v \int_0^\tau \int_0^{L_b}p_x \mu_x\,r^{2N-2} \dd x \dd s\nonumber\\
&= \int_0^{L_b}\Big(\frac{1}{2} (u_0+\v r_0^{N-1} \mu_{0x})^2+\frac{e(\rho_0,\bar{\rho})}{\rho_0}\Big)\, \dd x.
\end{align}

\smallskip
4. Plugging \eqref{6.9-8} back to the Eulerian coordinates, we have
\begin{align*}
&\int_a^{b}\Big(\frac{1}{2} \rho \big|u+\v  \frac{\mu_r}{\rho}\big|^2+ e(\rho,\bar{\rho})\Big)\,r^{N-1}\dd r
  +\v \int_0^\tau \int_a^{b} \frac{p_r}{\rho} \mu_r\,r^{N-1} \dd r \dd s\\
&= \int_a^{b}\Big(\frac{1}{2} \rho_0\big|u_0+\v \frac{\mu_{0r}}{\rho_0}\big|^2+ e(\rho_0,\bar{\rho})\Big)\,r^{N-1}\dd r,
\end{align*}
which, with \eqref{6.11}, leads to \eqref{6.12}.
$\hfill\Box$

\smallskip
\begin{lemma}\label{lem6.3}
For given $d$ and $D$ with $[d,D]\Subset [\d,b]$,
any smooth solution of \eqref{6.1}--\eqref{6.3} satisfies
\begin{align}\label{6.14}
\int_0^T\int_{K} \rho^{\g+1}(t,r)\, \dd r \dd t\leqq C(d,D,T,E_0),
\end{align}
where $K$ is any compact subset of $[d,D]$.
\end{lemma}

\noindent{\bf Proof.} We divide the proof into five steps.

\smallskip
1. Let $w(r)$ be a smooth compact support function with $\text{supp}\,w\subseteq[d,D]$ and $w(r)\equiv1$ for $r\in K$.
Multiplying $\eqref{6.1}_2$ by $w(r)$, we have
\begin{align}\label{3.32}
&(\r u w)_t+\big((\r u^2+p)w\big)_r+\f{N-1}r \r u^2 w\nonumber\\
&=\v\big((\mu+\lambda)(u_r+\f{N-1}r u) w\big)_r-\v \f{N-1}{r} u\mu_r w\nonumber\\
&\quad +\big(\r u^2+p-\v(\mu+\lambda)(u_r+\f{N-1}r u)\big)w_r.
\end{align}
Integrating \eqref{3.32}  over $[d,r)$ and multiplying the resultant equation by $\r w$,
we have
\begin{align}\label{3.33}
(\r^2 u^2+\r p)w^2
&=-\rho w \Big(\int_d^r\r u w \, \dd y\Big)_t -\rho w\int_d^r \f{N-1}y \r u^2 w\, \dd y\nonumber\\
&\quad+\rho w \int_d^r\Big(\r u^2+p-\v(\mu+\lambda)\big(u_y+\f{N-1}y u\big)\Big)w_y\, \dd y\nonumber\\
&\quad +\v \rho (\mu+\lambda)\big(u_r+\f{N-1}r u\big) w^2-\v\rho w \int_d^r \f{N-1}y u\mu_y w\, \dd y.
\end{align}
A direct calculation shows
\begin{align}\label{3.34}
\rho p w^2
&=-\Big(\r w\int_d^r\r u w\, \dd y\Big)_t-\Big(\r  uw\int_d^r\r u w\, \dd y\Big)_r+\r uw_r \int_d^r \r u w\, \dd y\nonumber\\
&\quad-\f{N-1}r \r u w \int_d^r \r u w\, \dd y -\rho w\int_d^r \f{N-1}y \r u^2 w\, \dd y\nonumber\\
&\quad+\rho w \int_d^r\Big(\r u^2+p-\v\frac{\mu+\lambda}{y}\big(yu_y+ (N-1)u\big)\Big)w_y\, \dd y\nonumber\\
&\quad -\v\rho w \int_d^r \f{N-1}y u\mu_y w\, \dd y +\v \rho (\mu+\lambda)\big(u_r+\f{N-1}r u\big) w^2\nonumber\\
&:=\sum_{j=1}^8 I_j.
\end{align}
To estimate the right-hand side of \eqref{3.34},
we first note from  \eqref{6.10-1} and \eqref{6.11} that
\begin{align}\label{3.35}
\int_d^D \r^\g\, r^{N-1} \dd r\leqq C(D, E_0).
\end{align}
Using \eqref{6.11} and  \eqref{3.35}, we see
\begin{align}
&\int_d^D\r\,  \dd r\leqq \f{C}{d^{N-1}}\int_d^D\r\, r^{N-1}\dd r
 \leqq C(d)\int_d^D (\r^{\g}+1)\,r^{N-1}\dd r\leqq C(d,D,E_0), \label{3.36}\\[1.5mm]
&\int_d^D \r u^2\, \dd r\leqq \f{C}{d^{N-1}}\int_d^D \rho u^2\,r^{N-1}\dd r
\leqq C(d,E_0).\label{3.37}
\end{align}

\medskip
2. Now it follows from \eqref{3.36}--\eqref{3.37} that
\begin{align}
&\left| \int_0^T\int_d^D I_1\, \dd r \dd t \right|
\leqq \int_d^D \Big(\Big|\Big(\r w\int_d^r\r u w\,\dd y\Big)(T,r)\Big|
  +\Big|\Big(\r w\int_d^r\r u w\, \dd y\Big)(0,r)\Big|\Big)\dd r\nonumber\\
&\qquad\qquad\qquad\quad\,\,\,\leqq C(d,D,T,E_0),\label{3.38}\\[1.5mm]
& \int_0^T\int_d^D I_2\, \dd r \dd t
  =\int_0^T\int_d^D \Big(\r u w\int_d^r\r u w\, \dd y\Big)_r\, \dd r \dd t =0, \label{3.39}\\[1.5mm]
&\left| \int_0^T\int_d^D I_3\, \dd r \dd t \right|
 =\left| \int_0^T\int_d^D \Big(\r u w_r\int_d^r\r u w\, \dd y\Big)\,  \dd r \dd t \right|
  \leqq C(d,D,T,E_0),
  \label{3.40}\\[1.5mm]
&\left| \int_0^T\int_d^D I_4\, \dd r \dd t \right|
\leqq C(d)\left| \int_0^T\int_d^D \Big(\r u\int_d^r\r u\,\dd y\Big)\,  \dd r\dd t \right|
\leqq C(d,D,T,E_0),
\label{3.41}\\[1.5mm]
&\left| \int_0^T\int_d^D I_5\, \dd r \dd t \right|
=\left| \int_0^T\int_d^D \Big(\r w \int_d^r\f{N-1}{y}\r u^2 w\, \dd y\Big)  \dd r \dd t \right|
\leqq C(d,D,T,E_0).
\label{3.42}
\end{align}

\medskip
3. We now estimate $I_6$. It follows from \eqref{6.11} that
\begin{align}
&\Big|\int_0^T\int_d^D \Big(\r w \int_d^r (\r u^2+p) w_y\,\dd y\Big) \dd r \dd t\Big|
\leqq C(d,D,T,E_0),\label{3.42-1}\\[2mm]
&\Big|\int_0^T\int_d^D \v \r w \Big( \int_d^r \frac{\rho+\alpha\d \rho^{\alpha}}{y}\big(yu_y+ (N-1)u\big) w_y\,\dd y\Big) \dd r \dd t \Big|\nonumber\\
&\leqq C(d,D,E_0) \Big\{\v\int_0^T \int_d^D (\rho+\delta \rho^\alpha) \big(u_r^2+ \frac{u^2}{y^2}+1\big)\,y^{N-1} \dd y \dd t\Big\}\nonumber\\
&\leqq C(d,D,T,E_0).\label{3.42-2}
\end{align}
Then it follows from \eqref{3.42-1}--\eqref{3.42-2} that
\begin{align}\label{3.43}
\Big|\int_0^T\int_d^D I_6\, \dd r \dd t\Big| \leqq C(d,D,T,E_0).
\end{align}

\medskip
4. For $I_7$, it follows from \eqref{6.11} and  integration by parts  that
\begin{align}
&\Big|\int_d^r \f{1}y u\mu_y w \,\dd y\Big|\nonumber\\
&\leqq \Big|\f1r (\mu u w)(t,r)\Big|
 +\Big|\int_d^r \f1y \mu \big(-\frac{1}{y}uw+ u_y w+ u w_y\big)(t,y)\,\dd y\Big|\nonumber\\
&\leqq \f1r \big((\rho +\d\rho^\alpha)|u w|\big)(t,r)
+C(d)\int_d^D \big(\rho u_r^2 +\delta  \rho^{\alpha} u_r^2\big)r^{N-1}\,\dd r
+C(d,D,E_0),\nonumber
\end{align}
which implies
\begin{align}\label{3.44}
&\Big|\int_0^T\int_d^DI_7\, \dd r \dd t\Big|\nonumber\\
&\leqq C(d,D,T,E_0)\Big(1+\v\int_0^T\int_d^D(\rho +\delta\rho^{\alpha})u_r^2\, r^{N-1}\dd r \dd t\Big)
+\v\int_0^T\int_d^D \rho^3 w^2\, \dd r \dd t\nonumber\\
&\leqq C(d,D,T,E_0)
+\v\int_0^T\int_d^D \rho^3 w^2\, \dd r \dd t,
\end{align}
where we have used $\alpha<1$.

\medskip
For $I_8$, it follows from \eqref{6.11} and the Cauchy inequality that
\begin{align}\label{3.46}
\Big|\int_0^T\int_d^D I_8\, \dd r \dd t\Big|
&\leqq \v\left|\int_0^T\int_d^D \r^2\big(u_r+\f{N-1}r u\big) w^2\, \dd r \dd t\right| \nonumber\\
&\quad + \v\delta\left|\int_0^T\int_d^D \r^{1+\alpha}\big(u_r+\f{N-1}r u\big) w^2\, \dd r \dd t\right|\nonumber\\
&\leqq C(d) \int_0^T\int_d^D \v (\rho+\delta\rho^\alpha)\big(u_r^2 +\frac{u^2}{r^2}\big)\,r^{N-1}\dd r \dd t \nonumber\\
&\quad  +\frac{\v}{2}\int_0^T\int_d^D \big(\r^3 +\rho^{2+\alpha}\big) w^2\, \dd r \dd t\nonumber\\
&\leqq C(d,D,T,E_0)
+\v\int_0^T\int_d^D \r^3 w^2\, \dd r \dd t.
\end{align}

To close the estimate, we still need to bound the last term on the right-hand sides of \eqref{3.44}--\eqref{3.46}.

We first consider the case: $\g\in(1,2]$. Notice that
\begin{align}\label{3.47}
&\v\int_0^T\int_d^D \r^3 w^2\, \dd r \dd t\nonumber\\
&\leqq \v \int_0^T \Big(\int_d^D \r^{\g}\, \dd r\Big) \sup_{r\in[d,D]}\big(\r^{3-\g} w^2\big)\,\dd t \nonumber\\
&\leqq C(d,D,E_0)
\int_0^T\v\sup_{r\in[d,D]}\big(\r^{3-\g} w^2\big) \dd t\nonumber\\
&\leqq \hat{C}(d,D,E_0)
\int_0^T\int_d^D \big(\v \r^{2-\g}|\r_r| w^2+\v \r^{3-\g} w |w_r|\big)\, \dd r \dd t,
\end{align}
where $\hat{C}(d,D,E_0)$ is a constant depending on $(d,D,E_0)$.
A direct calculation shows that
\begin{align}
\int_0^T\int_d^D \v \r^{2-\g}|\r_r| w^2\, \dd r \dd t
&\leqq \int_0^T\int_d^D \v\r^{\g-2} \r_r^2\,\dd r \dd t
  +\f\v2\int_{0}^T\int_d^D \r^{3(2-\g)} w^2\, \dd r \dd t\nonumber\\
&\leqq C(d,D,E_0)+\frac{\v}{2 \hat{C}(d,D,E_0)}
  \int_{0}^T\int_d^D \r^{3} w^2\, \dd r \dd t,\label{3.48}\\[2mm]
\int_0^T\int_d^D \v \r^{3-\g} w |w_r|\, \dd r \dd t
&\leqq \int_0^T\v \sup_{r}(\r w)(t,r) \Big(\int_d^D\r^{2-\g}|w _r|\,\dd r\Big)\, \dd t
\nonumber\\
&\leqq C(d,D,E_0)
\int_0^T\v\, \sup_{r}(\r w)(t,r)\,  \dd t\nonumber\\
&\leqq C(d,D,E_0)
\int_0^T \int_d^D \v\big(|\r_r| w+ \r |w_r|\big)\, \dd r \dd t
\nonumber\\&
\leqq C(d,D,E_0)
\Big(\int_0^T\int_d^D \big(\v\r^{\g-2} \r_r^2 + \r^{2-\g} w\big)\,\dd r \dd t
   +1\Big)\nonumber\\
&\leqq C(d,D,E_0).  \label{3.49}
\end{align}
Combining \eqref{3.47}--\eqref{3.49}, we have
\begin{align}\label{3.50}
\v\int_0^T\int_d^D \r^3 w^2\, \dd r \dd t\leqq C(d,D,E_0)\qquad\,\,\mbox{for $\g\in(1,2]$}.
\end{align}

For the case: $\g\in[2,3]$, notice that
\begin{align}\label{3.51}
&\v \int_0^T\int_d^D \r^3 w^2\, \dd r \dd t\nonumber\\
&\leqq \v \int_0^T\sup_{r\in[d,D]}\big(\r^2 w)\int_d^D \r w\, \dd r\, \dd t\nonumber\\
&\leqq C(d,D,E_0)
\int_0^T \int_d^D\big( \v\r|\r_r| w+\v\r^2 |w_r|\big)\, \dd r \dd t\nonumber\\
&\leqq C(d,D,E_0)
\int_0^T \int_d^D \big(\v^2\r^{\g-2}|\r_r|^2 w+\r^2 |w_r|+\r^{4-\g} w\big)\, \dd r \dd t\nonumber\\
&\leqq C(d,D,E_0).
\end{align}

For case $\gamma\in[3,\infty)$, it is direct to see
\begin{align}\label{3.51-1}
\int_0^T\int_d^D \rho^3 w^2\, \dd r \dd t\leqq C(d,D)\int_0^T\int_d^D\big( 1+ r^{N-1}\rho^\gamma\big)\,  \dd r \dd t
\leqq C(d,D,E_0).
\end{align}

\medskip
Now substituting \eqref{3.50}--\eqref{3.51-1}
into \eqref{3.44}--\eqref{3.46}, we obtain
\begin{align}\label{3.52}
\Big|\int_0^T\int_d^D (I_7+I_8)\, \dd r \dd t\Big|\leqq C(d,D,T,E_0).
\end{align}

\medskip
5. Integrating \eqref{3.34} over $[0,T]\times [d,D]$ and then using \eqref{3.38}--\eqref{3.42}, \eqref{3.43}, and \eqref{3.52},
we conclude \eqref{6.14}.
$\hfill\Box$

\section{Uniform Higher Integrability of the Approximate Solutions}

To employ the compensated compactness framework in \cite{Chen6},
we further require the higher integrability of the approximate solutions.

From now on, we denote
\begin{equation}\label{4.1a}
\begin{split}
&M_1:=E_0+\bar{\rho}+\bar{\rho}^{-1}+\delta^{-1}+\v^{-1}+\sup_{b\ge 1+\d^{-1}} E_2^{\v,\delta, b}<\infty,\\
&M_2:=M_1+\sup_{b\ge 1+\d^{-1}}\tilde{E}_0^{\v,\delta, b}<\infty,
\end{split}
\end{equation}
where
\begin{align}\label{4.1d}
\begin{split}
&E_2^{\v,\d,b}:=\int_\d^b  \rho_0 \Big(u_0^{2N} +
\big|\frac{\mu_{0r}}{\rho_0}\big|^{2N}\Big)\,r^{N-1}\dd r,\\
&\tilde{E}_0^{\v,\delta, b}:=\int_\d^b  \big(\f12 \rho_0 u_0^2+e(\rho_0,\bar{\r})\big)\,r^{2(N-1)+\vartheta}\dd r
\end{split}
\end{align}
for some $\vartheta\in (0,1)$.
From Lemma \ref{lem8.4}, we note that $E_2^{\v,\d,b}$ and $\tilde{E}_0^{\v,\delta, b}$ are uniformly bounded with respect to $b$,
while the upper bounds
may depend on $(\v,\delta)$,
so that $M_1$ and $M_2$ are finite for any fixed $(\v,\delta)$, independent of $b>0$.

\begin{proposition}\label{lem6.4}
Let $[d,D]\Subset [\d,b]$. Then the smooth solution of \eqref{6.1}--\eqref{6.3} satisfies
\begin{align}\label{6.15-1}
\int_0^T\int_d^D  \big(\rho|u|^3+\rho^{\g+\t}\big)(t,r)\,r^{N-1}\dd r\dd t
\leqq C(d,D,T, E_0)+C(T,M_2) b^{-\frac{\vartheta}{2}},
\end{align}
where $\vartheta\in (0,1)$ given in \eqref{4.1d}.
\end{proposition}

To prove \eqref{6.15-1}, we need to integrate the equations from the far-field, so that
the asymptotic behavior of $(\rho-\bar{\rho})(t,r)$ and $u(t,r)$ near boundary $r=b$ must be known.
Indeed, the key point for Proposition \ref{lem6.4} is that a decay rate of $(\rho-\bar{\rho})(t,r)$ and $u(t,r)$
can be derived,
and the positive constant $C(T, M_2)$
in \eqref{6.15-1} is independent of $b$ so that this term vanishes when $b\rightarrow\infty$.
In order to prove Proposition \ref{lem6.4}, it requires the following six lemmas.

To obtain the asymptotic behavior of $\rho(t,r)$ near boundary $r=b$,
we first need the lower and upper bounds of $\rho(t,r)$, which are independent of $b$.

\begin{lemma}[Upper Bound of the Density]\label{lem6.5}
There exists a constant $C(M_1)>0$ such that the smooth solution of \eqref{6.1}--\eqref{6.3} satisfies
\begin{align}\label{6.16}
0< \rho(t,r)\leqq C(M_1) \,\qquad\mbox{for $t\geqq 0$ and $r\in[\d,b]$}.
\end{align}
\end{lemma}

\noindent{\bf Proof.}
Notice that
\begin{align}
\begin{cases}
e(\rho,\bar{\r})\cong |\r-\bar{\r}|^2 \qquad \mbox{if $\r\in[\f{\bar{\r}}{2},  2\bar{\r}]$},\\[1mm]
e(\rho,\bar{\r})\cong |\r-\bar{\r}|^\gamma \qquad \mbox{if $\r\in\mathbb{R}_+\backslash[\f{\bar{\r}}{2},  2\bar{\r}]$}.\nonumber
\end{cases}
\end{align}
Denote
\begin{equation}\label{6.13}
A(t):=\{r\, :\, r\in[\d,b],~\rho(t,r)\geqq 2\bar{\r}\}
\end{equation}
with $A_1(t):=\{r\, :\, r\in[1,b],~ r\in A(t)\}\subset A(t)$  and $A_2(t):=A(t)\backslash  A_1(t)$.
It is easy to see that
\begin{align}\label{6.15}
e(\rho,\bar{\r})\geqq C(\bar{\r})^{-1}\qquad\,\, \mbox{for $r\in A(t)$},
\end{align}
which, along with \eqref{6.11}, yields
\begin{align}
E_0^{\v,\d,b}\geqq \int_\d^b e(\r,\bar{\r})\, r^{N-1} \dd r\geqq \int_{A_1(t)} e(\r,\bar{\r})\, \dd r
\geqq C(\bar{\r})^{-1} |A_1(t)|.\nonumber
\end{align}
Since $E_0^{\v,\d,b}\leqq C(E_0+1)$, we have
\begin{equation}\label{6.17}
|A(t)|\leqq |A_1(t)|+|A_2(t)|
\leqq C(\bar{\r},E_0).
\end{equation}
Since $\rho(t,r)$ is a continuous function on $[\delta,b]$, then, for any $r\in A(t)$, there exists $r_0\in A(t)$ such that
$\rho(t,r_0)=2\bar{\rho}$ and $|r-r_0|\leqq C(\bar{\r},E_0)$,
which implies
\begin{align}
\sqrt{\rho(t,r)}
&\leqq \sqrt{\rho(t,r_0)}+\int_{r_0}^r \frac{|\rho_y(t,y)|}{2\sqrt{\rho(t,y)}}\dd y\nonumber\\
&\leqq \sqrt{2\bar{\r}}+C(\bar{\r},E_0) \Big(\int_\d^b\frac{\rho_r^2}{\rho}\dd r\Big)^{\frac12}
\nonumber\\
&\leqq \sqrt{2\bar{\r}}+\frac{C(\bar{\r},E_0)}{\d^{\frac{N-1}{2}}\v}\nonumber\\
&\leqq C(\bar{\r},\v,\d,E_0).\nonumber
\end{align}
This completes the proof.
$\hfill\Box$

\begin{lemma}\label{lem6.7}
The smooth solution of \eqref{6.1}--\eqref{6.3} satisfies
\begin{align}\label{6.28}
\int_\d^b  \frac{\rho_r^{2N}}{\rho^{2N}}\,r^{N-1} \dd r
\leqq C(T,M_1)\qquad\,\, \mbox{for any $t\in [0,T]$}.
\end{align}
\end{lemma}

\noindent{\bf Proof.}
We divide the proof into three steps.

\medskip
1. We rewrite \eqref{6.9-4} as
\begin{align}\label{6.30}
(\v r^{N-1} \mu_x)_{\tau}=-u_{\tau}-r^{N-1}p_x
\end{align}
in the Lagrangian coordinates.
Integrating \eqref{6.30} over $[0,\tau]$ leads to
\begin{align}\label{6.31}
\v (r^{N-1} \mu_{x})(\tau,x)=\v (r^{N-1} \mu_{x})(0,x)-\big(u(x,\tau)-u_0(x)\big)-\int_0^\tau (r^{N-1}p_x)(s,x)\,\dd s.
\end{align}
Multiplying \eqref{6.31} by $(r^{N-1} \mu_{x})^{2N-1}$ and integrating the resultant equation yield
\begin{align*}
&\v\int_0^{L_b}\big|(r^{N-1} \mu_x)(\tau)\big|^{2N}\,\dd x\nonumber\\
&\leqq \left(\int_0^{L_b} \big|(r^{N-1} \mu_x)(\tau)\big|^{2N}\dd x\right)^{\frac{2N-1}{2N}}
\nonumber\\
&\quad\,\,\, \times
\Big\{
\big\|(u(\tau), u_0,(r^{N-1} \mu_x)(0))\big\|_{L^{2N}}
+C_T \big\|r^{N-1}(\rho^{\g})_x\big\|_{L^{2N}((0,\tau)\times (0, L_b))}\Big\},
\end{align*}
which leads to
\begin{align}\label{6.32-1}
&\int_0^{L_b}\big|(r^{N-1} \mu_x)(\tau)\big|^{2N}\dd x\nonumber\\
&\leqq C(\v)\Big\{
\big\|(u(\tau), u_0,(r^{N-1} \mu_x)(0))\big\|^{2N}_{L^{2N}}
+C_T\big\|r^{N-1}(\rho^{\g})_x\big\|^{2N}_{L^{2N}((0,\tau)\times (0, L_b))}\Big\}.
\end{align}
Notice that $|\mu_x|=\big|\big(\frac{1}{\alpha}\rho^{1-\alpha}+\delta\big)(\rho^\alpha)_x\big|\geqq \delta \big|(\rho^\alpha)_x\big|$
and $(\rho^\gamma)_x=\frac{\gamma}{\alpha}\rho^{\gamma-\alpha} (\rho^\alpha)_x$.
It follows from \eqref{6.16} and \eqref{6.32-1} that
\begin{align}\label{6.32-2}
&\int_0^{L_b}\big|\big(r^{N-1} (\rho^{\alpha})_x\big)(\tau)\big|^{2N}\dd x\nonumber\\
&\leqq C(T,\v,\d, E_0)\Big\{
\big\|(u(\tau),u_0,r^{N-1} \mu_x\big)(0))\big\|^{2N}_{L^{2N}}
\nonumber\\
&\qquad\qquad\qquad\,\,\quad\,
+\big\|r^{N-1}(\rho^{\alpha})_x\big\|^{2N}_{L^{2N}((0,\tau)\times (0, L_b))}\Big\}.
\end{align}
Plugging  \eqref{6.32-2} back to the Eulerian coordinates and noting $\alpha=\frac{2N-1}{2N}$, we see that, for $t\in[0,T]$,
\begin{align}\label{6.33}
&\int_\d^b \Big(\frac{\rho_r^{2N}}{\rho^{2N}}\Big)(t)\,r^{N-1} \dd r\nonumber\\
&\leqq C(T,\v,\d,E_0)\bigg\{ E_2^{\v,\d,b} +\int_\d^b  (\rho u^{2N})(t)\,r^{N-1}\dd r
\nonumber\\
&\qquad\qquad\qquad\quad\,\,\,
+\int_0^t\int_\d^b  \Big(\frac{\rho_r^{2N}}{\rho^{2N}}\Big)(s)\,r^{N-1}\dd r\dd s\bigg\}.
\end{align}

\medskip
2.
In order to close the above estimate, we need to bound $\int_\d^b \rho u^{2N}r^{N-1}\, \dd r$.
Multiplying $\eqref{6.1}_2$ by $r^{N-1} u^{2N-1}$ and then integrating by parts, we have
\begin{align}\label{6.33-1}
&\frac{1}{2N} \frac{\dd}{\dd t}\int_\d^b \rho u^{2N} \, r^{N-1}\dd r-\int_\d^b p (r^{N-1} u^{2N})_r\, \dd r\nonumber\\
&=-\v\int_\d^b\Big\{(\mu+\lambda) \big(u_r+\frac{N-1}{r} u\big) (r^{N-1} u^{2N-1})_r-(N-1) \mu (r^{N-2} u^{2N})_r\Big\}\, \dd r.
\end{align}
By similar arguments as in \eqref{6.11-4}--\eqref{6.11-5}, we obtain
\begin{align}\label{6.33-3}
&(\mu+\lambda) \big(u_r+\frac{N-1}{r} u\big) (r^{N-1} u^{2N-1})_r-(N-1) \mu (r^{N-2} u^{2N})_r\nonumber\\
&\geqq \rho u^{2N-2}\Big\{(2N-1)u_r^2+(N-1)\frac{u^{2}}{r^2}\Big\}r^{N-1}.
\end{align}

For the pressure term, it follows from \eqref{6.16} and the H\"{o}lder inequality that
\begin{align}\label{6.33-4}
&\left|\int_\d^b p (r^{N-1} u^{2N-1})_r \dd r\right|\nonumber\\
&=\left|\int_\d^b p \big((2N-1) r^{N-1} u^{2N-2}u_r +(N-1) r^{N-2} u^{2N-1}\big)\, \dd r\right|\nonumber\\
&\leqq \frac{\v}{2} \int_\d^b \rho u^{2N-2} \big( u_r^2+ \frac{u^2}{r^2}\big)\,r^{N-1}\dd r
   +C\int_\d^b  \rho^{2\gamma-1} u^{2N-2}\,r^{N-1}\dd r\nonumber\\
&\leqq \frac{\v}{2} \int_\d^b \rho u^{2N-2} \big(u_r^2+\frac{u^2}{r^2}\big)\,r^{N-1}\dd r
  +C(M_1) \Big(1+\int_\d^b \rho u^{2N}\,r^{N-1}\dd r  \Big).
\end{align}
Substituting \eqref{6.33-3}--\eqref{6.33-4} into \eqref{6.33-1}, we have
\begin{align}
\frac{\dd}{\dd t} \int_\d^b  \rho u^{2N}\,r^{N-1}\dd r
\leqq C(M_1)\Big(1+\int_\d^b \rho u^{2N}\, r^{N-1}\dd r  \Big),\nonumber
\end{align}
which, with the Gronwall inequality, implies
\begin{align}\label{6.33-7}
\int_\d^b  \rho u^{2N}\,r^{N-1}\dd r
\leqq C(T,M_1)\qquad\mbox{for $t\in[0,T]$}.
\end{align}

\medskip
3. Now substituting \eqref{6.33-7} into \eqref{6.33} yields
\begin{align}\label{6.33-8}
\int_\d^b \Big(\frac{\rho_r^{2N}}{\rho^{2N}}\Big)(t)\,r^{N-1}\dd r
\leqq C(T,M_1)\Big(1+\int_0^t\int_\d^b  \Big(\frac{\rho_r^{2N}}{\rho^{2N}}\Big)(s)\,r^{N-1}\dd r\dd s\Big).
\end{align}
Applying the Gronwall inequality to \eqref{6.33-8}, we conclude \eqref{6.28}.
$\hfill\Box$

\medskip
With the above preparation, we have the following lower bound of the density.

\begin{lemma}[Lower Bound of the Density]\label{lem6.8}
There exists $C(T,M_1)>0$ depending only on $(T,M_1)$ such that
the smooth solution of \eqref{6.1}--\eqref{6.3} satisfies
\begin{align}\label{6.34}
\rho(t,r) \geqq C(T,M_1)^{-1}>0\qquad\,\, \mbox{for $(t,r)\in [0,T]\times [\d,b]$}.
\end{align}
\end{lemma}

\noindent{\bf Proof.}
Define
\begin{equation}\label{2.76}
B(t):=\{ r\, :\, r\in[\d,b],\, 0\leqq \rho(t,r)\leqq \frac{\bar{\r}}{2}\}
\end{equation}
with $B_1(t):=\{r\,:\,r\in[1,b],\,r\in B(t)\}\subset B(t)$ and $B_2(t):=B(t)\backslash  B_1(t)$.
Similar to \eqref{6.15}--\eqref{6.17}, we have
\begin{equation}\label{6.37}
|B(t)|\leqq C(\bar{\r},E_0).
\end{equation}
Since $\rho(t,r)$ is a continuous function on $[\d,b]$, then, for any $r\in B(t)$, there exists $r_0\in B(t)$ such that
$\rho(t,r_0)=\frac{\bar{\rho}}{2}$ and $|r-r_0|\leqq C(\bar{\r},E_0)$.
Thus, for $\beta>0$,
\begin{align}
\rho(t,r)^{-\beta}&\leqq \rho(t,r_0)^{-\beta}
+\beta\Big|\int_{r_0}^r \rho^{-\beta-1}|\rho_r|\, \dd y\Big|\nonumber\\
&\leqq C(\bar{\r})+\beta\Big(\int_{\d}^b\frac{|\rho_r|^{2N}}{\rho^{2N}}\,\dd r\Big)^{\f{1}{2N}} \Big(\int_{B(t)} \rho^{-\f{2\beta N}{2N-1}}\,\dd r\Big)^{\f{2N-1}{2N}}\nonumber\\
&\leqq C(\bar{\r})+\beta  \hat{C}(T,M_1) \max_{r\in B(t)} \rho(t,r)^{-\beta},\nonumber
\end{align}
where \eqref{6.28} has been used in the last inequality.
Then we have
\begin{align}
\max_{r\in B(t)} \rho(t,r)^{-\beta}
\leqq C(\bar{\r})+\beta  \hat{C}(T,M_1)\max_{r\in B(t)}\rho(t,r)^{-\beta}.\nonumber
\end{align}
Taking $\beta>0$ small enough such that $\beta  \hat{C}(T,M_1)\leqq \f12$, we obtain
\begin{equation}\nonumber
\max_{r\in B(t)} \rho(t,r)^{-\beta}\leqq C(\bar{\r}).
\end{equation}
Therefore, we conclude
\begin{equation}\nonumber
\rho(t,r)\geqq C(\bar{\rho})^{-\frac1\beta}=C(T,M_1)^{-1}\qquad\,\, \mbox{for all $r\in B(t)$},
\end{equation}
which leads to
\eqref{6.34}.
$\hfill\Box$

\begin{remark}
Since $M_1$ is independent of $b$, the key point of Lemmas {\rm \ref{lem6.5}} and {\rm \ref{lem6.8}}
is that the lower and upper bounds of the density are independent of $b$.
\end{remark}

With the above lower and upper bounds of the density, even though they depend on $(\v, \d)$, we can have the following  weighted estimate:

\begin{lemma}\label{lem6.10}
Let $\vartheta\in(0,1)$ be some positive constant.	
Then the smooth solution of \eqref{6.1}--\eqref{6.3} satisfies
\begin{align}\label{6.42}
&\int_\d^b \big(\f12 \rho u^2+e(\rho,\bar{\r})\big)\, r^{2(N-1)+\vartheta}\dd r
 +\v\int_0^T\int_\d^b (\rho +\alpha\delta \rho^{\alpha}) u_r^2\, r^{2(N-1)+\vartheta}\dd r\dd s\nonumber\\[1.5mm]
& \leqq C(T,M_2).
\end{align}
\end{lemma}

\noindent{\bf Proof.} The proof consists of five steps.

\medskip
1. Let $L\in [0,N]$. Multiplying $\eqref{6.1}_2$ by $r^{N-1+L} u$ and then integrating by parts yield
\begin{align}\label{6.43}
&\frac{\dd}{\dd t} \int_\d^b \f12  \rho u^2\,r^{N-1+L} \dd r+\int_\d^b p_r u \,r^{N-1+L}\dd r\nonumber\\
&=\frac{L}{2} \int_\d^b  \rho u^3\,r^{N-2+L}\dd r\nonumber\\
&\quad-\v \int_\d^b (\mu+\lambda)\big(u_r+ \frac{N-1}{r} u\big) \big(u_r+\frac{N-1+L}{r} u\big)\,r^{N-1+L}\dd r\nonumber\\
&\quad+\v (N-1)\int_\d^b \mu u \big(2 u_r+\frac{N-2+L}{r} u\big)\,r^{N-2+L}\dd r.
\end{align}

\smallskip
2. It follows from integration by parts, \eqref{6.16}, and  \eqref{6.34} that
\begin{align}\label{6.44}
&\int_\d^b p_r u\, r^{N-1+L}\dd r\nonumber\\
&=\frac{\dd}{\dd t} \int_\d^b  e(\rho,\bar{\r})\,r^{N-1+L}\dd r
 -\frac{\kappa\gamma L}{\gamma-1} \int_\d^b  \rho u \big(\rho^{\gamma-1}-\bar{\r}^{\gamma-1}\big)\,r^{N-2+L}\dd r\nonumber\\
&\geqq -C(T,M_1) \int_\d^b \big(\rho u^2+e(\rho,\bar{\rho})\big)\,r^{N-2+L}\dd r
+\frac{\dd}{\dd t} \int_\d^b  e(\rho,\bar{\r})\,r^{N-1+L}\dd r.
\end{align}

Using the Sobolev inequality:
\begin{equation}\label{6.45-1}
\|u(t)\|_{L^\infty}\leqq  C \|u(t)\|_{L^2}^{\frac12}\,\|u_r(t)\|_{L^2}^{\frac12},
\end{equation}
we have
\begin{align}\label{6.45}
&\frac{L}{2}\left|\int_\d^b   \rho u^3\,r^{N-2+L}\dd r\right|\nonumber\\
&\leqq C\|u\|^{\f12}_{L^2}\|u_r\|^{\f12}_{L^2}\int_\d^b \rho u^2\,r^{N-2+L}\dd r \nonumber\\
&\leqq C(T,M_1) \bigg\{\int_\d^b  \rho u_r^2\, r^{N-1}\dd r+\Big(\int_\d^b  \rho u^2\,r^{N-2+L}\dd r\Big)^{\f43}\bigg\},
\end{align}
where we have used  \eqref{6.16}, \eqref{6.34}, and
$$
\|u\|_{L^2}\leqq C(T,M_1) \Big(\int_\d^b  \rho u^2\,r^{N-1}\dd r\Big)^{\frac12}
\leqq C(T,M_1).
$$

\medskip
3. For the viscous term, a direct calculation shows
\begin{align}\label{6.47}
&-(\mu+\lambda)\big(u_r+\frac{N-1}{r}u\big) \big( u_r+\frac{N-1+L}{r}u\big)r^{N-1+L} \nonumber\\
&\quad + (N-1) \mu u\big(2 u_r+ \frac{N-2+L}{r}u\big)r^{N-2+L}\nonumber\\
&\leqq -\frac12 (\rho + \alpha \delta  \rho^{\alpha})u_r^2 r^{N-1+L}
+ C(T,M_1) r^{N-3+L}\rho u^2.
\end{align}

\medskip
4. Substituting  \eqref{6.44} and  \eqref{6.45}--\eqref{6.47} into \eqref{6.43} yields
\begin{align}\label{6.48}
&\frac{\dd}{\dd t} \int_\d^b \big(\f12 \rho u^2+e(\rho,\bar{\r})\big)\,r^{N-1+L}\dd r
 +\frac{\v}{2} \int_\d^b (\rho +\alpha\delta \rho^{\alpha}) u_r^2\,r^{N-1+L}\dd r\nonumber\\[2mm]
&\leqq C(T,M_1)\Big\{\int_\d^b\big(\rho u^2+e(\rho,\bar{\rho})\big)\,r^{N-2+L}\dd r
 +\Big(\int_\d^b  \rho u^2r^{N-2+L}\,\dd r\Big)^{\f43}\nonumber\\
 &\qquad\qquad\qquad+\int_\d^b  \rho u_r^2\,r^{N-1}\dd r \Big\}.
\end{align}

\medskip
5. Taking $L=1$ in \eqref{6.48}, integrating the resultant inequality over $[0,t]$,
and using \eqref{6.11} yield
\begin{align*}
&\int_\d^b \big(\f12 \rho u^2+e(\rho,\bar{\r})\big)\, r^{N}\dd r
 +\frac{\v}{2}\int_0^t\int_\d^b (\rho +\alpha\delta \rho^{\alpha}) u_r^2\, r^{N}\dd r\dd s\\
&\leqq \int_\d^b \big(\f12 \rho_0 u_0^2+e(\rho_0,\bar{\r})\big)\,r^{N}\dd r
 +C(T,M_1)\\[1mm]
&\leqq C(T,M_2)\quad\qquad \mbox{for all $t\in[0,T]$}.
\end{align*}
Then,  taking $L=2,3, \cdots, N-1$, in \eqref{6.48} step by step, we have
\begin{align}\label{6.50}
&\int_\d^b \big(\f12 \rho u^2+e(\rho,\bar{\r})\big)\,r^{2N-2}\dd r
  +\frac{\v}{2}\int_0^t\int_\d^b (\rho +\alpha\delta \rho^{\alpha}) u_r^2\, r^{2N-2}\dd r\dd s\nonumber\\
&\leqq \int_\d^b \big(\f12 \rho_0 u_0^2+e(\rho_0,\bar{\r})\big)\,r^{2N-2}\dd r
  +C(T,M_2)\nonumber\\[1mm]
&\leqq C(T,M_2)
 \quad\qquad \mbox{for all $t\in[0,T]$}.
\end{align}
Finally, taking $L=N-1+\vartheta$ with $\vartheta\in(0,1)$ in \eqref{6.48}
and integrating it over $[0,t]$, then it follows from  \eqref{6.50} that
\begin{align*}
&\int_\d^b  \big(\f12 \rho u^2+e(\rho,\bar{\r})\big)\,r^{2N-2+\vartheta}\dd r
  +\v\int_0^t\int_\d^b (\rho +\alpha\delta \rho^{\alpha}) u_r^2\, r^{2N-2+\vartheta}\dd r\dd s\nonumber\\
&\leqq \int_\d^b \big(\f12 \rho_0 u_0^2+e(\rho_0,\bar{\r})\big)\,r^{2N-2+\vartheta}\dd r
   +C(T,M_2)\nonumber\\
&\leqq C(T,M_2)
 \quad\qquad \mbox{for all $t\in[0,T]$}.
\end{align*}
This completes the proof.
$\hfill\Box$

\medskip
\begin{lemma}[Decay Estimates]\label{lem6.11}
Any smooth solution of \eqref{6.1}--\eqref{6.3} satisfies that, for all $r\in[1,b]$,
\begin{align}
&|(\rho-\bar{\r})(t,r)|
\leqq C(T,M_2)r^{-\frac34(N-1)-\f\vartheta4},\label{6.52}\\[2mm]
&\int_0^T \big(|u(t,r)|+ |u(t,r)|^3\big)\,\dd t
\leqq C(T,M_2) r^{-N+1-\frac{\vartheta}{2}}
\quad\,\, \mbox{for any $T>0$}.
\label{6.53}
\end{align}
\end{lemma}

\noindent{\bf Proof.}
It follows from \eqref{6.12}, \eqref{6.16}, \eqref{6.34}, and \eqref{6.42} that, for all $t\in [0,T]$,
\begin{align}
&\int_1^b \big((|(\rho-\bar{\r})(t,r)|^2+|u(t,r)|^2) r^{N-1+\vartheta}+ |\rho_r(t,r)|^2\big)\,r^{N-1}\dd r\nonumber\\
&
+\int_0^T\int_1^b |u_r(t,r)|^2\,r^{2(N-1)+\vartheta}\dd r\dd t
\leqq C(T,M_2).
\label{6.55}
\end{align}
For any $r\in [n,n+1]\cap [1,b]$ with $n+1\leqq [b]$, it follows from \eqref{6.55} and the Sobolev inequality that
\begin{align}
|(\rho-\bar{\r})(t,r)|^2
&\leqq 2\Big(\int_n^{n+1}|(\rho-\bar{\r})(t,r)|^2\dd r\Big)^{\f12}\Big(\int_n^{n+1}|\rho_r(t,r)|^2\dd r\Big)^{\f12}\nonumber\\
&\quad+\int_n^{n+1}|(\rho-\bar{\r})(t,r)|^2\dd r\nonumber\\
&\leqq Cn^{-\frac32(N-1)-\f\vartheta2}\Big(\int_n^{n+1}|(\rho-\bar{\r})(t,r)|^2\,r^{2(N-1)+\vartheta}\dd r\Big)^{\f12}\nonumber\\
&\qquad\times  \Big(\int_n^{n+1}|\rho_r(t,r)|^2\,r^{N-1}\dd r\Big)^{\f12}\nonumber\\
&\quad +n^{-2(N-1)-\vartheta} \int_n^{n+1}|(\rho-\bar{\r})(t,r)|^2\,r^{2(N-1)+\vartheta}\dd r\nonumber\\
&\leqq C(T,M_2) r^{-\frac32(N-1)-\f\vartheta2}.\nonumber
\end{align}

\smallskip
Similarly, for $r\in [n,n+1]\cap [1,b]$ with $n+1\leqq [b]$, it follows from \eqref{6.45-1} and \eqref{6.55} that
\begin{align}
|u(t,r)|^2
&\leqq C r^{-2(N-1)-\vartheta}\int_n^{n+1}|u(t,r)|^2\, r^{2(N-1)+\vartheta}\dd r\nonumber\\
&\quad
+ Cr^{-2(N-1)-\vartheta} \bigg(\int_n^{n+1}|u(t,r)|^2\,r^{2(N-1)+\vartheta}\dd r\bigg)^{\f12}\nonumber\\
&\qquad\times  \bigg(\int_n^{n+1}|u_r(t,r)|^2\,r^{2(N-1)+\vartheta}\dd r\bigg)^{\f12}\nonumber\\
& \leqq C(T,M_2) r^{-2(N-1)-\vartheta}
  \bigg(\Big(\int_n^{n+1}|u_r(t,r)|^2\,r^{2(N-1)+\vartheta}\dd r\Big)^{\f12}+1\bigg),\nonumber
\end{align}
which yields
\begin{align}\label{6.57}
&|u(t,r)|+|u(t,r)|^{3}\nonumber\\
&\leqq C(T,M_2) r^{-N+1-\frac\vartheta2}
\bigg(\int_n^{n+1}|u_r(t,r)|^2\,r^{2(N-1)+\vartheta}\dd r +1\bigg).
\end{align}
Integrating \eqref{6.57} over $[0,T]$, we obtain
\begin{align}\nonumber
\int_0^T \big(|u(t,r)|+ |u(t,r)|^3\big)\,\dd t
\leqq C(T,M_2) r^{-N+1-\frac\vartheta2} \qquad \mbox{for any $r\in [1, [b]]$}.
\end{align}

Finally, we consider the case that $r\in [b-1,b]$. Then, by the same arguments as above,
we see that, for $r\in [b-1,b]$,
\begin{align*}
\begin{split}
&|(\rho-\bar{\rho})(r)| \leqq C(T, M_2) b^{-\frac{3}{2}(N-1)-\frac{\vartheta}{2}},\\
&\int_0^T \big(|u(t,r)|+ |u(t,r)|^3\big) \dd t
\leqq C(T,M_2) b^{-N+1-\frac\vartheta2}.
\end{split}
\end{align*}
Combining all the above estimates, we prove \eqref{6.52}--\eqref{6.53}.
This completes the proof.
$\hfill\Box$

\bigskip
Choosing $\psi(s)=\frac12 s|s|$ in \eqref{weakentropy} leads to the corresponding entropy pair as
\begin{align}\label{6.58}
\begin{cases}
\displaystyle
\eta^{\#}(\r,m)=\f12 \r \int_{-1}^1 (u+\r^{\t} s) |u+\r^{\t}s| [1-s^2]_+^{\hl}\,\dd s,\\[3mm]
\displaystyle
q^{\#}(\r, m)=\f12 \r \int_{-1}^1 (u+\theta\r^{\t}s)(u+\r^{\t} s) |u+\r^{\t}s|\,[1-s^2]_+^{\hl}\,\dd s,
\end{cases}
\end{align}
where $\hl=\frac{3-\gamma}{2(\gamma-1)}$, $\theta=\frac{\gamma-1}{2}$, and $m=\rho u$ as indicated earlier. Then a direct calculation shows
\begin{align}\label{6.59}
|\eta^{\#}(\r, m)|
\leqq C_\gamma \big(\r |u|^2+\r^{\g}\big), \quad
q^{\#}(\r, m)
\ge C_\gamma^{-1}
\big(\r |u|^3+\r^{\g+\t}\big),
\end{align}
where and whereafter $C_\gamma>0$ is a universal constant depending only on $\gamma>1$.

Moreover, notice that
\begin{align}\label{6.60}
\begin{split}
\partial_\r\eta^{\#}&=\int_{-1}^1 \big(-\f12 u+(\t+\f12)\r^{\t} s\big)|u+\r^{\t}s|\, [1-s^2]_+^{\hl}\, \dd s,\\
\partial_m\eta^{\#}&=\int_{-1}^1|u+\r^{\t}s|\,[1-s^2]_+^{\hl}\,\dd s.
\end{split}
\end{align}
Then
\begin{align}\label{6.61}
\begin{split}
&|\eta^{\#}_m|\leqq C_\gamma(|u|+\rho^\theta),\quad  |\eta^{\#}_\rho|\leqq C_\gamma(|u|^2+\rho^{2\theta}),\\
&\eta^{\#}_\rho(\rho,0)=0,\quad  \eta^{\#}_m(\rho,0)=2{\rho}^{\theta} \int_{0}^1 s[1-s^2]_+^{\hl}\dd s.
\end{split}
\end{align}
Now we define the relative entropy pair as
\begin{align}\label{6.77}
\begin{cases}
\displaystyle\tilde{\eta}(\rho,m)=\eta^{\#}(\r,m)-\eta^{\#}(\bar\r,0)- \eta^{\#}_m(\bar\r,0) m,\\[1.5mm]
\displaystyle\tilde{q}(\rho,m)=q^{\#}(\r,m)-q^{\#}(\bar\r,0)-\eta^{\#}_m(\bar\r,0) \big(\frac{m^2}{\rho}+p(\rho)-p(\bar{\rho})\big).
\end{cases}
\end{align}

With these,
we have the following useful lemma.

\begin{lemma}\label{lem6.14}
The relative entropy pair $(\tilde{\eta}, \tilde{q})$ satisfies
\begin{align}\label{6.78}
  m\partial_\rho \tilde{\eta}(\rho, m)+\frac{m^2}{\rho}\partial_m \tilde{\eta}(\rho,m)-\tilde{q}(\rho, m)
 \leqq C_\gamma(\bar{\rho}) \big(\frac{m^2}{\rho}+ e(\rho,\bar{\rho})\big),
\end{align}
where $C_\gamma(\bar{\rho})>0$ is a positive constant depending only on $(\gamma, \bar{\rho})$.
\end{lemma}

\noindent{\bf Proof}.
The estimate for \eqref{6.78} is very subtle,
which will be used to overcome the singularity from the far-field
in the M-D case, different from the 1-D case.
The proof is divided into three steps.

\medskip
1. {\it Claim}: $(\eta^{\#}, q^{\#})$ satisfies
\begin{align}\label{6.63}
& m \partial_\rho \eta^{\#}(\rho, m)+\frac{m^2}{\rho}\partial_m \eta^{\#}(\rho,m)-q^{\#}(\rho,m)\nonumber\\
&\leqq \min\big\{0,-q^{\#}(\rho,0)+C_{\gamma}\rho^{\theta-1} m^2\big\},
\end{align}
where $\displaystyle q^{\#}(\rho,0)=\theta\rho^{3\theta+1}\int_0^1s^3[1-s^2]_+^{\hl}\,\dd s$.

A direct calculation shows that
\begin{align}\label{6.64}
&m\partial_\rho \eta^{\#}(\rho,m)+\frac{m^2}{\rho}\partial_m \eta^{\#}(\rho,m)-q^{\#}(\rho,m)\nonumber\\
&=\frac{\theta}{2} \rho^{1+\theta} \int_{-1}^1 (u-\rho^{\theta}s)s |u+\rho^{\theta}s| [1-s^2]_+^{\hl}\,\dd s.
\end{align}
Now we divide the proof into three cases:

\smallskip
\noindent{\it Case} 1. $u\geqq 0$ and $|u|\geqq \rho^{\theta}$. For this case, $u+\rho^{\theta} s\geqq 0$ for $s\in[-1,1]$.
Then
\begin{align} \label{6.69}
 m\partial_\rho \eta^{\#}+\frac{m^2}{\rho}\partial_m \eta^{\#}-q^{\#}=0.
\end{align}
On the other hand, we have
\begin{align}\label{6.65}
&m\partial_\rho \eta^{\#}+ \frac{m^2}{\rho}\partial_m \eta^{\#}-q^{\#}
=0=-q^{\#}(\rho,0)+q^{\#}(\rho,0)\nonumber\\
&=-q^{\#}(\rho,0)+\theta \int_0^1s^3[1-s^2]_+^{\hl}\,\dd s\rho^{1+3\theta}\nonumber\\
&\leqq -q^{\#}(\rho,0)+C_\gamma \rho^{\theta-1} m^2,
\end{align}
where we have used  that $\rho^{\theta}\leqq |u|$ in the last inequality.

\medskip
\noindent{\it Case} 2. $u\geqq 0$ and $|u|< \rho^{\theta}$.
For this case,  $s_0:=-\frac{u}{\rho^{\theta}}\in (-1,0]$,
which implies that $u^2-s^2\rho^{2\theta}\leqq 0$ for $s\geqq|s_0|$.
Then
\begin{align}\label{6.66}
 m\partial_\rho \eta^{\#}+ \frac{m^2}{\rho} \partial_m \eta^{\#}-q^{\#}
&=\theta \rho^{1+\theta} \int_{|s_0|}^1(u^2-s^2\rho^{2\theta}) s [1-s^2]_+^{\hl}\,\dd s
\leqq 0.
\end{align}
On the other hand, we have
\begin{align}\label{6.67}
&m\partial_\rho \eta^{\#}+\frac{m^2}{\rho}\partial_m \eta^{\#}-q^{\#}\nonumber\\
&=\theta \rho^{1+\theta} \int_{|s_0|}^1(u^2-s^2\rho^{2\theta}) s [1-s^2]_+^{\hl}\,\dd s\nonumber\\
&=\theta \rho^{1+\theta} u^2\int_{|s_0|}^1s [1-s^2]_+^{\hl}\,\dd s-\theta \rho^{1+3\theta} \int_{0}^1s^3 [1-s^2]_+^{\hl}\,\dd s\nonumber\\
&\quad  +\theta \rho^{1+3\theta} \int_0^{|s_0|}s^3[1-s^2]_+^{\hl}\,\dd s\nonumber\\
&\leqq -q^{\#}(\rho,0)+C_\gamma \rho^{1+\theta}u^2+C_\gamma \rho^{1+3\theta} |s_0|^2\nonumber\\
&\leqq  -q^{\#}(\rho,0)+C_\gamma \rho^{\theta-1}m^2.
\end{align}

\medskip
\noindent{\it Case} 3. $u\leqq 0$. Similar to \eqref{6.69}--\eqref{6.67}, we
also obtain \eqref{6.63}.

\medskip
Combining Cases 1--3, we conclude the claim for
\eqref{6.63}.

\medskip
2. Claim: $(\eta^{\#}, q^{\#})$ satisfies
\begin{align}
	&\eta^{\#}_m(\bar{\rho},0)(p(\rho)-p(\bar{\rho}))-q^{\#}(\rho,0)+q^{\#}(\bar{\rho}, 0)\nonumber\\
	& =2\bar{\rho}^{\theta}\int_0^1s[1-s^2]_+^{\hl}\,\dd s\, \big(p(\rho)-p(\bar{\rho})-p'(\bar{\rho})(\rho-\bar{\rho})\big)\nonumber\\
	&\quad-\frac{4\theta^2}{3\gamma-1} \int_0^1s[1-s^2]_+^{\hl}\,\dd s
        \,\big(\rho^{1+3\theta}-\bar{\rho}^{1+3\theta}-(1+3\theta) \bar{\rho}^{3\theta} (\rho-\bar{\rho})\big),\label{6.72-1}\\[1.5mm]
	&\eta^{\#}_m(\bar{\rho},0) \big(p(\rho)-p(\bar{\rho})\big)+ q^{\#}(\bar{\rho}, 0)\nonumber\\
&     =\int_0^1s[1-s^2]_+^{\hl}\,\dd s
      \,\Big(2\bar{\rho}^{\theta}p(\rho)-\frac{4\theta^3}{\gamma(3\gamma-1)} \bar{\rho}^{\gamma+\theta}\Big).\label{6.72-2}
	\end{align}

\smallskip
A direct calculation shows that
\begin{align}\label{6.73}
\int_0^1s^3[1-s^2]_+^{\hl}\,\dd s
&=\frac12B(2,1+\hl)=\frac{1}{2(2+\hl)} B(1,1+\hl)\nonumber\\
&=\frac{1}{2+\hl}\int_0^1s[1-s^2]_+^{\hl}\,\dd s,
\end{align}
where we have used the properties of the beta function $B(\cdot,\cdot)$.
Using \eqref{6.73}, we have
\begin{align}\label{6.74}
&\eta^{\#}_m(\bar{\rho},0) \big(p(\rho)-p(\bar{\rho})\big)-q^{\#}(\rho,0)+q^{\#}(\bar{\rho}, 0)\nonumber\\[1mm]
&=2\bar{\rho}^{\theta} \int_0^1s[1-s^2]_+^{\hl}\,\dd s \big( p(\rho)-p(\bar{\rho})-p'(\bar{\rho})(\rho-\bar{\rho}) \big)\nonumber\\
&\quad-\frac{\theta}{2+\hl}\int_0^1s[1-s^2]_+^{\hl}\,\dd s
  \Big(\rho^{1+3\theta}-\bar{\rho}^{1+3\theta}-\gamma \kappa \frac{2(2+\hl)}{\theta}  \bar{\rho}^{3\theta} (\rho-\bar{\rho})\Big).
\end{align}
Combining
$2+\hl=\frac{3\gamma-1}{4\theta}$ and $\gamma \kappa \frac{2(2+\hl)}{\theta}=1+3\theta$
with \eqref{6.74}, we conclude \eqref{6.72-1}.

\smallskip
For \eqref{6.72-2}, we note that
\begin{align}
&\eta^{\#}_m(\bar{\rho},0) \big(p(\rho)-p(\bar{\rho})\big)+q^{\#}(\bar{\rho}, 0)\nonumber\\
&=\int_0^1s[1-s^2]_+^{\hl}\,\dd s\,\Big(2 \bar{\rho}^{\theta} p(\rho)+(\frac{\theta}{2+\hl}-2\kappa) \bar{\rho}^{\gamma+\theta} \Big)\nonumber\\
&=\int_0^1s[1-s^2]_+^{\hl}\,\dd s\,\Big(2 p(\rho)-\frac{4\theta^3}{\gamma(3\gamma-1)} \bar{\rho}^{\gamma} \Big) \bar{\rho}^{\theta},\nonumber
\end{align}
which implies \eqref{6.72-2}.

\smallskip
3. Noting \eqref{6.10-1} and \eqref{6.13}, we have
\begin{align}\label{6.124}
e(\rho,\bar{\rho}) I_{A(t)}(r)&\geqq C_\gamma \rho (\rho^{\theta}-\bar{\rho}^{\theta})^2 \,I_{A(t)}(r)\nonumber\\
&\geqq C_\gamma \rho (1-\frac1{2^{\theta}})^2\rho^{2\theta}\, I_{A(t)} (r)\nonumber\\
&\geqq C_\gamma p(\rho)\, I_{A(t)} (r).
\end{align}

If $r\in A(t)$, then it follows from \eqref{6.63} and \eqref{6.72-2} that
\begin{align}\label{6.79}
&\big(m\partial_\rho \tilde{\eta}(\rho, m)+\frac{m^2}{\rho}\partial_m \tilde{\eta}(\rho,m)-\tilde{q}(\rho,m)\big) I_{A(t)}(r)\nonumber\\
&= \Big\{\big(m \partial_\rho \eta^{\#}(\rho,m)+\frac{m^2}{\rho} \partial_m \eta^{\#}(\rho,m)-q^{\#}(\rho, m)\big)\nonumber\\
&\qquad +\eta^{\#}_m(\bar{\rho},0) \cdot (p(\rho)-p(\bar{\rho}))+q^{\#}(\bar{\rho}, 0)\Big\}I_{A(t)}(r)\nonumber\\
&\leqq \int_0^1s[1-s^2]_+^{\hl}\,\dd s\,\Big(2\bar{\rho}^{\theta}p(\rho)-\frac{4\theta^3}{\gamma(3\gamma-1)} \bar{\rho}^{\gamma+\theta}\Big)I_{A(t)}(r)\nonumber\\[1mm]
&\leqq C_{\gamma}(\bar{\rho})  p(\rho) I_{A(t)}(r)\nonumber\\[1mm]
&\leqq C_{\gamma}(\bar{\rho}) e(\rho,\bar{\rho})\, I_{A(t)}(r),
\end{align}
where \eqref{6.124} has been used in the last inequality.

On the other hand, for $r\in A^c(t)=[\delta, b]\setminus A(t)$, it follows from \eqref{6.63} and \eqref{6.72-1} that
\begin{align}\label{6.80}
&\big( m\partial_\rho \tilde{\eta}(\rho, m)+\frac{m^2}{\rho}\partial_m \tilde{\eta}(\rho,m)-\tilde{q}(\rho,m)\big)I_{A^c(t)}(r)\nonumber\\
&=\Big\{ \big(m\partial_\rho \eta^{\#}(\rho,m)+\frac{m^2}{\rho}\partial_m \eta^{\#}(\rho,m)-q^{\#}(\rho,m)\big)\nonumber\\
&\qquad+\eta^{\#}_m(\bar{\rho},0) \, (p(\rho)-p(\bar{\rho}))+q^{\#}(\bar{\rho}, 0)\Big\}I_{A^c(t)}(r) \nonumber\\
&\leqq \Big\{q^{\#}(\bar{\rho}, 0)+\eta^{\#}_m(\bar{\rho},0) \,(p(\rho)-p(\bar{\rho}))-q^{\#}(\rho,0)\Big\}I_{A^c(t)}(r)
  +C_{\gamma}\rho^{\theta-1} m^2I_{A^c(t)}(r) \nonumber\\
&= 2\bar{\rho}^{\theta}\int_0^1s[1-s^2]_+^{\hl}\,\dd s\,\big(p(\rho)-p(\bar{\rho})-p'(\bar{\rho})(\rho-\bar{\rho})\big)I_{A^c(t)}(r)\nonumber\\
&\quad-\frac{4\theta^2}{3\gamma-1} \int_0^1s[1-s^2]_+^{\hl}\,\dd s\,\big(\rho^{1+3\theta}-\bar{\rho}^{1+3\theta}-(1+3\theta) \bar{\rho}^{3\theta} (\rho-\bar{\rho})\big)I_{A^c(t)}(r)\nonumber\\
&\quad +C_{\gamma}\rho^{\theta-1} m^2I_{A^c(t)}(r)\nonumber\\
&\leqq 2\bar{\rho}^{\theta}\int_0^1s[1-s^2]_+^{\hl}\,\dd s\,\big(p(\rho)-p(\bar{\rho})-p'(\bar{\rho})(\rho-\bar{\rho})\big)I_{A^c(t)}(r)\nonumber\\
&\quad + C_{\gamma}\rho^{\theta-1} m^2I_{A^c(t)}(r)\nonumber\\[1mm]
&\leqq C_\gamma \big(\rho^{\theta-1} m^2+e(\rho,\bar{\rho})\big)I_{A^c(t)}(r)\nonumber\\
&\leqq C_\gamma(\bar{\rho}) \big(\frac{m^2}{\rho}+e(\rho,\bar{\rho})\big)I_{A^c(t)}(r),
\end{align}
where we have used \eqref{6.10} and  $\rho^{\theta}(t,r)\leqq (2\bar{\rho})^{\theta}$ for $r\in A^c(t)$.
Combining \eqref{6.79} with \eqref{6.80}, we conclude \eqref{6.78}.
$\hfill\Box$

\smallskip
Now we are in the position to prove the key estimate, Proposition \ref{lem6.4}.
\vspace{2.5mm}

\noindent{\bf Proof of Proposition \ref{lem6.4}.}
We divide the proof into six steps.

\smallskip
1. For $\tilde{\eta}(\rho, m)$ defined in \eqref{6.77},
we multiply $\eqref{6.1}_1$ by $r^{N-1}\partial_\rho \tilde{\eta}(\rho, m)$
and $\eqref{6.1}_2$ by $r^{N-1}\partial_m \tilde{\eta}(\rho, m)$ to obtain
\begin{align}\label{6.81}
&(r^{N-1} \tilde{\eta})_t+(r^{N-1}\tilde{q})_r
  +(N-1)r^{N-2}\big(-\tilde{q}+ m\partial_\rho \tilde{\eta}+\frac{m^2}{\rho}\partial_m \tilde{\eta}\big)\nonumber\\
&=\v r^{N-1} \partial_m \tilde{\eta} \Big\{\big((\r+\alpha \delta\r^{\alpha})(u_r+\frac{N-1}{r}u)\big)_r-\frac{N-1}{r}  (\rho+\d \rho^{\alpha})_r u \Big\}.
\end{align}
Let $y\in[b-1,b]$ and $r\in[d,D]$.  Integrating \eqref{6.81} over $[r,y]$ leads to
\begin{align}\label{6.82}
&\tilde{q}(t,r)r^{N-1}\nonumber\\
&=\frac{\dd}{\dd t} \int_r^y  \tilde{\eta}(t,z)\,z^{N-1}\dd z+\tilde{q}(t,y)y^{N-1} \nonumber\\
&\quad
+(N-1)\int_r^y \big(-\tilde{q}+ m\partial_\rho \tilde{\eta}+\frac{m^2}{\rho}\partial_m \tilde{\eta}\big)(t,z)\,z^{N-2}\dd z\nonumber\\
&\quad -\v \int_r^y  \partial_m \tilde{\eta}\,\Big\{ \big((\r+\alpha \delta\r^{\alpha})(u_z+\frac{N-1}{z}u)\big)_z
       -\frac{N-1}{z}  (\rho+\d \rho^{\alpha})_z u \Big\}\,z^{N-1}\dd z.
\end{align}
Integrating \eqref{6.82} over $[0,T]\times[b-1,b]\times [d,D]$, we have
\begin{align}\label{6.83}
&\int_0^T\int_d^D\tilde{q}(t,r)\,r^{N-1}\dd r\dd t\nonumber\\
&=(N-1)\int_0^T\int_{b-1}^b\int_d^D\int_r^y \big(m\partial_\rho \tilde{\eta}+\frac{m^2}{\rho}\partial_m \tilde{\eta}-\tilde{q}\big)(t,z)\,z^{N-2}\dd z\dd r\dd y\dd t\nonumber\\
&\quad+\int_{b-1}^b\int_d^D\int_r^y  \big(\tilde{\eta}(T,z) -\tilde{\eta}(0,z)\big)\,z^{N-1} \dd z \dd r \dd y\nonumber\\
&\quad +(D-d)\int_0^T\int_{b-1}^b \tilde{q}(t,y)\,y^{N-1}\dd y\dd t\nonumber\\
&\quad-\v\int_0^T\int_{b-1}^b\int_d^D\int_r^y  \partial_m \tilde{\eta}\,\Big\{ \big((\r+\alpha \delta\r^{\alpha})(u_z+\frac{N-1}{z}u)\big)_z\nonumber\\
&\qquad\qquad\qquad\qquad\qquad\qquad\,\,\,\,\,\,
  -\frac{N-1}{z}  (\rho+\d \rho^{\alpha})_z u \Big\}\,z^{N-1}\dd z\dd r\dd y\dd t\nonumber\\
&=:\sum_{j=1}^4 J_j.
\end{align}

\medskip
2. For $J_1$ in \eqref{6.83}, it follows from \eqref{6.11} and Lemma \ref{lem6.14}  that
\begin{align}\label{6.89}
J_1 &\leqq C_\gamma(\bar{\rho})\frac{D}{d}\int_0^T\int_d^b \big(\rho u^2+ e(\rho,\bar{\rho})\big)(t,z)\,z^{N-1} \dd z\dd t\nonumber\\[1mm]
&\leqq C_\gamma(\bar{\rho})\frac{DT}{d}(E_0+1).
\end{align}

\medskip
3. For $J_2$ in \eqref{6.83}, we first note that
$|\partial_{mm}\eta^{\#}(\rho, m)|\leqq \frac2\rho \int_0^1[1-s^2]_+^{\hl}\, \dd s$.
This, combining \eqref{6.60} and \eqref{6.61} with the Taylor expansion of $\eta^{\#}(\rho,m)$ at $m=0$, yields
\begin{align}\label{2.85}
\eta^{\#}(\rho, m)=2\int_0^1s[1-s^2]_+^{\hl}\, \dd s\, \rho^{\theta}m +R_1(\rho, m)
\end{align}
with
\begin{equation}\label{2.86}
|R_1(\rho,m)|\leqq C_\gamma \frac{m^2}{\rho}.
\end{equation}
Then it follows from \eqref{6.10-1}, \eqref{6.61}--\eqref{6.77}, and \eqref{2.85}--\eqref{2.86} that
\begin{align}
|\tilde{\eta}(\rho,m)|&\leqq 2\int_0^1s[1-s^2]_+^{\hl}\, \dd s \, |m(\rho^{\theta}-\bar{\rho}^{\theta})|+|R_1(\rho,m)|
\leqq C_\gamma \big(\frac{m^2}{\rho}+e(\rho,\bar{\rho})\big),\nonumber
\end{align}
which, along with \eqref{6.11}, implies
\begin{align}\label{6.87}
|J_2|=\left|\int_{b-1}^b\int_d^D\int_r^y \big(\tilde{\eta}(T,z) -\tilde{\eta}(0,z)\big)\,z^{N-1}\dd z \dd r \dd y\right|
\leqq C_\gamma D(E_0+1).
\end{align}

\medskip
4. For the third term $J_3$ in \eqref{6.83}, we need to use the decay properties obtained in Lemma \ref{lem6.11}. A direct calculation shows
that
$$
|q^{\#}(\rho, m)-q^{\#}(\rho,0)|\leqq C_\gamma \big(\frac{|m|^3}{\rho^2}+\rho^{2\theta}|m|\big),
$$
which, with \eqref{6.72-1}, yields
\begin{align}\label{6.91}
\tilde{q}(\rho,m)
&=\frac{4\theta^2}{3\gamma-1} \int_0^1s[1-s^2]_+^{\hl}\,\dd s \, \big(\rho^{1+3\theta}-\bar{\rho}^{1+3\theta}-(1+3\theta) \bar{\rho}^{3\theta} (\rho-\bar{\rho})\big)\nonumber\\[2mm]
&\quad-2\bar{\rho}^{\theta}\int_0^1s[1-s^2]_+^{\hl}\,\dd s\, \big(\rho u^2+p(\rho)-p(\bar{\rho})-p'(\bar{\rho})(\rho-\bar{\rho})\big)\nonumber\\
&\quad+\big(q^{\#}(\rho,\rho u)-q^{\#}(\rho,0)\big)\nonumber\\[2mm]
&\leqq C(T,M_2)\,\big(|\rho-\bar{\rho}|^2+|u|^3+|u|\big),
\end{align}
where we have used the Taylor expansion, \eqref{6.16}, and \eqref{6.34} in the last inequality.
Now it follows from \eqref{6.91} and Lemma \ref{lem6.11} that
\begin{align}\label{6.92}
|J_3|&\leqq C(T,M_2) D
  \int_{b-1}^b\int_0^T  \big(|\rho-\bar{\rho}|^2+|u|^3+|u|\big)(t,y)\,y^{N-1}\dd y\dd t\nonumber\\[1mm]
&  \leqq C(T,M_2) b^{-\frac{\vartheta}{2}}.
\end{align}

\medskip
5. For $J_4$ in \eqref{6.83}, we regard $\tilde\eta_m(\rho,\rho u)$ as a function of $(\rho, u)$ to obtain
\begin{equation}\label{6.93}
\begin{split}
&|\partial_{m}\tilde{\eta}(\rho,\rho u)|\leqq C_\gamma\big(|u|+|\rho^{\theta}-\bar{\rho}^{\theta}|\big),\\
&|\partial_{mu}\tilde{\eta}(\rho,\rho u)|\leqq C_\gamma, \,\,\,
|\partial_{m\rho}\tilde{\eta}(\rho,\rho u)|\leqq C_\gamma \rho^{\theta-1},
\end{split}
\end{equation}
which, with integration by parts, leads to
\begin{align}\label{6.94}
|J_4|
&\leqq\v\int_0^T\int_{b-1}^b\int_{d}^D\bigg|\int_r^y z^{N-1} \partial_m \tilde{\eta}\,\Big\{\big((\r+\alpha\delta\rho^{\alpha}) u_z\big)_z
  +(\r+\alpha\delta\rho^{\alpha}) (\frac{N-1}{z}u)_z\nonumber\\
&\qquad\qquad\qquad\qquad\qquad\qquad\quad\quad\,\,\,\,\,
   +(\alpha-1) \delta(\rho^{\alpha})_z\frac{N-1}{z} u \Big\}\,\dd z\bigg|\,\dd y\dd r\dd t\nonumber\\
&\leqq C\v\int_0^T\int_{b-1}^b\int_{d}^D \int_r^y \Big((\r+\delta\rho^{\alpha}) \big(\big|u_z (z^{N-1} \partial_m\tilde{\eta})_z\big|
                +\big|\partial_m\tilde{\eta}\,(\frac{u}{z})_z\big| z^{N-1}\big)\nonumber\\
 &\qquad\qquad\qquad\qquad\qquad\,\,\,\, +\delta \rho^{\alpha}\big|(z^{N-1} \partial_m\tilde{\eta})_z \frac{u}{z}\big|\Big)\,\dd z \dd r\dd y \dd t\nonumber\\
&\quad +C\v \int_0^T\int_{d}^D\Big( \big|\big(r^{N-1}(\r+\delta\rho^{\alpha})\partial_m \tilde{\eta} u_r\big)(t,r)\big|\nonumber\\
&\qquad\qquad\qquad\qquad  +\d\big|\big(r^{N-2}\rho^{\alpha} u \partial_m\tilde{\eta}\big)(t,r)\big|\Big)\, \dd r\dd t\nonumber\\
&\quad+CD\v\int_0^T\int_{b-1}^b\Big( \big|\big(y^{N-1}  (\r+\delta\rho^{\alpha})
  \partial_m \tilde{\eta} u_y\big)(t,y)\big|\nonumber\\
  &\qquad\qquad\qquad\qquad\quad +\d\big|\big(y^{N-2}\rho^{\alpha} u \partial_m\tilde{\eta}\big)(t,y)\big|\Big)\, \dd y\dd t.
\end{align}

In order to estimate the terms on the right-hand side of \eqref{6.94},  we notice that
\begin{equation}\label{2.134}
e(\rho,\bar{\rho})\,I_{B(t)}(r)\geqq C(\bar{\rho})^{-1},
\end{equation}
where $B(t)$ is defined in \eqref{2.76}.
Then combining \eqref{2.134} with  \eqref{6.10-1}, \eqref{6.11}, and \eqref{6.37} implies that
\begin{align}\label{6.95}
&\int_d^b \big(\rho^{\alpha} (\rho^{\theta}-\bar{\rho}^{\theta})^2\big)(t,r)\,r^{N-1}\dd r\nonumber\\
&\leqq C(\bar{\rho}) \int_d^b  I_{B^c(t)}(r) \big(\rho (\rho^{\theta}-\bar{\rho}^{\theta})^2\big)(t,r)\,r^{N-1}\dd r\nonumber\\
&\quad\, +\int_d^b I_{B(t)}(r) \big(\rho^{\alpha} (\rho^{\theta}-\bar{\rho}^{\theta})^2\big)(t,r)\,r^{N-1}\dd r \nonumber\\
&\leqq C(\bar{\rho})\int_d^b  e(\rho,\bar{\rho})(t,r)\,r^{N-1}\dd r +C(\bar{\rho}) \int_d^b I_{B(t)}(r)\, r^{N-1}\dd r\nonumber\\
&\leqq C(\bar{\rho})\int_d^b  e(\rho,\bar{\rho})(t,r)\,r^{N-1} \dd r
  +C(\bar{\rho}) \int_d^b  I_{B(t)}(r)\, e(\rho,\bar{\rho})(t,r)\, r^{N-1}\dd r\nonumber\\
&  \leqq C(\bar{\rho},E_0).
\end{align}
Combining \eqref{6.93} and \eqref{6.95} with \eqref{6.11}, \eqref{6.12}, and the Cauchy inequality,
we conclude that
the first term on the right-hand side of \eqref{6.94} are bounded by
\begin{align}\label{6.96}
& C(d)\int_0^T\int_{b-1}^b\int_{d}^D \int_r^y \Big\{\v z^{2}
   (\rho+\delta \rho^\alpha)(u_z^2+\rho^{\gamma-3}\rho_z^2)
  +\v \d  \rho^{\alpha} u^2     \nonumber\\
&\qquad\qquad\qquad\qquad\quad\quad\,\,
  +\rho^{\alpha}(\rho^{\theta}-\bar{\rho}^{\theta})^2 +z^{2}\big(\rho u^2+e(\rho,\bar{\rho})\big)\Big\}\,z^{N-3}\dd z\dd r\dd y\dd t\nonumber\\
&\leqq C(\bar{\rho},d,D,T,E_0).
\end{align}
Using \eqref{6.11}, \eqref{6.12}, \eqref{6.95}, and the Cauchy inequality,
we can bound the second term on the right-hand side of \eqref{6.94} by
\begin{align}\label{6.97}
& C(D,d)\int_0^T\int_{d}^D\Big\{\v(\rho+\delta \rho^\delta) u_r^2
   +\v \d \rho^{\alpha} \frac{u^2}{r^2}
+(\rho u^2+e(\rho,\bar{\rho})) \nonumber\\
&\qquad\qquad\qquad\quad\,\,\,\,\,
+\rho^{\alpha}(\rho^{\theta}-\bar{\rho}^{\theta})^2\Big\}\,r^{N-1}\dd r\dd t\nonumber\\
&\leqq C(\bar{\rho},d,D,T,E_0).
\end{align}
Using \eqref{6.11}, \eqref{6.12}, \eqref{6.16}, \eqref{6.34}, \eqref{6.95}, the Cauchy inequality,
and Lemma \ref{lem6.11}, the last term on the right-hand side of \eqref{6.94} can be bounded by
\begin{align}\label{6.98}
&\frac{CD}{d} \int_0^T\int_{b-1}^b\Big( \v (\r+\delta\rho^{\alpha}) |u_y|^2 +\big(\rho u^2+e(\rho,\bar{\rho})\big)\Big)\,y^{N-1}\dd y \dd t\nonumber\\[1.5mm]
&\quad+C(M_2, T)\int_0^T\int_{b-1}^b  |u(t,y)|^2\,y^{N-1} \dd y\dd t\nonumber\\
&\leqq C(\bar{\rho}, d,D,T,E_0)+C(T,M_2) b^{-\frac{\vartheta}{2}}.
\end{align}

\medskip
6. Substituting \eqref{6.89}, \eqref{6.87}, \eqref{6.92},  and \eqref{6.94}--\eqref{6.98}, we have
\begin{align}\label{6.83-1}
\int_0^T\int_d^Dr^{N-1}\tilde{q}(t,r)\,\dd r\dd t
\leqq C(\bar{\rho},d,D,T,E_0)+C(T,M_2) b^{-\frac{\vartheta}{2}}.
\end{align}
Then \eqref{6.15-1} follows
from \eqref{6.11}, \eqref{6.59}, and \eqref{6.83-1}.
This completes the proof.
$\hfill\Box$

\bigskip
Employing Proposition \ref{lem6.4}, we can obtain the following higher integrability estimate
up to the origin.

\begin{lemma}\label{lem6.15}
The smooth solution of \eqref{6.1}--\eqref{6.3} satisfies
\begin{align}\label{6.99}
\int_0^T\int_\d^1   \big(\rho|u|^3+\rho^{\g+\t}\big)(t,r)\,r^{N-1} \dd r\dd t
\leqq C(T,E_0)+C(T,M_2)\, b^{-\frac{\vartheta}{2}}.
\end{align}
\end{lemma}

\noindent{\bf Proof.} Let $w(r)$
be a smooth non-negative cut-off function with $\mbox{supp}\, w\subset [0,2]$  and $w(r)\equiv1$ for $r\in[0,1]$.
Multiplying $\eqref{6.1}_1$ by $w\partial_\rho \eta^{\#}(\rho,m)r^{N-1}$ and $\eqref{6.1}_2$ by $w\partial_m \eta^{\#}(\rho,m)r^{N-1}$, we have
\begin{align}\label{6.100}
&( w\eta^{\#}r^{N-1})_t+(wq^{\#}r^{N-1})_r-w_r q^{\#}r^{N-1}\nonumber\\
&\quad+(N-1)w \big(-q^{\#}+ m \partial_\rho\eta^{\#}+\frac{m^2}{\rho}\partial_m \eta^{\#}\big)r^{N-2}\nonumber\\
&=\v w \partial_m \eta^{\#}\,\Big\{\big((\r+\alpha \delta\r^{\alpha})(u_r+\frac{N-1}{r}u)\big)_r-\frac{N-1}{r}  (\rho+\d \rho^{\alpha})_r u \Big\}r^{N-1}.
\end{align}
Integrating \eqref{6.100} over $[r,2]$ with $r\leqq 2$, and then integrating the resultant equation over $[0,T]\times [\d,2]$ and using \eqref{6.63},
we have
\begin{align}\label{6.102}
&\int_0^T \int_{\d}^2 w(r)q^{\#}(t,r)\,r^{N-1} \dd r\nonumber\\
&\leqq  \left|\int_{\d}^2 \int_r^2  w(y)\eta^{\#}(T,y)\,y^{N-1} \dd y\dd r -\int_{\d}^2 \int_r^2 y^{N-1} w(y)\eta^{\#}(0,y)\,y^{N-1} \dd y\dd r\right|\nonumber\\
&\quad+\int_0^T\int_{\d}^2 \int_r^2 w_y(y) q^{\#}(t,y)\,y^{N-1}\dd y\dd r \dd t\nonumber\\
&\quad-\v\int_0^T\int_{\d}^2 \int_r^2 w(y) \partial_m \eta^{\#}\, \big((\r+\alpha \delta\r^{\alpha})u_y\big)_y\,y^{N-1}\dd y \dd r \dd t \nonumber\\
&\quad-(N-1)\v\int_0^T\int_{\d}^2 \int_r^2 w(y) \partial_m \eta^{\#}\, (\r+\alpha \delta\r^{\alpha})\,(\frac{u}{y})_y\,y^{N-1} \dd y \dd r \dd t \nonumber\\
&\quad -(N-1)(\alpha-1)\v\d\int_0^T\int_{\d}^2 \int_r^2 w(y) \partial_m \eta^{\#}\, (\r^{\alpha})_y\, \frac{u}{y}\,y^{N-1} \dd y \dd r \dd t\nonumber\\
&:=\sum_{j=1}^5 I_j.
\end{align}

For $I_1$, it follows from \eqref{6.59} and  Lemma \ref{lem6.1} that
\begin{align}\label{6.103}
I_1&\leqq C\int_{\d}^2\big(\rho|u|^2+\rho^\gamma\big)(T,y)\, y^{N-1} \dd y
 +\int_{\d}^2  \big(\rho_0|u_0|^2+\rho_0^\gamma\big)(y)\,y^{N-1}\dd y\nonumber\\
&\leqq C\int_{\d}^2  \big(1+\frac12 \rho |u|^2+e(\rho,\bar{\rho})\big)(t,y)\,y^{N-1}\dd y\nonumber\\
&\quad+C\int_{\d}^2  \big(1+\frac12 \rho_0 |u_0|^2+e(\rho_0,\bar{\rho})\big)(y)\,y^{N-1}\dd y\nonumber\\[1mm]
&\leqq C\big(E_0+1).
\end{align}

For $I_2$, we use Proposition \ref{lem6.4} with $d=1$ and $D=2$ to obtain
\begin{align}\label{6.104}
I_2\leqq C\int_0^T\int_1^2 q^{\#}(t,y)\,y^{N-1}\dd y \dd t
\leqq C(T,E_0)+C(T,M_2)\, b^{-\frac{\vartheta}{2}}.
\end{align}

For $I_3$, we integrate by parts to obtain
\begin{align}\label{6.105}
I_3&=(N-1)\v\int_0^T\int_{\d}^2\int_r^2  (\r+\alpha \delta\r^{\alpha})u_y\, \partial_m \eta^{\#} w(y)\,y^{N-2}\dd y \dd r \dd t\nonumber\\
&\quad +\v\int_0^T\int_{\d}^2 \int_r^2  (\r+\alpha \delta\r^{\alpha})u_y\,\partial_m \eta^{\#} w_y(y)\, y^{N-1} \dd y \dd r \dd t \nonumber\\
&\quad+\v\int_0^T\int_{\d}^2 \int_r^2  (\r+\alpha \delta\r^{\alpha})u_y\, (\partial_m \eta^{\#})_y w(y)\,y^{N-1} \dd y \dd r \dd t\nonumber\\
&\quad +\v\int_0^T\int_{\d}^2 (\r+\alpha \delta\r^{\alpha})u_r\,  \partial_m \eta^{\#}w(r)\,r^{N-1} \dd r \dd t\nonumber\\
&:=\sum_{j=1}^4 I_{3j}.
\end{align}
We regard $\eta_m^{\#}(\rho,\rho u)$ as a function of $(\rho,u)$ to see that
$$|\partial_{mu}\eta^{\#}(\rho,\rho u)| + \rho^{1-\theta}|\partial_{m\rho}\eta^{\#}(\rho,\rho u)|\leqq C_\gamma,$$
which, with \eqref{6.61} and Lemmas \ref{lem6.1}--\ref{lem6.2}, leads to
\begin{align}\label{6.106}
\sum_{j=2}^4I_{3j}
&\leqq C\int_0^T\int_\d^2 \v (\rho+\delta \r^{\alpha})\big(|u_y|^2 + \rho^{\gamma-3}\rho_y^2\big)\,y^{N-1}\dd y \dd t\nonumber\\
&\quad +\int_0^T\int_\d^2 \big(\v (\rho+\delta \r^{\alpha})|u|^2 +(\rho^{\gamma}+\r^{\alpha+\gamma-1})\big)\,y^{N-1}\dd y \dd t\nonumber\\[1.5mm]
&\leqq C(T,E_0).
\end{align}
To estimate $I_{31}$, we have to be more careful, since the weight is $y^{N-2}$  that may not be enough.
Fortunately, we can gain a weight $y$ by changing the order of integration:
\begin{align}\label{6.107}
I_{31}&=(N-1)\v\int_0^T \int_\d^2  (\r+\alpha \delta\r^{\alpha})u_y\, \partial_m \eta^{\#}w(y)\,(y-\delta)y^{N-2}\dd y\dd t\nonumber\\
&\leqq C\v\int_0^T \int_\d^2  (\r+\alpha \delta\r^{\alpha}) |u_y|(|u|+\rho^\theta )\, y^{N-1}\dd y\dd t\nonumber\\[1.5mm]
&\leqq  C(T,E_0).
\end{align}
Combining \eqref{6.105}--\eqref{6.106} with \eqref{6.107} yields
\begin{equation}\label{6.108}
I_3\leqq C(T,E_0).
\end{equation}

For $I_4$, using \eqref{6.61} and changing the order of integration as in \eqref{6.107}, we have
\begin{align}\label{6.109}
I_4
&\leqq C\v\int_0^T \int_{\d}^2 \int_r^2  (|u|+\rho^\theta)(\r+\alpha \delta\r^{\alpha})\, \big(|u_y|+\frac{|u|}{y}\big)\,y^{N-2}
    \dd y \dd r \dd t\nonumber\\
&\leqq C(T,E_0)+C\int_0^T\int_\d^2 \big(\v (\rho+\delta\rho^\alpha)\frac{|u|^2}{y^2}
  +(\rho^\gamma +\r^{\alpha+\gamma-1})\big)\,y^{N-1} \dd y \dd t\nonumber\\
&\leqq C(T,E_0).
\end{align}

Finally, for $I_5$, we first integrate by parts and then change the order of integration
as in \eqref{6.107} to obtain
\begin{align}\label{6.110}
I_5
&\leqq \, C(T,E_0)+C\v\d\int_0^T\int_{\d}^2 \rho^\alpha \big(|u_y|^2+\rho^{\gamma-3}|\rho_y|^2+ \frac{u^2}{y^2}+\rho^{\gamma-1}\big)
     \,y^{N-1}\dd y\dd t\nonumber\\
&\leqq\, C(T,E_0).
\end{align}

Substituting \eqref{6.103}--\eqref{6.104} and \eqref{6.108}--\eqref{6.110} into \eqref{6.102},
and using \eqref{6.59}, we conclude \eqref{6.99}.
$\hfill\Box$

\medskip
We now prove a lemma which is needed when we take the limit: $b\rightarrow \infty$.

\begin{lemma}\label{lem4.1}
The smooth solution of \eqref{6.1}--\eqref{6.3} satisfies that, for any $t\in [0,T]$,
\begin{align}\label{4.1}
&\|u_r(t)\|^2_{L^2}+\int_0^T\big(\|u_t(t)\|_{L^2}^2+\|u_{rr}(t)\|^2_{L^2}\big)\,\dd t
\leqq C(T,\|u_{0r}\|_{L^2},M_2).
\end{align}
\end{lemma}

\noindent{\bf Proof.}
It follows from $\eqref{6.1}_1$ that
\begin{align}\label{4.2}
-\v \big((\mu+\lambda)u_r\big)_r+\rho u_t=H,
\end{align}
where
$H:=-\rho u u_r-p_r+\v (\mu+\lambda)\big(\frac{N-1}{r} u\big)_r+\v \frac{N-1}{r} u \lambda_r$.
Multiplying \eqref{4.2} by $u_t$ and  integrating it over $[\d,b]$, we have
\begin{align}\label{4.4}
\frac{\v}{2} \frac{\dd}{\dd t}\int_\d^b (\mu+\lambda) |u_r|^2\, \dd r
  +\int_\d^b\rho u_t^2\, \dd r
  =\frac{\v}{2} \int_\d^b (\mu+\lambda)_t |u_r|^2\, \dd t+\int_\d^b H u_t\, \dd r.
\end{align}
Using \eqref{6.11}, \eqref{6.12}, \eqref{6.16}, \eqref{6.34}, and the Sobolev inequality:
$$
\|u_r\|_{L^\infty}\leqq C\big(\|u_r\|_{L^2}+\|u_r\|_{L^2}^{\frac12}\,\|u_{rr}\|_{L^2}^{\frac12}\big),
$$
we obtain
\begin{align}
&\frac{\v}{2} \int_\d^b (\mu+\lambda)_t |u_r|^2\, \dd r\nonumber\\
&\leqq C(T,M_2) \int_\d^b \big(|\rho_r u|+|u_r|+|u|\big) |u_r|^2\, \dd r\nonumber\\
&\leqq C(T,M_2)\Big\{\|u\|_{L^2}^{\frac12}\|\rho_r\|_{L^2}\,\big(\|u_r\|_{L^2}^2\|u_{rr}\|^{\frac12}_{L^2}+\|u_{r}\|^{\frac52}_{L^2}\big)
\nonumber\\
&\qquad\qquad\quad \,\,\,\,\,+\|u\|_{L^2}^{\frac12}\|u_{r}\|^{\frac52}_{L^2} +\|u_{r}\|^{2}_{L^2}\,\big(\|u_{r}\|_{L^2}+\|u_r\|_{L^2}^{\frac12}\|u_{rr}\|^{\frac12}_{L^2}\big)\Big\}\nonumber\\
&\leqq C(T,M_2)\Big\{\big(\|u_r\|_{L^2}^2+\|u_r\|_{L^2}^{\frac52}\big)\|u_{rr}\|^{\frac12}_{L^2}+\|u_{r}\|^{3}_{L^2}+1\Big\},\label{4.5}\\[1mm]
&\big|\int_\d^b H u_t\, \dd r\big|\nonumber\\
&\leqq \frac{1}8 \int_\d^b \rho |u_t|^2\, \dd r+C\int_\d^b\rho^{-1}|H|^2\, \dd r\nonumber\\
&\leqq C(T,M_2)\Big\{ \big(\|u\|^2_{L^\infty}+1\big)\|(\r_r, u_r)\|_{L^2}^2+\|u\|_{L^2}^2\Big\}
  +\frac{1}8 \int_\d^b \rho |u_t|^2\, \dd r\nonumber\\
&\leqq \frac{1}8 \int_\delta^b \rho |u_t|^2 \dd r
  +C(T,M_2)\big(\|u_r\|^3_{L^2}+1\big). \label{4.6}
\end{align}
To close the above estimate, we combine \eqref{4.2} with \eqref{6.11}, \eqref{6.12},  \eqref{6.16}, and \eqref{6.34} to obtain
\begin{align}
\|u_{rr}\|^2_{L^2}
&\leqq C(T,M_2)
  \Big\{ \|\sqrt{\rho} u_t\|_{L^2}^2+ \|\rho_r\|^2_{L^2} \|u_r\|_{L^2} \|u_{rr}\|_{L^2} +\|H\|_{L^2}^2 \Big\}\nonumber\\
&\leqq C(T,M_2)\Big\{ \|\sqrt{\rho} u_t\|_{L^2}^2+ \|u_r\|_{L^2} \|u_{rr}\|_{L^2}+\|u_r\|_{L^2}^3+1 \Big\}\nonumber\\
&\leqq C(T,M_2)\Big\{ \|\sqrt{\rho} u_t\|_{L^2}^2+\|u_r\|_{L^2}^3+1 \Big\}.\label{4.7}
\end{align}
Combining \eqref{4.4}--\eqref{4.7}, we obtain
\begin{align}
\frac{\dd}{\dd t}\int_\d^b (\mu+\lambda) |u_r|^2\, \dd r+\int_\d^b\rho u_t^2\, \dd r
\leqq C(T,M_2)
 \Big\{ 1+\|u_r\|_{L^2}^2 \int_\d^b(\mu+\lambda) |u_r|^2\,\dd r\Big\}.\nonumber
\end{align}
Applying the Gronwall inequality, we have
\begin{align}
\int_\delta^b (\mu+\lambda) |u_r|^2\, \dd r+\int_0^t\int_\d^b\rho u_t^2\, \dd r\dd s
\leqq C(T,\|u_{0r}\|_{L^2},M_2),\nonumber
\end{align}
which,  with \eqref{4.7}, implies \eqref{4.1}.
$\hfill\Box$

\section{Limits of the Approximate Solutions for the Navier-Stokes Equations}
In this section, we first take the limit, $b\rightarrow\infty$, to obtain
global strong solutions $(\rho^{\v,\d}, u^{\v,\delta})$ of the Navier-Stokes equations with some uniform bounds
with respect to $(\v, \d)$.
Then we take the limit, $\delta\rightarrow0+$, to obtain global, spherically symmetric weak solutions of the
Navier-Stokes equations \eqref{1.1} with some desired uniform bounds with respect to $\v>0$ on $[0,T]\times [0, \infty)$,
which are essential for us to employ the compensated compactness framework in \S 6.

\subsection{Passage the limit: $b\rightarrow \infty$}
$\,\,$ In this subsection,
we fix parameters $(\v, \d)$ and denote the solution of \eqref{6.1}--\eqref{6.3}
as $(\rho^{\v,\d,b}, u^{\v,\d,b})$.
It follows from \eqref{data-b}--\eqref{data-d} and
Lemmas \ref{lem8.2}--\ref{lem8.4} in the appendix
that there exist sequences of smooth approximate initial data functions $(\rho_0^{\v,\d,b}, u_0^{\v,\d,b})$
and $(\rho_0^{\v,\d}, u_0^{\v,\d})$ satisfying \eqref{6.3a} and the properties:
\begin{align}\label{5.6}
\begin{cases}
\displaystyle(\rho^{\v,\d,b}_0, m^{\v,\d,b}_0)(r)\rightarrow (\rho_0^{\v,\delta}, m_0^{\v,\d})(r) \,
  &\mbox{in $L^1_{\rm loc}([\delta, \infty); r^{N-1}{\rm d}r)\, $
  as $b\rightarrow \infty$},\\[2mm]
\displaystyle (E_0^{\v,\d,b},  E_1^{\v,\d,b})\rightarrow (E_0^{\v,\d}, E_1^{\v,\d})\, &\mbox{as $b\rightarrow \infty$},\\[2mm]
\displaystyle E_2^{\v,\d,b}+ \tilde{E}_0^{\v,\d,b}+ \|u^{\v,\d,b}_{0r}\|_{L^2} \,   &\mbox{is uniform bounded in $b$},
\end{cases}
\end{align}
where
\begin{align}
E_0^{\v,\d}:&=\int_\d^\infty\bar{\eta}^{\ast}(\rho_0^{\v,\d},m_0^{\v,\d})\,r^{N-1} \dd r<\infty,\label{1.53-1a}\\
E_1^{\v,\d}:&= \v^2\int_\d^\infty \big(1+2\alpha\d (\rho_0^{\v,\delta})^{\alpha-1}+\alpha^2\d^2 (\rho_0^{\v,\delta})^{2\alpha-2}\big)
  \big|\big(\sqrt{\rho_{0}^{\v,\delta}}\big)_r\big|^2\, r^{N-1} \dd r <\infty.\label{1.53-1}
\end{align}

\smallskip
From \eqref{6.11}, \eqref{6.12}, \eqref{6.16}, \eqref{6.34}, and \eqref{4.1},
there exists a positive constant $\tilde{C}>0$ that may depend on $(\v, \delta, T)$, but is independent of $b$, so that
\begin{align}
0<\tilde{C}^{-1}\leqq \rho^{\v,\d,b}(t,r) &\leqq \tilde{C}, \label{5.1-1}\\[1.5mm]
\sup_{t\in[0,T]}\Big(\big\|(\rho^{\v,\d,b}-\bar{\rho},  u^{\v,\d,b})\big\|^2_{H^1([\d,b])}
  &+\big\|\rho^{\v,\d,b}_t\big\|^2_{L^2([\d,b])}\Big)(t)\nonumber\\
 +\int_0^T\big\|(u^{\v,\d,b}_t, u^{\v,\d,b}_{rr})\big\|^2_{L^2([\d,b])}(t)\,\dd t &\leqq \tilde{C}.\label{5.1}
\end{align}
We extend $\rho^{\v,\d,b}(t,r)$ and $u^{\v,\d,b}(t,r)$ to $[0,T]\times[\d,\infty)$
by defining $\rho^{\v,\d,b}(t,r)=\bar{\rho}$ and $u^{\v,\d,b}(t,r)=0$ for all $r\in [0,T]\times(b,\infty)$.
Then it follows from \eqref{5.1} and  the  Aubin-Lions lemma that
\begin{align}\nonumber
(\rho^{\v,\d,b}, u^{\v,\d,b}) \qquad \mbox{is compact in $C([0,T]; L^p_{\rm loc}[\d,\infty))$ with $ p\in[1,\infty)$}.
\end{align}

More precisely, we have
\begin{lemma}\label{lem5.1}
There exist functions $(\rho^{\v,\delta},u^{\v,\delta})(t,r)$ so that,  as $b\rightarrow\infty$ $($up to a subsequence$)$,
\begin{align}\nonumber
(\rho^{\v,\d,b}, u^{\v,\d,b})\rightarrow (\rho^{\v,\delta},u^{\v,\delta})  \quad \mbox{strongly in $C([0,T]; L^p_{\rm loc}[\d,\infty))$ for all $p\in[1,\infty)$.}
\end{align}
In particular,  as $b\to\infty$ $($up to a subsequence$)$,
\begin{align}\nonumber
(\rho^{\v,\d,b}, u^{\v,\d,b})\rightarrow (\rho^{\v,\delta},u^{\v,\delta})  \qquad \mbox{a.e.\, $(t,r)\in  [0,T]\times[\d,\infty)$}.
\end{align}
\end{lemma}

Using Lemma \ref{lem5.1}, it can immediately be proved that $(\rho^{\v,\delta},u^{\v,\delta})$ is a weak solution of the initial-boundary value problem (IBVP)
of the Navier-Stokes equations \eqref{6.1}:
\begin{align}\label{5.3}
\begin{cases}
(\rho, u)(0,r)=(\rho_0^{\v,\delta}, u_0^{\v,\delta})(r) \qquad \mbox{for $r\in [\d, \infty)$},\\[1mm]
u|_{r=\d}=0 \qquad\mbox{for $t\geqq 0$}.
\end{cases}
\end{align}
Moreover,  it follows from \eqref{5.1-1}--\eqref{5.1} and the lower semicontinuity that
\begin{align}
&0<\tilde{C}^{-1}\leqq \rho^{\v,\delta}(t,r) \leqq \tilde{C},\label{5.5}\\[1.5mm]
&\sup_{t\in[0,T]}\Big(\big\|(\rho^{\v,\delta}-\bar{\rho},  u^{\v,\delta})\big\|^2_{H^1([\d,\infty))}
 + \|\rho_t^{\v,\d}\|^2_{L^2([\d,\infty))}\Big)(t)\nonumber\\
&\qquad +\int_0^T\big\|(u_t^{\v,\delta}, u_{rr}^{\v,\delta})\big\|^2_{L^2([\d,\infty))}(t)\,\dd t\leqq \tilde{C}.\label{5.4}
\end{align}
These facts yield that the weak solution $(\rho^{\v,\delta}, u^{\v,\delta})$ of \eqref{5.3}
is indeed a strong solution. The uniqueness of this strong solution $(\rho^{\v,\delta},u^{\v,\delta})$
is ensured by properties \eqref{5.5}--\eqref{5.4},
the corresponding version of Lemmas 3.1--3.2 ({\it i.e.}, \eqref{5.7}--\eqref{5.8} below),
and the basic $L^2$--energy estimate as in \S 3.
This implies that the whole sequence $(\rho^{\v,\d,b},u^{\v,\d,b})$
converges to $(\rho^{\v,\delta},u^{\v,\delta})$ as $b\to\infty$.

\vspace{1.5mm}
Then it is direct to know that
$(\rho^{\v,\d}, \mathcal{M}^{\v,\d})(t,\textbf{x})=(\rho^{\v,\d}(t,r),m^{\v,\d}(t,r)\, \frac{\textbf{x}}{r})$
with $\rho^{\v,\d}(t,\textbf{x})>0$
is a strong solution of the initial-boundary problem of system \eqref{1.1} with $(h,g)$ determined by \eqref{6.1-2}
for  $(t,\textbf{x})\in [0,\infty)\times \big(\mathbb{R}^{N}\backslash B_{\d}(\mathbf{0})\big)$
with the following initial-boundary data:
\begin{align}\label{5.5-2}
\begin{cases}
(\rho^{\v,\d}, \mathcal{M}^{\v,\d})(0,\textbf{x})=(\rho^{\v,\d}_0(r), m^{\v,\d}_0(r)\,\frac{\textbf{x}}{r}),\\[1.5mm]
\mathcal{M}^{\v,\d}(t,\textbf{x})|_{\textbf{x}\in \partial B_\d(\mathbf{0})}=0.
\end{cases}
\end{align}

\medskip
From Lemma \ref{lem5.1},  \eqref{6.11}, \eqref{6.12}, \eqref{6.14}, \eqref{6.15-1}--\eqref{6.16}, \eqref{6.99}, \eqref{5.6},
Fatou's lemma, and the lower semicontinuity, we have

\begin{proposition}\label{prop5.1}
Under assumption \eqref{5.6}, for any fixed $(\v, \delta)$,  there exists a unique strong solution $(\rho^{\v,\delta}, u^{\v,\delta})$
of IBVP \eqref{5.3}. Moreover, $(\rho^{\v,\delta}, u^{\v,\delta})$ satisfies \eqref{5.5} and, for $t\in (0, T]$,
\begin{align}
&\int_\d^\infty  \Big(\frac12\rho^{\v,\delta} |u^{\v,\delta}|^2+e(\rho^{\v,\delta},\bar{\r}) \Big)(t,r)\,r^{N-1}\dd r\nonumber\\
&\quad+\v\int_0^T\int_\d^\infty \Big(\rho^{\v,\delta}|u_r^{\v,\delta}|^2+\rho^{\v,\delta}\frac{|u^{\v,\delta}|^2}{r^2}\Big)(s,r)\,r^{N-1}\dd r \dd s\nonumber\\
&\quad +c_N\v\delta \int_0^T\int_\d^\infty  \Big((\rho^{\v,\delta})^{\alpha}\big(|u_r^{\v,\delta}|^2+ \frac{|u^{\v,\delta}|^2}{r^{2}}\big)\Big)(s,r)\,r^{N-1}\dd r \dd s\nonumber\\[1mm]
&\leqq E_0^{\v,\delta}\leqq C(E_0+1),
  \label{5.7}\\[3mm]
&\v^2 \int_\d^\infty\Big(\big|(\sqrt{\rho^{\v,\d}})_r\big|^2 +\d (\rho^{\v,\delta})^{\alpha-2}|\rho_r^{\v,\delta}|^2
  +\d^2 (\rho^{\v,\delta})^{2\alpha-3}|\rho_r^{\v,\delta}|^2\Big)(t,r)\, r^{N-1}\dd r\nonumber\\
&\quad+\v\int_0^T\int_\d^\infty \left(\big|\big((\rho^{\v,\d})^{\frac{\gamma}{2}}\big)_r\big|^2
  +\d (\rho^{\v,\delta})^{\g+\alpha-3} |\rho_r^{\v,\delta}|^2 \right)(s,r)\, r^{N-1}\dd r\dd s\nonumber\\[1mm]
  & \leqq C (E_0+1),
\label{5.8}\\[2mm]
&\int_0^T\int_{d}^D (\rho^{\v,\delta})^{\g+1}(t,r)\,\dd r \dd t\leqq C(d,D,T,E_0), \label{5.9}\\[3mm]
&\int_0^T\int_\d^D  \big(\rho^{\v,\delta}|u^{\v,\delta}|^3+(\rho^{\v,\delta})^{\gamma+\theta}\big)(t,r)\,r^{N-1}\dd r\dd t
 \leqq C(D,T,E_0)  \label{5.10}
\end{align}
for any fixed $T>0$ and any compact subset $[d,D]$ of $(\d,\infty)$, where $c_N>0$ is some constant
depending only on $N$ determined in Lemma {\rm \ref{lem6.1}}.
\end{proposition}

\subsection{Passage the limit: $\delta\rightarrow0+$}
In this subsection, for fixed $\v>0$, we consider the limit, $\delta\rightarrow0+$,
to obtain the weak solution of the Navier-Stokes equations.
It follows from Lemma \ref{lem8.4} in the appendix that
\begin{align}\label{5.11}
\begin{cases}
(\rho_0^{\v,\delta}, m_0^{\v,\d})(r)\rightarrow (\rho_0^\v, m_0^\v)(r) \quad &\mbox{in $L^1_{\rm loc}([0,\infty); r^{N-1}\dd r)\,$
   as $\delta\rightarrow 0+$},\\[1.5mm]
(E_0^{\v,\d}, E_1^{\v,\d}) \rightarrow (E_0^\v, E_1^\v) \qquad &\mbox{as $\d\rightarrow0$}.
\end{cases}
\end{align}
To take the limit, we have to be careful since the weak solution may involve the vacuum state.
We use similar compactness arguments as in \cite{MV,Guo-Jiu-Xin-2} to
consider the limit: $\delta\rightarrow0+$.
We first extend our solution $(\rho^{\v,\delta}, u^{\v,\delta})$ as the zero extension of $(\rho^{\v,\delta}, u^{\v,\delta})$
outside  $[0,T]\times[\d,\infty)$.

\begin{lemma}\label{lem5.2}
There exists a function $\rho^{\v}(t,r)$ such that, as $\delta\rightarrow0+$ $($up to a subsequence$)$,
\begin{align}
(\rho^{\v,\delta},\sqrt{\rho^{\v,\delta}}) \rightarrow (\rho^{\v},\sqrt{\rho^{\v}})
\,\, \mbox{a.e. and strongly in $C(0,T; L^{q}_{\rm loc})$ for any $q\in [1,\infty)$},\label{5.12}
\end{align}
where $L^q_{\rm loc}$ means $L^q(K)$ for any $K\Subset (0,\infty)$.
\end{lemma}

\noindent{\bf Proof.}
It follows from \eqref{5.7}--\eqref{5.8} that
$$
\sqrt{\rho^{\v,\delta}}\in L^\infty(0,T; H^{1}_{\rm loc})\hookrightarrow L^\infty(0,T; L^{\infty}_{\rm loc}) \quad \mbox{uniformly}.
$$
Notice that, for fixed $\v>0$,
the solution sequence $(\rho^{\v,\d}, u^{\v,\d})$ satisfies  $\eqref{6.1}$ for $(t,r)\in [0,\infty)\times [\delta,\infty)$.
Using \eqref{5.7} and the mass equation $\eqref{6.1}_1$, we see that
\begin{align}
\partial_t\sqrt{\rho^{\v,\d}}
=-(\sqrt{\rho^{\v,\d}} u^{\v,\d})_r+\frac12 \sqrt{\rho^{\v,\d}} u_r^{\v,\d}-\frac{N-1}{2r}\sqrt{\rho^{\v,\d}} u^{\v,\d}
\end{align}
is uniformly bounded in $L^{2}(0,T; H^{-1}_{\rm loc})$,
which, using the Aubin-Lions lemma, implies that
\begin{align}
\sqrt{\rho^{\v,\delta}} \qquad \mbox{is compact in $C(0,T; L^q_{\rm loc})$ for any $q\in [1,\infty)$}.\nonumber
\end{align}

Since $\sqrt{\rho^{\v,\delta}}$ and $\sqrt{\rho^{\v,\delta}} u^{\v,\d}$
are uniformly  bounded in $L^\infty(0,T; L^{\infty}_{\rm loc})$ and $L^\infty(0,T; L^{2}_{\rm loc})$
respectively,
we see that
\begin{align}\label{5.18}
\rho^{\v,\delta} u^{\v,\d}=\sqrt{\rho^{\v,\delta}}\big(\sqrt{\rho^{\v,\delta}} u^{\v,\d}\big)
\quad \mbox{is uniformly bounded in $L^\infty(0,T; L^{2}_{\rm loc})$}.
\end{align}
Then it follows from the mass equation $\eqref{6.1}_1$ that
\begin{align}\nonumber
\partial_t \rho^{\v,\d}=-(\rho^{\v,\d} u^{\v,\d})_r-\frac{N-1}{r} \rho^{\v,\d} u^{\v,\d}
\,\,\,\,\,\mbox{is uniformly bounded in $L^{\infty}(0,T; H^{-1}_{\rm loc})$}.
\end{align}
Moreover, we obtain that
\begin{align}\nonumber
\rho_r^{\v,\d}=2 \sqrt{\rho^{\v,\d}} (\sqrt{\rho^{\v,\d}})_r \quad \mbox{is uniformly bounded in $L^{\infty}(0,T; L^{2}_{\rm loc})$}.
\end{align}
Then the Aubin-Lions lemma implies that
\begin{align*}
\rho^{\v,\d}\qquad  \mbox{is compact in $C(0,T; L^{q}_{\rm loc})$ with $q\in [1,\infty)$}.
\end{align*}
$\hfill\Box$

\medskip
\begin{corollary}\label{cor5.3}
The pressure function sequence $p(\rho^{\v,\d})$ is uniformly bounded in $L^\infty(0,T; L^q_{\rm loc})$
for all $q\in[1,\infty]$ and, as $\delta\to 0+$ $($up to a subsequence$)$,
\begin{align}\label{5.16}
p(\rho^{\v,\d})\rightarrow p(\rho^{\v}) \qquad  \mbox{strongly in $L^q(0,T; L^q_{\rm loc})$ for all $q\in[1,\infty)$}.
\end{align}
\end{corollary}

\vspace{1mm}
\begin{lemma}\label{lem5.5}
As $\delta\rightarrow 0+$ $($up to a subsequence$)$,
$m^{\v,\delta}$
converges strongly in $L^2(0,T; L^q_{\rm loc})$
to some function $m^{\v}(t,r)$ for all $q\in[1,\infty)$, which implies that
\begin{align}
m^{\v,\d}(t,r)=(\rho^{\v,\delta} u^{\v,\delta})(t,r)\rightarrow m^{\v}(t,r) \qquad  \mbox{a.e. in $[0,T]\times(0,\infty)$}. \nonumber
\end{align}
\end{lemma}

\noindent{\bf Proof.}   A direct calculation shows
that
\begin{align}\label{5.19}
m^{\v,\d}_r
=2\big(\sqrt{\rho^{\v,\d}}\big)_r\, \big(\sqrt{\rho^{\v,\d}}u^{\v,\d}\big)
 +\sqrt{\rho^{\v,\d}}\, \big(\sqrt{\rho^{\v,\d}} u_r^{\v,\d}\big)
\end{align}
is uniformly bounded in $L^{2}(0,T; L^1_{\rm loc})$.
Thus, it follows from  \eqref{5.18}--\eqref{5.19} that
\begin{align}\label{5.20}
m^{\v,\d} \qquad  \mbox{is uniformly bounded in $L^{2}(0,T; W^{1,1}_{\rm loc})$}.
\end{align}

It follows from \eqref{5.7} and \eqref{5.16} that
$\partial_r\big((\sqrt{\rho^{\v,\d}} u^{\v,\d})^2\big)$,
$\frac{N-1}{r} \big(\sqrt{\rho^{\v,\d}} u^{\v,\d}\big)^2$,
and $\partial_r p(\rho^{\v,\d})$
are uniformly bounded in $L^{\infty}(0,T;W_{\rm loc}^{-1,1})$,
$L^{\infty}(0,T; L_{\rm loc}^{1})$, and $L^{2}(0,T;H_{\rm loc}^{-1})$,
respectively.

From \eqref{5.7}, we see that
$$
\sqrt{(\rho^{\v,\d})^\alpha}\big(\sqrt{\d\,(\rho^{\v,\d})^\alpha}(u^{\v,\d}_r+\frac{N-1}{r} u^{\v,\d})\big)
\,\,\,\,\mbox{and}\,\,\,\,
\sqrt{\rho^{\v,\d}}\big(\sqrt{\rho^{\v,\d}}(u^{\v,\d}_r+\frac{N-1}{r} u^{\v,\d})\big)
$$
are uniformly bounded in $L^2(0,T; L^2_{\rm loc})$.

Since
\begin{align*}
&\big(\mu(\rho^{\v,\d})+\lambda(\rho^{\v,\d})\big)\big(u_r^{\v,\d}+\frac{N-1}{r} u^{\v,\d}\big)\\
&=\Big(\alpha\d\sqrt{(\rho^{\v,\d})^\alpha}+  \sqrt{(\rho^{\v,\d})^{2-\alpha}}\Big) \,\Big(\sqrt{(\rho^{\v,\d})^\alpha}u^{\v,\d}_r
 +\frac{N-1}{r} \sqrt{(\rho^{\v,\d})^\alpha} u^{\v,\d}\Big),
\end{align*}
we conclude that
$$
\partial_r\Big(\big(\mu(\rho^{\v,\d})+\lambda(\rho^{\v,\d})\big)\big(u_r^{\v,\d}+\frac{N-1}{r} u^{\v,\d}\big)\Big)
$$
is uniformly bounded in $L^{2}(0,T; H^{-1}_{\rm loc})$.
Also, it follows from \eqref{5.7}--\eqref{5.8} that
\begin{align*}
\frac{N-1}{r}\partial_r\mu(\rho^{\v,\d}) u^{\v,\d}
=\frac{2(N-1)}{r}\big((\sqrt{\rho^{\v,\d}})_r
+\alpha \delta(\rho^{\v,\d})^{\alpha-\frac32}\rho_r^{\v,\d}\big)\big(\sqrt{\rho^{\v,\d}}u^{\v,\d}\big)
\end{align*}
is uniformly bounded in $L^2(0,T; L^1_{\rm loc})$.
Then
we conclude that
\begin{align}
\partial_t m^{\v,\d}  \qquad\mbox{is uniformly bounded in $L^{2}(0,T; W^{-2,\frac43}_{\rm loc})$},\nonumber
\end{align}
which,  with \eqref{5.20} and the Aubin-Lions lemma, implies that
\begin{align}
m^{\v,\d} \qquad \mbox{is compact in $L^{2}(0,T; L^p_{\rm loc})$ for all $p\in[1,\infty)$}.\nonumber
\end{align}
$\hfill\Box$

\medskip
\begin{lemma}\label{lem5.6}
$m^{\v}(t,r)=0$ a.e. on $\{(t,r)\, :\, \rho^{\v}(t,r)=0\}$. Furthermore,
there exists a function $u^{\v}(t,r)$ so that $m^{\v}(t,r)=\rho^{\v}(t,r) u^{\v}(t,r)$ a.e.,
$\, u^\v(t,r)=0$ a.e. on $\{(t,r)\, :\, \rho^{\v}(t,r)=0\}$,
and
\begin{align}
& m^{\v,\d} \rightarrow m^{\v} &&\mbox{strongly in $L^2(0,T; L^p_{\rm loc})$ for $p\in[1,\infty)$},\nonumber\\
&  \frac{m^{\v,\d}}{\sqrt{\rho^{\v,\d}}}  \rightarrow \frac{m^\v}{\sqrt{\rho^\v}}=\sqrt{\rho^\v}u^\v &&\mbox{strongly in $L^2(0,T; L^2_{\rm loc})$}.\nonumber
\end{align}
\end{lemma}

\noindent{\bf Proof.}
Since
$\frac{m^{\v,\d}}{\sqrt{\rho^{\v,\d}}}r^{\frac{N-1}{2}}$
is uniformly bounded in $L^{\infty}(0,T;L^2)$, then Fatou's lemma implies
\begin{align}
&\int_0^T\int_0^\infty \liminf\limits_{\delta\to 0+} \frac{|m^{\v,\d}(t,r)|^2}{\rho^{\v,\d}(t,r)}\,r^{N-1}\dd r\dd t\nonumber\\
&\leqq \lim_{\d\rightarrow0+}\int_0^T\int_0^\infty \frac{|m^{\v,\d}(t,r)|^2}{\rho^{\v,\d}(t,r)}\,r^{N-1}\dd r\dd t<\infty.\nonumber
\end{align}
Thus, $m^{\v}(t,r)=0$ {\it a.e.} on $\{(t,r) \,:\, \rho^{\v}(t,r)=0\}$.
Then, if the limit velocity $u^{\v}(t,r)$ is defined by setting
$u^{\v}(t,r):=\frac{m^{\v}(t,r)}{\rho^{\v}(t,r)}$ {\it a.e.} on $\{(t,r) \,:\, \rho^{\v}(t,r)\neq 0\}$
and $u^{\v}(t,r)=0$ {\it a.e.} on $\{(t,r) \,:\, \rho^{\v}(t,r)=0\}$,
we have
\begin{align}
&m^\v(t,r)=\rho^{\v}(t,r) u^{\v}(t,r) \,\,\,\, a.e.,\nonumber\\[1mm]
&\int_0^T\int_0^\infty\big|\frac{m^\v}{\sqrt{\rho^\v}}\big|^2 r^{N-1}\,\dd r\dd t
=\int_0^T\int_0^\infty\rho^\v |u^\v|^2 r^{N-1}\dd r\dd t <\infty.
\nonumber
\end{align}
Moreover, it follows from \eqref{5.10} and Fatou's lemma that, for $[d,D]\Subset(0, \infty)$,
\begin{align}\label{5.28}
\int_0^T\int_{d}^{D} \rho^\v |u^\v|^3\,\dd r\dd t
\leqq \lim_{\d\rightarrow0+}\int_0^T\int_{d}^{D} \frac{|m^{\v,\d}|^3}{(\rho^{\v,\d})^2}\,\dd r\dd t\leqq C(d,D,T,E_0)<\infty.
\end{align}

Next, since $m^{\v,\d}$ and $\rho^{\v,\d}$ converge {\it a.e.},
it is direct to know that sequence $\sqrt{\rho^{\v,\d}} u^{\v,\d}=\frac{m^{\v,\d}}{\sqrt{\rho^{\v,\d}}}$ converges
{\it a.e.} to $\sqrt{\rho^{\v}} u^{\v}=\frac{m^{\v}}{\sqrt{\rho^{\v}}}$ on $\{(t,r)\,:\, \rho^{\v}(t,r)\neq0\}$.
Moreover, for any given positive constant $R>0$, it follows from Lemmas 5.3 and 5.6 that
\begin{align}\label{5.29}
\sqrt{\rho^{\v,\d}} u^{\v,\d} I_{|u^{\v,\d}|\leqq R} \rightarrow \sqrt{\rho^{\v}} u^{\v} I_{|u^{\v}|\leqq R} \qquad a.e.
\end{align}

For $R\geqq1$, we cut the $L^2$--norm as follows:
\begin{align}\label{5.30}
&\int_0^T\int_d^D \big|\sqrt{\rho^{\v,\d}} u^{\v,\d}-\sqrt{\rho^\v} u^{\v}\big|^2\,\dd r\dd t\nonumber\\
&\leqq \int_0^T\int_d^D \big|\sqrt{\rho^{\v,\d}} u^{\v,\d} I_{|u^{\v,\d}|\leqq R}
          - \sqrt{\rho^{\v}} u^{\v} I_{|u^{\v}|\leqq R}\big|^2\,\dd r\dd t\nonumber\\
&\quad + 2\int_0^T\int_d^D \big|\sqrt{\rho^{\v,\d}} u^{\v,\d} I_{|u^{\v,\d}|\geqq R}\big|^2\,\dd r\dd t\nonumber\\
&\quad  +2\int_0^T\int_d^D\big|\sqrt{\rho^{\v}} u^{\v} I_{|u^{\v}|\geqq R}\big|^2\, \dd r\dd t.
\end{align}
It is direct to know that $\sqrt{\rho^{\v,\d}} u^{\v,\d} I_{|u^{\v,\d}|\leqq R}$ is uniformly bounded
in $L^\infty(0,T; L^p_{\rm loc})$ for all $p\in[1,\infty)$.   Then it follows from \eqref{5.29} that
\begin{align}\label{5.31}
\int_0^T\int_d^D \big|\sqrt{\rho^{\v,\d}} u^{\v,\d} I_{|u^{\v,\d}|\leqq R}
- \sqrt{\rho^{\v}} u^{\v} I_{|u^{\v}|\leqq R}\big|^2\,\dd r\dd t\rightarrow 0 \qquad\,\,
\mbox{as $\delta\rightarrow0+$}.
\end{align}
Using \eqref{5.28}, we have
\begin{align}
&\int_0^T\int_d^D \big(\big|\sqrt{\rho^{\v,\d}} u^{\v,\d} I_{|u^{\v,\d}|\geqq R}\big|^2
    +\big|\sqrt{\rho^{\v}} u^{\v} I_{|u^{\v}|\geqq R}\big|^2\big)\,\dd r\dd t\nonumber\\
&\leqq \frac{1}{R}\int_0^T\int_d^D \big(\rho^{\v,\d} |u^{\v,\d} |^3
     +\rho^{\v} |u^{\v} |^3\big)\,\dd r\dd t
\leqq C(d,D,T,E_0)R^{-1}.\label{5.33}
\end{align}

Substituting \eqref{5.31}--\eqref{5.33} into \eqref{5.30} leads to
\begin{align}\nonumber
\lim_{\delta\rightarrow0+} \int_0^T\int_d^D \big|\sqrt{\rho^{\v,\d}} u^{\v,\d}-\sqrt{\rho^\v} u^{\v}\big|^2\,\dd r\dd t
\leqq C(d,D,T,E_0)R^{-1} \quad \mbox{for all $R>0$}.
\end{align}
Then the lemma follows by taking $R\rightarrow \infty$.
$\hfill\Box$

\medskip
Let $(\rho^\v, m^\v)$ be the limit obtained above. By using Fatou's lemma and
the lower semicontinuity and Proposition \ref{prop5.1}, it is direct to obtain

\begin{proposition}\label{prop5.2}
Under assumption \eqref{5.11}, for any fixed $\v$ and  $T>0$, the limit functions $(\rho^\v, m^\v)=(\rho^\v, \rho^\v u^\v)$
satisfy
\begin{align}
&\rho^{\v}(t,r)\geqq0 \,\,\,\, a.e.,\\
&u^\v(t,r)=0, \,\,\big(\frac{m^\v}{\sqrt{\rho^\v}}\big)(t,r)=\sqrt{\rho^\v}(t,r)u^\v(t,r)=0 \,\, a.e.\,\, \mbox{on $\{(t,r)\,:\, \rho^\v(t,r)=0\}$}, \label{5.7-0}\\[1mm]
&\int_0^\infty \Big(\frac12\Big|\frac{m^{\v}}{\sqrt{\rho^\v}}\Big|^2+e(\rho^{\v},\bar{\r}) \Big)(t,r)\, r^{N-1} \dd r
  +\v\int_{\mathbb{R}_+^2} \Big|\frac{m^{\v}}{\sqrt{\rho^\v}}\Big|^2(s,r)\,r^{N-3}\dd r \dd s\nonumber\\[1mm]
  &\quad\leqq E_0^\v\leqq E_0+1\quad\,\, \mbox{for $t\geqq0$}, \label{5.7-1}\\[3mm]
&\v^2\int_0^\infty \big|\big(\sqrt{\rho^{\v}(t,r)}\big)_r\big|^2\,r^{N-1}\dd r
  +\v\int_{\mathbb{R}_+^2} \big|\big((\rho^{\v}(s,r))^{\frac{\gamma}{2}}\big)_r\big|^2\,r^{N-1}\dd r\dd s\nonumber\\[1mm]
&\quad    \leqq C(E_0+1)\quad\,\,\mbox{for $t\geqq 0$},\label{5.8-1}\\[3mm]
&\int_0^T\int_{d}^D (\rho^{\v})^{\g+1}(t,r)\,\dd r\dd t\leqq C(d,D,T,E_0),\label{5.9-1}\\[1mm]
&\int_0^T\int_0^D  \big(\rho^{\v}|u^\v|^3+(\rho^{\v})^{\gamma+\theta}\big)(t,r)\,r^{N-1} \dd r\dd t
  \leqq C(D,T,E_0),\label{5.10-1}
\end{align}
where $[d,D]\Subset (0,\infty)$.
\end{proposition}

\medskip
We now show that
\begin{align}\label{3.35-1}
(\rho^\v, \M^\v)(t,\mathbf{x})=(\rho^\v(t,r), m^\v(t,r)\frac{\mathbf{x}}{r})
\end{align}
is a weak solution of the Cauchy problem \eqref{1.1} and \eqref{initial-data}
in $\mathbb{R}^N$ in the sense of Definition \ref{definition-NS}.

\begin{lemma}\label{lem5.8}
Let $0\leqq t_1<t_2\leqq T$,
and let $\zeta(t,\mathbf{x})\in C^1([0,T]\times\mathbb{R}^N)$ be any smooth function with compact support.
Then
\begin{align}\label{5.39}
&\int_{\mathbb{R}^N} \rho^{\v}(t_2,\mathbf{x}) \zeta(t_2,\mathbf{x})\, \dd\mathbf{x}\nonumber\\
&=\int_{\mathbb{R}^N} \rho^{\v}(t_1,\mathbf{x}) \zeta(t_1,\mathbf{x})\, \dd\mathbf{x}
 +\int_{t_1}^{t_2} \int_{\mathbb{R}^N} \big(\rho^{\v} \zeta_t + \M^{\v}\cdot\nabla\zeta\big)\,\dd\mathbf{x}\dd t.
 \end{align}
 \end{lemma}

\noindent{\bf Proof.}
Notice that $(\rho^{\v,\d}, \mathcal{M}^{\v,\d})$ is a strong solution of \eqref{1.1} and \eqref{5.5-2}
over $[0,\infty)\times \big(\mathbb{R}^N\setminus B_\delta(\mathbf{0})\big)$.
It follows from $\eqref{1.1}_1$ and a direct calculation that
 \begin{align}\label{5.40}
0&=\int_{t_1}^{t_2} \int_{\mathbb{R}^N\backslash B_\d(\mathbf{0})}
  \big((\rho^{\v,\d})_t +\mbox{div} \M^{\v,\d} \big) \zeta(t,\mathbf{x})\, \dd\mathbf{x}\dd t\nonumber\\
&=\int_{\mathbb{R}^N\backslash B_\d(\mathbf{0})} \rho^{\v,\d} \zeta\, d\mathbf{x}\Big|_{t_1}^{t_2}
-\int_{t_1}^{t_2} \int_{\mathbb{R}^N\backslash B_\d(\mathbf{0})} \big(\rho^{\v,\d} \zeta_t + \M^{\v,\d}\cdot\nabla\zeta\big)\, \dd\mathbf{x}\dd t\nonumber\\
&=\int_{\mathbb{R}^N} \rho^{\v,\d} \zeta\, \dd\mathbf{x}\Big|_{t_1}^{t_2}
-\int_{t_1}^{t_2} \int_{\mathbb{R}^N} \big(\rho^{\v,\d} \zeta_t + \M^{\v,\d} \cdot\nabla\zeta\big)\, \dd\mathbf{x}\dd t,
\end{align}
where we have used the fact that $(\rho^{\v,\d},m^{\v,\d})$ is extended by zero in $[0,T]\times[0,\d)$.

Notice that, for $i=1,2$,
\begin{align}\label{5.41}
&\Big| \int_{\mathbb{R}^N} \big(\rho^{\v,\d}-\rho^{\v}\big)(t_i,\mathbf{x}) \zeta(t_i,\mathbf{x})\, \dd\mathbf{x} \Big|\nonumber\\
&\leqq \Big| \int_{\mathbb{R}^N\backslash B_\s(\mathbf{0})} \big(\rho^{\v,\d}-\rho^{\v}\big)(t_i,\mathbf{x}) \zeta(t_i,\mathbf{x})\, \dd\mathbf{x} \Big|
\nonumber\\
&\quad\, +\Big| \int_{B_\s(\mathbf{0})} \big(\rho^{\v,\d}-\rho^{\v}\big)(t_i,\mathbf{x}) \zeta(t_i,\mathbf{x})\, \dd\mathbf{x} \Big|.
 \end{align}
Denote
\begin{align}\label{3.39-1}
\phi(t,r):=\int_{\partial B_1(\mathbf{0})}\zeta(t,r\omega)\, \dd\omega\in C^1_{0}([0,T]\times[0,\infty)).
\end{align}
Then, with \eqref{5.12}, for any fixed $\sigma>0$, we have
\begin{align}\label{5.42}
&\lim_{\d\rightarrow0+}\Big|\int_{\mathbb{R}^N\backslash B_\s(\mathbf{0})} \big(\rho^{\v,\d}-\rho^{\v}\big)(t_i,\mathbf{x}) \zeta(t_i,\mathbf{x})\, \dd\mathbf{x} \Big|\nonumber\\
&=\omega_N\lim_{\d\rightarrow0+}\Big|\int_{\sigma}^\infty \big(\rho^{\v,\d}-\rho^{\v}\big)(t_i,r) \phi(t_i,r)\, r^{N-1}\dd r\Big|=0.
\end{align}

Using \eqref{5.7} and \eqref{5.7-1}, we obtain
\begin{align}\label{5.42-1}
&\Big|\int_{B_\s(\mathbf{0})} (\rho^{\v,\d}-\rho^{\v})(t_i,\mathbf{x}) \zeta(t_i,\mathbf{x})\, \dd\mathbf{x} \Big|\nonumber\\
&\leqq C\|\zeta\|_{L^\infty}\Big\{\int_0^\sigma \big((\rho^{\v,\d})^{\gamma}+(\rho^{\v})^{\gamma}\big)\, r^{N-1} \dd r \Big\}^{\frac1\gamma}
  \sigma^{N(1-\frac1\gamma)}
  \nonumber\\
&\leqq C(E_0) \|\zeta\|_{L^\infty} \sigma^{N(1-\frac1\gamma)}\rightarrow0\qquad\,\,  \mbox{as $\sigma \rightarrow0$},
\end{align}
which, along with  \eqref{5.41} and \eqref{5.42}, yields
\begin{align}\label{5.43}
\lim_{\d\rightarrow0+}\int_{\mathbb{R}^N} \rho^{\v,\d}(t_i,\mathbf{x}) \zeta(t_i,\mathbf{x})\, \dd\mathbf{x}
=\int_{\mathbb{R}^N} \rho^{\v}(t_i,\mathbf{x}) \zeta(t_i,\mathbf{x})\, \dd\mathbf{x}\qquad\,\,\mbox{for $i=1,2$}.
\end{align}

From \eqref{3.39-1}, a direct calculation shows
\begin{align}\label{3.43-1}
\phi_r=\int_{\partial B_1(\mathbf{0})}\omega\cdot \nabla\zeta(t,r\omega)\, \dd\omega
\end{align}
which, with \eqref{5.12} and Lemma \ref{lem5.6}, implies
 \begin{align}\label{5.44}
 &\lim_{\delta\rightarrow0+}\int_{t_1}^{t_2} \int_{\mathbb{R}^N\backslash B_\s(\mathbf{0})}
   \big(\rho^{\v,\d} \zeta_t + \M^{\v,\d}\cdot\nabla\zeta\big)\, \dd\mathbf{x}\dd t\nonumber\\
& =\omega_N \lim_{\delta\rightarrow0+}\int_{t_1}^{t_2}\int_\sigma^\infty \big(\rho^{\v,\d} \phi_t + m^{\v,\d} \phi_r\big)\, r^{N-1}\dd r\dd t\nonumber\\
 &=\omega_N\int_{t_1}^{t_2}\int_\sigma^\infty \big(\rho^{\v} \phi_t + m^{\v} \phi_r\big)\, r^{N-1}\dd r\dd t\nonumber\\
& =\int_{t_1}^{t_2} \int_{\mathbb{R}^N\backslash B_\s(\mathbf{0})} \big(\rho^{\v} \zeta_t + \M^{\v} \cdot\nabla\zeta\big)\, \dd\mathbf{x}\dd t.
 \end{align}
Similar to that in \eqref{5.42-1}, we also have
\begin{align*}
&\Big|\int_{t_1}^{t_2} \int_{B_\s(\mathbf{0})}  (\rho^{\v,\d}-\rho^{\v}) \zeta_t\, \dd\mathbf{x}\dd t\Big|
\leqq C(T,E_0) \|\zeta_t\|_{L^\infty}\, \sigma^{N(1-\frac1\gamma)},\\[1.5mm]
&\Big|\int_{t_1}^{t_2} \int_{B_\s(\mathbf{0})}  \big(\M^{\v,\d}-\M^{\v}\big) \cdot\nabla \zeta\, \dd\mathbf{x}\dd t\Big|\nonumber\\
&\leqq C\|\nabla\zeta\|_{L^\infty}\Big\{\int_{t_1}^{t_2}\int_0^\sigma \big(\frac{|m^{\v,\d}|^2}{\rho^{\v,\d}}
  +\frac{|m^{\v}|^2}{\rho^{\v}}\big)(t,r) r^{N-1}\,\dd r \dd t\Big\}^{\frac12}\nonumber\\
&\qquad\times \Big\{\int_{t_1}^{t_2}\int_0^\sigma \big(\rho^{\v,\d}+\rho^{\v}\big)(t,r) r^{N-1}\,\dd r \dd t\Big\}^{\frac12}\nonumber\\
&\leqq C(T,E_0) \|\nabla\zeta\|_{L^\infty} \sigma^{\frac{N}{2}(1-\frac1\gamma)},\nonumber
\end{align*}
which,  with \eqref{5.44},
yields
 \begin{align}\label{5.47}
\lim_{\delta\rightarrow0+}\int_{t_1}^{t_2} \int_{\mathbb{R}^N}
   \big(\rho^{\v,\d} \zeta_t + \M^{\v,\d} \cdot\nabla\zeta\big)\, \dd\mathbf{x}\dd t
=\int_{t_1}^{t_2} \int_{\mathbb{R}^N} \big(\rho^{\v} \zeta_t + \M^{\v} \cdot\nabla\zeta\big)\, \dd\mathbf{x}\dd t.
\end{align}
Combining \eqref{5.43} and \eqref{5.47} with \eqref{5.40}, we conclude \eqref{5.39}.
$\hfill\Box$

\medskip
\begin{lemma}\label{lem5.9}
Let  $\psi(t,\mathbf{x})\in \left(C^2_0([0,\infty)\times \mathbb{R}^N)\right)^N$ be any smooth function
with $\mbox{supp}\, \psi \Subset [0, T)\times \mathbb{R}^N$ for some fixed $T\in (0, \infty)$.
Then
\begin{align}
&\int_{\mathbb{R}_+^{N+1}} \Big\{\M^{\v}\cdot\partial_t\psi
+\frac{\M^{\v}}{\sqrt{\rho^{\v}}} \cdot \big(\frac{\M^{\v}}{\sqrt{\rho^{\v}}}\cdot \nabla\big)\psi + p(\rho^{\v})\,\mbox{\rm div}\,\psi\Big\}\, \dd\mathbf{x}\dd t\nonumber\\
&\quad   +\int_{\mathbb{R}^N} \M_0^\v(\mathbf{x})\cdot \psi(0,\mathbf{x})\,\dd\mathbf{x} \nonumber\\
&=-\v\int_{\mathbb{R}_+^{N+1}}
\Big\{\frac{1}{2}\M^{\v}\cdot \big(\Delta \psi+\nabla\mbox{\rm div}\,\psi \big)
 + \frac{\M^{\v}}{\sqrt{\rho^{\v}}} \cdot \big(\nabla\sqrt{\rho^{\v}}\cdot \nabla\big)\psi \nonumber\\
&\qquad \qquad \qquad+ \nabla\sqrt{\rho^{\v}}  \cdot \big(\frac{\M^{\v}}{\sqrt{\rho^{\v}}}\cdot \nabla\big)\psi\Big\}  \dd\mathbf{x}\dd t \label{5.48}\\
&=\sqrt{\v}\int_{\mathbb{R}_+^{N+1}}
\sqrt{\rho^{\v}} \Big\{V^{\v}  \frac{\mathbf{x}\otimes\mathbf{x}}{r^2}
  +\frac{\sqrt{\v}}{r}\frac{m^\v}{\sqrt{\rho^\v}}\big(I_{N\times N}-\frac{\mathbf{x}\otimes\mathbf{x}}{r^2}\big)\Big\}: \nabla\psi\, \dd\mathbf{x}\dd t,\label{5.48-5}
\end{align}
where $V^{\v}(t,r)\in L^2(0,T; L^2(\mathbb{R}^N))$ is a function such that
$$
\displaystyle\int_0^T\int_{\mathbb{R}^N} |V^{\v}(t,\mathbf{x})|^2  \dd\mathbf{x}\dd t\leqq CE_0
$$
for some $C>0$, independent of $T>0$.
\end{lemma}

\noindent{\bf Proof.}
Let $\psi=(\psi_1,\cdots,\psi_N) \in \big(C_0^2([0,\infty)\times\mathbb{R}^N)\big)^N$ be a smooth
function with $\mbox{supp}\, \psi \Subset [0, T)\times \mathbb{R}^N$.
For any given $\sigma\in (0,1]$, let $\chi_\sigma(r)\in C^\infty(\mathbb{R})$ be a cut-off function satisfying
\begin{align}\label{3.49-1}
\begin{split}
&\displaystyle\chi_\sigma(r)=0\,\,\,  \mbox{for $r\leqq \sigma$},\quad
\displaystyle\chi_\sigma(r)=1\,\,  \mbox{for $r\geqq 2\sigma$},\\
&\displaystyle|\chi_\sigma'(r)|\leqq \frac{C}{\sigma},\quad  |\chi_\sigma''(r)|\leqq \frac{C}{\sigma^2}.
\end{split}
\end{align}
Denote $\Psi_\sigma(t,\textbf{x}):=\psi(t,\textbf{x}) \chi_\sigma(|\textbf{x}|)$.

Taking $\d$ small enough so that  $\d\leqq \sigma$, then it follows from $\eqref{1.1}_2$ and integration by parts that
 \begin{align}\label{5.48-1}
&\int_{\mathbb{R}_+^{N+1}}
\Big\{\M^{\v,\d} \cdot\partial_t\Psi_\sigma
  +\frac{\M^{\v,\d}}{\sqrt{\rho^{\v,\d}}} \cdot \big(\frac{\M^{\v,\d}}{\sqrt{\rho^{\v,\d}}}\cdot \nabla\big)\Psi_{\sigma}
  + p(\rho^{\v,\d})\,\mbox{div}\,\Psi_{\sigma}\Big\}\,  \dd\mathbf{x}\dd t\nonumber\\
&\quad+\int_{\mathbb{R}^N} \M_0^{\v,\d}(\mathbf{x})\cdot \Psi_\s(0,\mathbf{x})\,\dd\mathbf{x}\nonumber\\[1mm]
&=:J_1^{\v,\d}+J_2^{\v,\d},
\end{align}
 where
 \begin{align}
 J_1^{\v,\d}:=&\,\d\v\int_{\mathbb{R}_+^{N+1}}
 (\rho^{\v,\d})^\alpha \Big\{ D(\frac{\M^{\v,\d}}{\rho^{\v,\d}}): \nabla\Psi_{\sigma}
 +(\alpha-1) \mbox{div}\big(\frac{\M^{\v,\d}}{\rho^{\v,\d}}\big)\,\mbox{div}\Psi_{\sigma}\Big\}\, \dd\mathbf{x}\dd t,\label{5.49}\\
 J_2^{\v,\d}:&=
 -\v\int_{\mathbb{R}_+^{N+1}}
 \Big\{
  \frac{1}{2}\M^{\v,\d} \cdot \big(\Delta \Psi_{\sigma}+\nabla\mbox{div}\Psi_{\sigma} \big)
  +\frac{\M^{\v,\d}}{\sqrt{\rho^{\v,\d}}}\cdot \big(\nabla\sqrt{\rho^{\v,\d}}\cdot \nabla\big)\Psi_{\sigma}\nonumber\\
 &\qquad\qquad\quad\,\,\quad + \nabla\sqrt{\rho^{\v,\d}}\cdot \big(\frac{\M^{\v,\d}}{\sqrt{\rho^{\v,\d}}}\cdot\nabla\big)\Psi_{\sigma}\Big\}
  \dd\mathbf{x} \dd t \nonumber\\
 &=\v\int_{\mathbb{R}_+^{N+1}}
 \sqrt{\rho^{\v,\d}} \,\sqrt{\rho^{\v,\d}} D(\frac{\M^{\v,\d}}{\rho^{\v,\d}}): \nabla\Psi_{\sigma}\, \dd\mathbf{x} \dd t. \label{5.48-3}
 \end{align}

A direct calculation leads to
\begin{align}\label{5.51}
\partial_i \big(\frac{\mathcal{M}^{\v,\d}_j}{\rho^{\v,\d}}\big)=u^{\v,\d}_r \frac{x_ix_j}{r^2}+\frac{u^{\v,\d}}{r} \big(\d_{ij}-\frac{x_i x_j}{r^2}\big).
\end{align}
Using \eqref{5.7}, there exists a function $V^\v(t,r)$ so that
 \begin{align}\label{5.52-1}
\sqrt{\v}\sqrt{\rho^{\v,\d}} D(\frac{\mathcal{M}^{\v,\d}_j}{\rho^{\v,\d}})
\rightharpoonup V^{\v} \frac{\mathbf{x}\otimes\mathbf{x}}{r^2}
  +\frac{\sqrt{\v}}{r}\frac{m^\v}{\sqrt{\rho^\v}}\big(I_{N\times N}-\frac{\mathbf{x}\otimes\mathbf{x}}{r^2}\big)
\end{align}
in $ L^2(0,T; (L^2(\mathbb{R}^N\backslash B_\sigma(\mathbf{0})))^{N\times N}) $ as  $\d\rightarrow0+$ for any given $\sigma>0$.
Moreover, we have
\begin{align}\label{5.52-2}
\int_0^T\int_{\mathbb{R}^N} \left| V^\v\right|^2 \dd\mathbf{x} \dd t \leqq C E_0.
\end{align}

It follows from  \eqref{5.7} and \eqref{5.51} that
\begin{align}\label{5.52}
|J_1^{\v,\d}|&\leqq C(\|\psi\|_{C^1},\mbox{supp}\,\psi, \sigma)\sqrt{\v\d}
\Big\{\int_{{\rm supp}\,\Psi_\sigma}  (\rho^{\v,\d})^{\alpha}\,r^{N-1}\dd r\dd t\Big\}^{\frac12}\nonumber\\
&\qquad\times  \bigg\{\d\v\int_{{\rm supp}\,\Psi_\sigma} (\rho^{\v,\d})^{\alpha} \Big(|u^{\v,\d}_r|^2+\frac{ |u^{\v,\d}|^2}{r^2} \Big)\, r^{N-1}\dd r\dd t\bigg\}^{\frac12} \nonumber\\[1mm]
&\leqq C(\|\psi\|_{C^1}, {\rm supp}\,\psi,\sigma, E_0)\sqrt{\v\d}\rightarrow0\qquad\,\,   \mbox{as $\d\rightarrow0+$}.
\end{align}

Denote
\begin{align}
\phi_{1\sigma}(t,r):=\int_{\partial B_1(\mathbf{0})}\big(\omega\cdot(\Delta \Psi_{\sigma})(t,r\omega)+\omega\cdot(\nabla\mbox{div}\Psi_{\sigma})(t,r\omega) \big)\, \dd\omega.\nonumber
\end{align}
Then it is clear that $\phi_{1\sigma}\in C_0^2([0,T]\times (0,\infty))$.
Thus, using Lemma \ref{lem5.6},
we have
\begin{align}\label{3.55-1}
&\int_{\mathbb{R}_+^{N+1}} \M^{\v,\d}\cdot \big(\Delta \Psi_{\sigma}+\nabla\mbox{div}\Psi_{\sigma} \big) \,\dd\mathbf{x}\dd t\nonumber\\
&=\omega_N\int_{\mathbb{R}_+^{2}} m^{\v,\d} \phi_{1\sigma}\, r^{N-1}\, \dd r \dd t\nonumber\\
&\rightarrow \omega_N\int_{\mathbb{R}_+^{2}} m^\v\phi_{1\sigma}\, r^{N-1}\,\dd r \dd t\nonumber\\
&=\int_{\mathbb{R}_+^{N+1}} \M^\v\cdot \big(\Delta \Psi_{\sigma}+\nabla\mbox{div}\Psi_{\sigma} \big)\,  \dd\mathbf{x}\dd t
\qquad  \mbox{as $\d\rightarrow0$}.
\end{align}
Similarly, using Lemmas \ref{lem5.2} and \ref{lem5.6},
we can prove
\begin{align}\label{3.56-1}
&\int_{\mathbb{R}_+^{2}}
\Big\{\frac{\M^{\v,\d}}{\sqrt{\rho^{\v,\d}}} \cdot \big(\nabla\sqrt{\rho^{\v,\d}}\cdot \nabla\big)
+ \nabla\sqrt{\rho^{\v,\d}}  \cdot \big(\frac{\M^{\v,\d}}{\sqrt{\rho^{\v,\d}}}\cdot \nabla\big)\Big\}\Psi_{\sigma}\, \dd\mathbf{x} \dd t\nonumber\\
&\rightarrow \int_{\mathbb{R}_+^{2}}
\Big\{ \frac{\M^{\v}}{\sqrt{\rho^{\v}}} \cdot \big(\nabla\sqrt{\rho^{\v}}\cdot \nabla\big)
+ \nabla\sqrt{\rho^{\v}}\cdot \big(\frac{\M^{\v}}{\sqrt{\rho^{\v}}}\cdot \nabla\big)\Big\}\Psi_{\sigma}\, \dd\mathbf{x}\dd t\,\,\,  \mbox{as $\d\rightarrow0$}.
\end{align}
Combining \eqref{5.52-1} with \eqref{3.55-1}--\eqref{3.56-1}, we obtain that, as $\d\rightarrow0$,
\begin{align}
J_2^{\v,\d}\rightarrow
& -\v\int_{\mathbb{R}_+^{N+1}}
\Big\{\frac{1}{2}\M^{\v}\cdot \big(\Delta \Psi_{\sigma}+\nabla\mbox{div}\Psi_{\sigma} \big)
+\frac{\M^{\v}}{\sqrt{\rho^{\v}}}\cdot \big(\nabla\sqrt{\rho^{\v}}\cdot \nabla\big)
\Psi_{\sigma} \nonumber\\
&\qquad\qquad \quad\, + \nabla\sqrt{\rho^{\v}}\cdot \big(\frac{\M^{\v}}{\sqrt{\rho^{\v}}}\cdot \nabla\big)\Psi_{\sigma}\Big\}\dd\mathbf{x}\dd t\nonumber\\
&=\sqrt{\v}\int_{\mathbb{R}_+^{N+1}}
\sqrt{\rho^{\v}} \Big\{V^{\v}  \frac{\mathbf{x}\otimes\mathbf{x}}{r^2}
 +\frac{\sqrt{\v}}{r}\frac{m^\v}{\sqrt{\rho^\v}}
  \big(I_{N\times N}-\frac{\mathbf{x}\otimes\mathbf{x}}{r^2}\big)\Big\}: \nabla\Psi_{\sigma}\, \dd\mathbf{x}
  \dd t. \label{3.57-2}
\end{align}
Also, by similar arguments as in \eqref{3.55-1}, applying  Lemma \ref{lem5.2}, Corollary \ref{cor5.3}, and Lemma \ref{lem5.6},
we have
 \begin{align}
&\int_{\mathbb{R}_+^{N+1}}
\Big\{\M^{\v,\d} \cdot\partial_t\Psi_\sigma
    +\frac{\M^{\v,\d}}{\sqrt{\rho^{\v,\d}}} \cdot \big(\frac{\M^{\v,\d}}{\sqrt{\rho^{\v,\d}}}\cdot \nabla\big)\Psi_{\sigma}
    +p(\rho^{\v,\d})\,\mbox{div} \Psi_{\sigma}\Big\}\, \dd\mathbf{x}\dd t\nonumber\\
&\quad  +\int_{\mathbb{R}^N} \M_0^{\v,\d}(\mathbf{x})\cdot \Psi_\s(0,\mathbf{x})\, \dd\textbf{x}\nonumber\\[1.5mm]
&\rightarrow \int_{\mathbb{R}_+^{N+1}}
\Big\{\M^{\v} \cdot\partial_t\Psi_\sigma
   +\frac{\M^{\v}}{\sqrt{\rho^{\v}}}\cdot \big(\frac{\M^{\v}}{\sqrt{\rho^{\v}}}\cdot \nabla\big)\Psi_{\sigma} +p(\rho^{\v}) \mbox{div} \Psi_{\sigma}
   \Big\}\, \dd\mathbf{x}\dd t \nonumber\\
&\quad  +\int_{\mathbb{R}^N} \M_0^{\v}(\mathbf{x})\cdot \Psi_\s(0,\mathbf{x})\, \dd\mathbf{x}\nonumber
\end{align}
as $\d\rightarrow0$, which, with
\eqref{3.57-2}, yields
\begin{align}\label{3.59-1}
&\int_{\mathbb{R}_+^{N+1}}
\Big\{\M^{\v} \cdot\partial_t\Psi_\sigma
   + \frac{\M^{\v}}{\sqrt{\rho^{\v}}} \cdot \big(\frac{\M^{\v}}{\sqrt{\rho^{\v}}}\cdot \nabla\big)\Psi_{\sigma}
  + p(\rho^{\v})\,\mbox{div} \Psi_{\sigma}
    \Big\}\,\dd\mathbf{x}\dd t \nonumber\\
 &\quad +\int_{\mathbb{R}^N} \M_0^\v(\mathbf{x}) \cdot \Psi_\s(0,\mathbf{x})\, \dd\mathbf{x} \nonumber\\
&=
 -\v\int_{\mathbb{R}_+^{N+1}}
 \Big\{\frac{1}{2} \M^{\v} \cdot \big(\Delta \Psi_{\sigma}+\nabla\mbox{div}\Psi_{\sigma} \big)
 +\frac{\M^{\v}}{\sqrt{\rho^{\v}}} \cdot \big(\nabla\sqrt{\rho^{\v}}\cdot \nabla\big)\Psi_{\sigma}\nonumber\\
&\qquad\qquad\qquad + \nabla\sqrt{\rho^{\v}}\cdot \big(\frac{\M^{\v}}{\sqrt{\rho^{\v}}}\cdot \nabla\big)\Psi_{\sigma}\Big\}\, \dd\mathbf{x}\dd t\nonumber\\
&=\sqrt{\v}\int_{\mathbb{R}_+^{N+1}}
\sqrt{\rho^{\v}} \Big\{V^{\v}  \frac{\mathbf{x}\otimes\mathbf{x}}{r^2}
 +\frac{\sqrt{\v}}{r}\frac{m^\v}{\sqrt{\rho^\v}}
   \big(I_{N\times N}-\frac{\mathbf{x}\otimes\mathbf{x}}{r^2}\big)\Big\}: \nabla\Psi_{\sigma}\, \dd\mathbf{x}\dd t.
\end{align}

Next, we consider the limit, $\sigma\rightarrow0$, in \eqref{3.59-1}. We first define
\begin{align}\label{3.63-2}
\varphi(t,r):&=\int_{\partial B_1(\mathbf{0})} \omega\cdot\psi(t,r\omega)\, \dd\omega\nonumber\\
&=\frac{1}{r^{N-1}}\int_{\partial B_r(\mathbf{0})} \omega\cdot\psi(t,\mathbf{y})\, \dd S_{\mathbf{y}}\nonumber\\
&=\frac{1}{r^{N-1}}\int_{B_r(\mathbf{0})} \mbox{div}\,\psi(t,\mathbf{y})\, \dd\mathbf{y},
\end{align}
which implies
\begin{align}\label{3.63-1}
|\varphi(t,r)|\leqq C(\|\psi\|_{C^1}) r;
\end{align}
also see \cite{Song Jiang, Schrecker}.
Using \eqref{3.63-2}, Lebesgue's dominated convergence theorem, and Proposition \ref{prop5.2}, we have
\begin{align}\label{3.61}
&\lim_{\s\rightarrow0}\left\{\int_{\mathbb{R}_+^{N+1}}
\M^{\v}\cdot\partial_t\Psi_\sigma\, \dd\mathbf{x}\dd t
  + \int_{\mathbb{R}^N} \M_0^\v(\mathbf{x}) \cdot \Psi_\s(0,\mathbf{x})\,\dd\mathbf{x}\right\}\nonumber\\
&=\omega_N\lim_{\s\rightarrow0}\left\{\int_{\mathbb{R}_+^{2}}
m^{\v} \,\partial_t\varphi \, \chi_\sigma(r)\,r^{N-1}\dd r\dd t
  +\int_0^\infty m_0^\v(r) \varphi(0,r) \chi_\sigma(r)\,r^{N-1}\dd r\right\}\nonumber\\
&=\omega_N\int_{\mathbb{R}_+^{2}}
  m^{\v} \,\partial_t\varphi \, r^{N-1} \dd r\dd t
  +\omega_N\int_0^\infty m_0^\v(r)\varphi(0,r)\,r^{N-1}\dd r\nonumber\\
&=\int_{\mathbb{R}_+^{N+1}}
\M^{\v}\cdot\partial_t\psi\, \dd\mathbf{x}\dd t
  + \int_{\mathbb{R}^N} \M_0^\v(\mathbf{x})\cdot \psi(0,\mathbf{x})\,\dd\mathbf{x}.
\end{align}

Using \eqref{3.63-1} and Proposition \ref{prop5.2}, we have
\begin{align}
&\left|\int_{\mathbb{R}_+^{N+1}} \big(\frac{|m^\v|^2}{\rho^{\v}}+p(\rho^\v)\big)
   \big(\psi\cdot \frac{\mathbf{x}}{r}\big)\chi'_\sigma(r)
  \,\dd\mathbf{x}\dd t\right|\nonumber\\
&  \leqq C \int_0^T\int_{\sigma}^{2\sigma} \big(\frac{|m^\v|^2}{\rho^{\v}}+p(\rho^\v)\big) |\varphi(t,r) \chi'_\sigma(r)|\,r^{N-1}\dd r\dd t\nonumber\\
& \leqq C\int_0^T\int_{\sigma}^{2\sigma} \big(\frac{|m^\v|^2}{\rho^{\v}}+p(\rho^\v)\big)\,r^{N-1}\dd r\dd t\rightarrow 0\quad\,\,   \mbox{as $\sigma\rightarrow0$}, \label{3.64}\\[1.5mm]
&\left|\v\int_{\mathbb{R}_+^{N+1}} \frac{m^\v}{\sqrt{\rho^{\v}}} \big(\sqrt{\rho^{\v}}\big)_r
\big(\psi\cdot \frac{\mathbf{x}}{r}\big) \chi'_\sigma(r)\,\dd\mathbf{x}\dd t\right|\nonumber\\
&\leqq C\v\int_0^T\int_{\sigma}^{2\sigma} \big|\frac{m^\v}{\sqrt{\rho^{\v}}}\big(\sqrt{\rho^{\v}}\big)_r\big|
   \,|\varphi(t,r) \chi'_\sigma(r)|\,r^{N-1}\dd r\dd t\nonumber\\
&\leqq C\int_0^T\int_{\sigma}^{2\sigma}\Big(\frac{|m^\v|^2}{\rho^{\v}}
  + \v^2\big|\big(\sqrt{\rho^{\v}}\big)_r\big|^2\Big)\, r^{N-1}\dd r\dd t\rightarrow 0
\quad\,\,  \mbox{as $\sigma\rightarrow0$},\label{3.66}\\[1.5mm]
&\left|\int_{\mathbb{R}_+^{N+1}} \chi'_\sigma(r)\sqrt{\rho^{\v}} \Big\{V^{\v}  \frac{\mathbf{x}\otimes\mathbf{x}}{r^2}
  +\frac{\sqrt{\v}}{r}\frac{\m^\v}{\sqrt{\rho^{\v}}}
  \big(I_{N\times N}-\frac{\mathbf{x}\otimes\mathbf{x}}{r^2}\big)\Big\}: \big(\psi\otimes \frac{\mathbf{x}}{r}\big)\, \dd\mathbf{x} \dd t\right|\nonumber\\
&=\left|\int_{\mathbb{R}_+^{N+1}} \chi'_\sigma(r)\sqrt{\rho^{\v}}V^{\v}   \big(\psi\cdot  \frac{\mathbf{x}}{r}\big)\, \dd\mathbf{x}\dd t\right|\nonumber\\
& \leqq C\left|\int_0^T\int_{\sigma}^{2\sigma} \sqrt{\rho^{\v}}V^{\v} \, r^{N-1} \dd r \dd t\right|\rightarrow 0 \qquad \mbox{as $\sigma\rightarrow0$}. \label{3.67}
\end{align}
Using \eqref{3.64}--\eqref{3.67}, Lebesgue's dominated convergence theorem, and Proposition \ref{prop5.2}, we obtain
\begin{align}
&\lim_{\s\rightarrow0}\int_{\mathbb{R}_+^{N+1}}
\Big\{ \frac{\M^{\v}}{\sqrt{\rho^{\v}}}\cdot (\frac{\M^{\v}}{\sqrt{\rho^{\v}}}\cdot \nabla)\Psi_{\sigma}
   +p(\rho^{\v})\,\mbox{div} \Psi_{\sigma}\Big\}\, \dd\mathbf{x}\dd t \nonumber\\
&\,\,\,\,=\int_{\mathbb{R}_+^{N+1}}
\Big\{\frac{\M^{\v}}{\sqrt{\rho^{\v}}}\cdot (\frac{\M^{\v}}{\sqrt{\rho^{\v}}}\cdot \nabla)\psi
  +p(\rho^{\v})\,\mbox{div} \psi\Big\}\, \dd\mathbf{x}\dd t,\label{3.68}\\[1mm]
&\lim_{\s\rightarrow0}\int_{\mathbb{R}_+^{N+1}}
\Big\{\frac{\M^{\v}}{\sqrt{\rho^{\v}}}\cdot (\nabla\sqrt{\rho^{\v}}\cdot \nabla)\Psi_{\sigma}
+ (\nabla\sqrt{\rho^{\v}})  \cdot (\frac{\M^{\v}}{\sqrt{\rho^{\v}}}\cdot \nabla)\Psi_{\sigma}\Big\}\, \dd\mathbf{x}\dd t\nonumber\\
&\,\,\,\, =\int_{\mathbb{R}_+^{N+1}}
\Big\{\frac{\M^{\v}}{\sqrt{\rho^{\v}}}\cdot (\nabla\sqrt{\rho^{\v}}\cdot \nabla)\psi
+ (\nabla\sqrt{\rho^{\v}})\cdot (\frac{\M^{\v}}{\sqrt{\rho^{\v}}}\cdot \nabla)\psi\Big\}\, \dd\mathbf{x} \dd t,  \label{3.69}\\[1mm]
&\lim_{\s\rightarrow0}\int_{\mathbb{R}_+^{N+1}}
\sqrt{\rho^{\v}}
 \Big\{V^{\v}  \frac{\mathbf{x}\otimes\mathbf{x}}{r^2}
  +\frac{\sqrt{\v}}{r}\frac{m^\v}{\sqrt{\rho^\v}}\big(I_{N\times N}-\frac{\mathbf{x}\otimes\mathbf{x}}{r^2}\big)\Big\}: \nabla\Psi_{\sigma}\,\dd\mathbf{x} \dd t\nonumber\\
&\,\,\,\,=\int_{\mathbb{R}_+^{N+1}}
\sqrt{\rho^{\v}} \Big\{V^{\v}  \frac{\mathbf{x}\otimes\mathbf{x}}{r^2}
  +\frac{\sqrt{\v}}{r}\frac{m^\v}{\sqrt{\rho^\v}}\big(I_{N\times N}-\frac{\mathbf{x}\otimes\mathbf{x}}{r^2}\big)\Big\}: \nabla\psi\,  \dd\mathbf{x}\dd t.\label{3.70}
\end{align}

We notice that
\begin{align}\label{3.70-1}
\begin{split}
\Delta (\Psi_{\sigma})_i=&\,\chi_\sigma(r) \Delta \psi_i +2\nabla\psi_i \cdot \nabla \chi_\sigma(r)+\psi_i \Delta \chi_\sigma(r),\\
\partial_i\mbox{div}\Psi_{\sigma}=&\, \chi_\sigma(r)\,\partial_i\mbox{div}\psi+\mbox{div}\psi\,\partial_i \chi_\sigma(r) +\partial_i \psi\cdot \nabla \chi_\sigma(r) \\
&\,+\frac{x_i}{r}\chi''_\sigma(r)\,\psi\cdot \frac{\mathbf{x}}{r}
+\chi'_{\sigma}(r)\psi\cdot \big(\frac{\nabla x_i}{r}- \frac{x_i}{r^2}\frac{\mathbf{x}}{r}\big).
\end{split}
\end{align}
It follows from \eqref{3.63-1} and Proposition \ref{prop5.2} that
\begin{align}\label{3.71}
&\left|\sum_{i=1}^N\v\int_{\mathbb{R}_+^{N+1}}
m^{\v} \frac{x_i}{r}
\Big\{ 2\nabla\psi_i \cdot \nabla \chi_\sigma+\psi_i \Delta \chi_\sigma+\mbox{div}\psi\, \partial_i \chi_\sigma(r)+\partial_i \psi\cdot \nabla \chi_\sigma(r) \right.\nonumber\\
&\qquad\qquad\qquad\qquad\,\,\, \left.+\frac{x_i}{r}\chi''_\sigma(r) \big(\psi\cdot \frac{\mathbf{x}}{r}\big)+\chi'_{\sigma}(r)
             \big(\psi\cdot \frac{\nabla x_i}{r}-\big(\psi \cdot \frac{\mathbf{x}}{r}\big)\frac{x_i}{r^2}\big) \Big\}\,  \dd\mathbf{x}\dd t\right|\nonumber\\
&\leqq C(\|\psi\|_{C^1})\int_0^T\int_{\sigma}^{2\sigma} \v|m^{\v}|
   \Big(|\chi'_\sigma(r)|+\frac1{r}\varphi(r)|\chi'_\sigma(r)|+\varphi(r)|\chi''_\sigma(r)| \Big)\,r^{N-1}\dd r \dd t\nonumber\\
&\leqq C(\|\psi\|_{C^1})\int_0^T\int_{\sigma}^{2\sigma} \v|m^\v|\, r^{N-2} \dd r \dd t\nonumber\\
&\leqq C(\|\psi\|_{C^1})\left\{\int_0^T\int_{\sigma}^{2\sigma} \rho^\v \, r^{N-1} \dd r \dd t \right\}^{\frac12}
  \left\{\v\int_0^T\int_{\sigma}^{2\sigma} \frac{|m^\v|^2}{\rho^\v} \, r^{N-3} \dd r \dd t \right\}^{\frac12}\nonumber\\[1.5mm]
&\rightarrow 0\qquad\mbox{as $\sigma\rightarrow0$}.
\end{align}

Using \eqref{3.70-1}--\eqref{3.71},
Lebesgue's dominated convergence theorem, and Proposition \ref{prop5.2}, we have
\begin{align}\label{3.72}
&\lim_{\sigma\rightarrow0}\v\int_{\mathbb{R}_+^{N+1}}
\M^{\v}\cdot \big(\Delta \Psi_{\sigma}+\nabla\mbox{div}\Psi_{\sigma} \big)\,  \dd\mathbf{x}\dd t\nonumber\\
&=\v\int_{\mathbb{R}_+^{N+1}}
\M^{\v} \cdot \big(\Delta \psi+\nabla\mbox{div}\psi\big)\, \dd\mathbf{x}\dd t.
\end{align}
Substituting \eqref{3.61}, \eqref{3.68}--\eqref{3.70}, and \eqref{3.72} into \eqref{3.59-1},
we conclude \eqref{5.48}--\eqref{5.48-5}.
$\hfill\Box$

\begin{remark}
It is not so clear to show that the right-hand side terms of \eqref{5.48} vanish  as $\v\rightarrow0$ by direct arguments. However,
we can prove the vanishing of these terms by using \eqref{5.48-5}, which is the main reason why the form of \eqref{5.48-5} is important to us.
\end{remark}

\bigskip
We also need the $H_{\rm loc}^{-1}$--compactness of weak entropy dissipation
measures of $(\rho^\v, m^\v)$.

\begin{lemma}[$H_{\rm loc}^{-1}$--compactness]\label{lem7.1}
Let $(\eta, q)$ be a weak entropy pair defined in \eqref{weakentropy}
for any smooth compact supported function $\psi(s)$ on $\mathbb{R}$.
Then, for $\v\in (0, 1]$,
\begin{align}\label{7.1}
\partial_t\eta(\rho^\v,m^\v)+\partial_rq(\rho^\v,m^\v) \qquad   \mbox{is compact in $ H^{-1}_{\rm loc}(\mathbb{R}^2_+)$}.
\end{align}
\end{lemma}

\smallskip
\noindent{\bf Proof.}
To obtain \eqref{7.1},
we have to be careful since $(\rho^\v, \mathcal{M}^{\v})$ is a weak solution
of the Navier-Stokes equations \eqref{1.1}.
In fact, we first have to study the equation for $\partial_t\eta(\rho^\v,m^\v)+\partial_rq(\rho^\v,m^\v) $ in
the distributional sense, which is much complicated than that in \cite{Chen6,Chen7}.
We divide the proof into five steps.

\smallskip
1. Since
\begin{align}
\begin{split}
\displaystyle\eta(\rho, m)&=\rho \int_{-1}^1 \psi(u+\rho^{\theta}s)[1-s^2]_+^{\hl}\,\dd s,\\
\displaystyle q(\rho, m)&=\rho \int_{-1}^1 (u+\theta\rho^{\theta}s)\psi(u+\rho^{\theta}s)[1-s^2]_+^{\hl}\,\dd s,\nonumber
\end{split}
\end{align}
then it follows from \cite[Lemma 2.1]{Chen6} that
\begin{align}
&|\eta(\rho, m)|+|q(\rho,m)|\leqq C_{\psi} \rho \qquad\,\, \mbox{for $\gamma\in (1,3]$},\label{7.2}\\
&|\eta(\rho, m)|\leqq C_{\psi} \rho, \,\,\,\, |q(\rho, m)|\leqq C_{\psi} \big(\rho+\rho^{1+\theta}\big)\quad\,\, \mbox{for $\gamma\in (3,\infty)$},\label{7.2-1}\\
&|\partial_\rho\eta (\rho, m)|\leqq C_{\psi} \big(1+\rho^{\theta}\big),
\quad  |\partial_m \eta(\rho, m)|\leqq C_{\psi}.\label{7.6}
\end{align}
Moreover, if $\partial_m \eta(\rho, \rho u)$ is regarded as a function of $(\rho, u)$, then
\begin{align}\label{7.3}
|\partial_{m\rho}\eta|\leqq C_{\psi} \rho^{\theta-1}, \qquad |\partial_{mu}\eta|\leqq C_{\psi}.
\end{align}

\medskip
2. Denote $(\eta^{\v,\d}, q^{\v,\d}):=(\eta, q)(\rho^{\v,\d}, m^{\v,\d})$
and $(\eta^{\v}, q^{\v}):=(\eta, q)(\rho^{\v}, m^{\v})$ for simplicity.
Multiply $\eqref{6.1}_1$ by $\eta_\rho^{\v,\d}$,  $\eqref{6.1}_2$ by $\eta_m^{\v,\d}$, and add them together to obtain
\begin{align}\label{7.4}
&\partial_t \eta^{\v,\d}+\partial_r q^{\v,\d}\nonumber\\
&=-\frac{N-1}{r} m^{\v,\d}\big(\eta_\rho^{\v,\d}+u^{\v,\d} \eta_m^{\v,\d}\big)\nonumber\\
&\quad+ \v \partial_m \eta^{\v,\d} \Big\{ \big(\rho^{\v,\d} (u_r^{\v,\d}+\frac{N-1}{r} u^{\v,\d})\big)_r
  -\frac{N-1}{r} \rho^{\v,\d}_r u^{\v,\d} \Big\}\nonumber\\
&\quad +\v \partial_m \eta^{\v,\d} \Big\{ \alpha\d\big((\rho^{\v,\d})^{\alpha} (u_r^{\v,\d}+\frac{N-1}{r} u^{\v,\d})\big)_r
   -\d\frac{N-1}{r} \big((\rho^{\v,\d})^{\alpha}\big)_r u^{\v,\d} \Big\}.
\end{align}
Let $\phi(t,r)\in C_0^\infty(\mathbb{R}^2_+)$, and let $\d\ll1$ so that $\mbox{\rm supp} (\phi(t,\cdot))\Subset(\d,\infty)$.
Multiplying \eqref{7.4}  by $\phi$ and integrating by parts, we have
\begin{align}\label{7.7}
&\int_{\mathbb{R}^2_+} \big(\partial_t \eta^{\v,\d}+\partial_r q^{\v,\d}\big) \phi\, \dd r\dd t\nonumber\\
&=-\int_{\mathbb{R}^2_+} \frac{N-1}{r} m^{\v,\d}\big(\eta_\rho^{\v,\d}+u^{\v,\d} \eta_m^{\v,\d}\big) \phi\,\dd r\dd t\nonumber\\
&\quad-\v \int_{\mathbb{R}^2_+}\rho^{\v,\d} (\partial_m \eta^{\v,\d})_r \, \big(u_r^{\v,\d}+\frac{N-1}{r} u^{\v,\d}\big) \phi\, \dd r\dd t\nonumber\\
&\quad -\v \int_{\mathbb{R}^2_+}  \rho^{\v,\d}  \partial_m\eta^{\v,\d}  \, \big(u_r^{\v,\d}+\frac{N-1}{r} u^{\v,\d}\big) \phi_r\, \dd r\dd t\nonumber\\
&\quad-\v \int_{\mathbb{R}^2_+} \partial_m \eta^{\v,\d} \frac{N-1}{r} \rho^{\v,\d}_r u^{\v,\d}  \phi\, \dd r\dd t\nonumber\\
&\quad-\alpha\v\d \int_{\mathbb{R}^2_+} (\rho^{\v,\d})^\alpha (\partial_m \eta^{\v,\d})_r \,\big(u_r^{\v,\d}+\frac{N-1}{r} u^{\v,\d}\big) \phi\, \dd r\dd t\nonumber\\
&\quad -\alpha\v\d \int_{\mathbb{R}^2_+} \Big( (\rho^{\v,\d})^\alpha  \partial_m\eta^{\v,\d}  \, \big(u_r^{\v,\d}+\frac{N-1}{r} u^{\v,\d}\big) \phi_r \nonumber\\
&\qquad\qquad\qquad\,\,  +
\partial_m \eta^{\v,\d} \frac{N-1}{r} (\rho^{\v,\d})^{\alpha-1}\rho^{\v,\d}_r u^{\v,\d}  \phi\Big)\, \dd r\dd t\nonumber\\
&:= \sum_{j=1}^6 I_j^{\v,\d}.
\end{align}

\medskip
3. It is direct to see that
\begin{align}\label{7.8}
\eta^{\v,\d}\rightarrow \eta^{\v} \qquad  \mbox{{\it a.e.} in $\{(t,r)\,:\, \rho^{\v}(t,r)\neq0\}\,$ as $\d\rightarrow0+$}.
\end{align}
In $\{(t,r)\,:\, \rho^{\v}(t,r)=0\}$,
\begin{align}\label{7.9}
|\eta^{\v,\d}|\leqq C_{\psi} \rho^{\v,\d}\rightarrow0=\eta^{\v}\qquad  \mbox{as $\d\rightarrow0+$}.
\end{align}
Thus, combining  \eqref{7.8} with \eqref{7.9}, we have
\begin{align}\label{7.10}
\eta^{\v,\d}\rightarrow \eta^{\v} \qquad  \mbox{{\it a.e.} as $\d\rightarrow0+$}.
\end{align}
Similarly, we have
\begin{align}\label{7.11}
q^{\v,\d}\rightarrow q^{\v} \qquad \mbox{{\it a.e.} as $\d\rightarrow0+$}.
\end{align}
Let $K\Subset (0,\infty)$ be any compact subset.  For $\gamma\in(1,3]$,
it follows from \eqref{5.9} and \eqref{7.2} that
\begin{align}\label{7.12}
\int_0^T\int_{K} \big(|\eta^{\v,\d}|+|q^{\v,\d}|\big)^{\gamma+1}\,\dd r\dd t
&\leqq C_{\psi}\int_0^T\int_{K} |\rho^{\v,\d}|^{\gamma+1}\,\dd r\dd t\nonumber\\[1mm]
&\leqq C_{\psi}(K,T,E_0).
\end{align}
For $\gamma\in(3,\infty)$, it follows from \eqref{5.10} and \eqref{7.2-1} that
\begin{align}\label{7.13}
\int_0^T\int_{K} \big(|\eta^{\v,\d}|+|q^{\v,\d}|\big)^{\frac{\gamma+\theta}{1+\theta}}\,\dd r\dd t
&\leqq C_{\psi}\int_0^T\int_{K} \big(|\rho^{\v,\d}|^{\frac{\gamma+\theta}{1+\theta}}+|\rho^{\v,\d}|^{\gamma+\theta}\big)\,\dd r\dd t\nonumber\\[1mm]
&\leqq C_{\psi}(K,T,E_0).
\end{align}
We take  $p_1=\gamma+1>2$ when $\gamma\in(1,3]$, and $p_1=\frac{\gamma+\theta}{1+\theta}>2$ when $\gamma\in (3,\infty)$.
Then it follows from \eqref{7.12}--\eqref{7.13} that
\begin{align}\label{7.14}
(\eta^{\v,\d},  q^{\v,\d})  \qquad   \mbox{is uniformly bounded in $L^{p_1}_{\rm loc}(\mathbb{R}^2_+)$},
\end{align}
which, with \eqref{7.10}--\eqref{7.11}, implies that, up to a subsequence,
\begin{align}\nonumber
(\eta^{\v,\d}, q^{\v,\d}) \rightarrow (\eta^{\v}, q^{\v}) \qquad\,\, \mbox{in $L^2_{\rm loc}(\mathbb{R}^2_+)\,$ as $\delta\rightarrow0+$}.
\end{align}
Thus, for any $\phi\in C^1_0(\mathbb{R}^2_+)$, we see that,
as $\delta\rightarrow 0+$ (up to a subsequence),
\begin{align}\label{7.16}
&\int_{\mathbb{R}^2_+} \big(\partial_t \eta^{\v,\d}+\partial_r q^{\v,\d}\big) \phi\, \dd r\dd t
=-\int_{\mathbb{R}^2_+} \big(\eta^{\v,\d}\partial_t \phi+ q^{\v,\d} \partial_r \phi\big)\, \dd r\dd t
  \nonumber\\
&\rightarrow-\int_{\mathbb{R}^2_+} \big(\eta^{\v}\partial_t \phi+ q^{\v} \partial_r \phi\big)\, \dd r\dd t
 =\int_{\mathbb{R}^2_+} \big(\partial_t \eta^{\v}+\partial_r q^{\v}\big) \phi\,\dd r\dd t.
\end{align}
Furthermore,
$(\eta^{\v}, q^{\v})$  is  uniformly bounded in $L^{p_1}_{\rm loc}(\mathbb{R}^2_+)$ for some $p_1>2$,
which implies that
\begin{align}\label{7.17-1}
\partial_t\eta^{\v}+\partial_rq^{\v} \qquad  \mbox{is uniformly bounded in $\v>0$ in $W^{-1,p_1}_{\rm loc}(\mathbb{R}_+^2)$}.
\end{align}

\medskip
4. Now we estimate the terms on the right-hand side of \eqref{7.7}.
For $I_1^{\v,\d}$, a direct calculation shows that
$|\eta_\rho+u\eta_m|\leqq C_{\psi} \big(1+\rho^{\theta}\big)$,
which,  together Lemma \ref{lem5.6} and  similar arguments in \eqref{7.8}--\eqref{7.10},  leads to
\begin{align}\label{7.18}
\frac{N-1}{r} m^{\v,\d}\big(\eta_\rho^{\v,\d}+u^{\v,\d} \eta_m^{\v,\d}\big)
\,\rightarrow\, \frac{N-1}{r} m^{\v} \big(\eta_\rho^{\v}+u^{\v} \eta_m^{\v}\big) \qquad \mbox{{\it a.e.}$\,$ as  $\d\rightarrow0+$}.
\end{align}
Then it follows from \eqref{5.9}--\eqref{5.10}  that
\begin{align}\label{7.19}
&\int_0^T\int_{K} \Big|\frac{N-1}{r} m^{\v,\d}\big(\eta_\rho^{\v,\d}+u^{\v,\d} \eta_m^{\v,\d}\big)\Big|^{\frac76}\dd r\dd t\nonumber\\
&\leqq C(K)\int_0^T\int_{K}  \big(\rho^{\v,\d}|u^{\v,\d}|^2 +\rho^{\v,\d} +(\rho^{\v,\d})^{\gamma}\big)^{\frac76}\,  \dd r\dd t\nonumber\\
&\leqq
\begin{cases}
\displaystyle C(K)\Big(1+\int_0^T\int_{K} \rho^{\v,\d} |u^{\v,\d}|^3\, \dd r\dd t \Big)^{\frac79}
 \Big(\int_0^T\int_{K}\big(1+|\rho^{\v,\d}|^{\gamma+1}\big)\dd r\dd t \Big)^{\frac29}\\
\qquad\qquad \qquad\qquad \qquad\qquad \qquad\qquad\qquad\qquad \qquad\qquad  \mbox{for $\gamma\in(1,3]$},\\[3mm]
\displaystyle C(K)\Big(1+\int_0^T\int_{K}\rho^{\v,\d}|u^{\v,\d}|^3\, \dd r\dd t \Big)^{\frac79}
  \Big(\int_0^T\int_{K}\big(1+|\rho^{\v,\d}|^{\gamma+\theta}\big)\dd r\dd t \Big)^{\frac29}\\
\qquad\qquad \qquad\qquad \qquad\qquad \qquad\qquad \qquad\qquad \qquad\quad \mbox{for $\gamma\in(3,\infty)$},
\end{cases}
\nonumber\\[1.5mm]
&\leqq C(K,T,E_0).
\end{align}
Using \eqref{7.18}--\eqref{7.19},  we have
\begin{align}
&I^{\v,\d}_1 \rightarrow -\int_{\mathbb{R}^2_+} \frac{N-1}{r} m^{\v} \big(\eta_\rho^{\v}+u^{\v} \eta_m^{\v}\big) \phi\, \dd r\dd t
\,\,  \mbox{as $\d\rightarrow0+$ (up to a subsequence)},\label{7.20}\\
&\int_0^T\int_{K}\Big|\frac{N-1}{r} m^{\v} \big(\eta_\rho^{\v}+u^{\v} \eta_m^{\v}\big) \Big|^{\frac76}\dd r\dd t
\leqq C(K,T,E_0).\label{7.20-1}
\end{align}

For $I_2^{\v,\d}$, $I_4^{\v,\d}$, and $I_5^{\v,\d}$, it follows from \eqref{5.7}--\eqref{5.8} and  \eqref{7.6}--\eqref{7.3} that
\begin{align}
&\int_0^T\int_{K}  \Big|\v \rho^{\v,\d} (\partial_m \eta^{\v,\d})_r \, \big(u_r^{\v,\d}+\frac{N-1}{r} \frac{m^{\v,\d}}{\rho^{\v,\d}}\big) \Big|\,\dd r\dd t\nonumber\\
&\leqq C_{\psi}(K) \int_0^T\int_{K}\big(\v \rho^{\v,\d}|u_r^{\v,\d}|^2
     +\v(\rho^{\v,\d})^{\gamma-2}|\rho_r^{\v,\d}|^2+ \rho^{\v,\d}|u^{\v,\d}|^2 \big)\,\dd r\dd t\nonumber\\
&\leqq C_{\psi}(K,T,E_0),\nonumber\\[2mm]
&\int_0^T\int_{K}  \Big|\v\frac{N-1}{r}\partial_m \eta^{\v,\d}  \rho^{\v,\d}_r u^{\v,\d}\Big|\, \dd r\dd t\nonumber\\
&\leqq C_{\psi}(K) \Big(\v^2\int_0^T\int_{K} \frac{|\rho_r^{\v,\d}|^2}{\rho^{\v,\d}}\, \dd r\dd t\Big)^{\frac12}
\Big(\int_0^T\int_{K} \rho^{\v,\d}|u^{\v,\d}|^2\, \dd r\dd t\Big)^{\frac12}\nonumber\\
&\leqq C_{\psi}(K,T,E_0),\nonumber\\[2mm]
&\int_0^T\int_{K}  \Big| \v\d (\rho^{\v,\d})^\alpha (\partial_m \eta^{\v,\d})_r \,
   \big(u_r^{\v,\d}+\frac{N-1}{r} u^{\v,\d}\big)\big|\,\dd r\dd t\nonumber\\
&\leqq C_{\psi}(K) \int_0^T\int_{K}\v\d (\rho^{\v,\d})^\alpha
   \big( |u_r^{\v,\d}|^2+(\rho^{\v,\d})^{\gamma-3}|\rho_r^{\v,\d}|^2
   +|u^{\v,\d}|^2\big)\,\dd r\dd t\nonumber\\
&\leqq C_{\psi}(K,T,E_0).\nonumber
\end{align}
Thus, there exist local bounded Radon measures $\mu_1^\v, \mu_2^\v$, and $\mu_3^\v$ on $\mathbb{R}^2_+$ so that,
as $\d\rightarrow0+$ (up to a subsequence),
\begin{align}
&-\v \rho^{\v,\d} (\partial_m \eta^{\v,\d})_r \, \big(u_r^{\v,\d}+\frac{N-1}{r} u^{\v,\d}\big) \rightharpoonup \mu_1^{\v},\nonumber\\
&-\v \partial_m \eta^{\v,\d} \frac{N-1}{r} \rho^{\v,\d}_r u^{\v,\d} \rightharpoonup
\mu_2^{\v},\nonumber\\
&-\alpha \v \d(\rho^{\v,\d})^\alpha (\partial_m \eta^{\v,\d})_r \, \big(u_r^{\v,\d}+\frac{N-1}{r} u^{\v,\d}\big) \rightharpoonup
\mu_3^{\v}.\nonumber
\end{align}
In addition,
\begin{equation}\label{7.21-2}
\mu_i^\v((0,T)\times V)\leqq C_{\psi}(K,T,E_0)\qquad\,\,   \mbox{for $i=1,2,3$},
\end{equation}
for each open subset $V\subset K$.
Then we have
\begin{align}\label{7.23}
I_{2}^{\v,\d}+I_{4}^{\v,\d}+I_{5}^{\v,\d} \rightarrow \langle \mu_1^{\v}+\mu_2^{\v}+\mu_3^{\v},\, \phi\rangle
\quad  \mbox{as $\d\rightarrow0+$ (up to a subsequence)}.
\end{align}

For $I_3^{\v,\d}$, we notice from \eqref{5.7} that
\begin{align}
& \int_0^T\int_{K} \Big|\sqrt{\v}\rho^{\v,\d}  \partial_m\eta^{\v,\d}  \,  \big(u_r^{\v,\d}+\frac{N-1}{r} u^{\v,\d}\big) \Big|^{\frac43}\dd r\dd t\nonumber\\
&\leqq C_{\psi}(K)\int_0^T\int_{K} \Big|\sqrt{\v}\rho^{\v,\d}  (|u_r^{\v,\d}|+|u^{\v,\d}|) \Big|^{\frac43}\dd r\dd t\nonumber\\
&\leqq C_{\psi}(K) \Big(\v\int_0^T\int_{K}\big(\rho^{\v,\d}|u_r^{\v,\d}|^2+\rho^{\v,\d}|u^{\v,\d}|^2\big)\,\dd r\dd t\Big)^{\frac23}
\Big(\int_0^T\int_{K} |\rho^{\v,\d}|^2\dd r\dd t\Big)^{\frac13}\nonumber\\
&\leqq C_{\psi}(K,T,E_0).\nonumber
\end{align}
Then there exists a function $f^\v$ such that, as $\delta\rightarrow 0+$ (up to a subsequence),
\begin{align}
&\sqrt{\v}\rho^{\v,\d}  \partial_m\eta^{\v,\d}  \, \big(u_r^{\v,\d}+\frac{N-1}{r} u^{\v,\d}\big) \rightharpoonup f^\v
\quad  \mbox{weakly in $L^{\frac43}_{\rm loc}(\mathbb{R}^2_+)$},\label{7.25}\\
&\int_0^T\int_{K} |f^{\v}|^{\frac43}\,\dd r\dd t\leqq C_{\psi}(K,T,E_0). \label{7.26}
\end{align}
It follows from \eqref{7.25} that, as $\delta\rightarrow 0+$ (up to a subsequence),
\begin{align}\label{7.27}
I_3^{\v,\d}\rightarrow \sqrt{\v}\int_0^T\int_{K} f^{\v} \phi_r\,\dd r\dd t.
\end{align}

For $I_6^{\v,\d}$,
it follows from \eqref{5.7}--\eqref{5.8} and \eqref{7.6} that
\begin{align}\label{7.28}
|I_6^{\v,\d}|
\leqq&\, C_{\psi}({\rm supp}\,\phi ) \v\d \int_{\mathbb{R}^2_+}
\Big( (\rho^{\v,\d})^\alpha  \big(|u_r^{\v,\d}|+| u^{\v,\d}|\big) \phi_r
    +\big|(\rho^{\v,\d})^{\alpha-1}\rho^{\v,\d}_r u^{\v,\d}  \phi\big|\Big)\, \dd r\dd t\nonumber\\
\leqq&\, C_{\psi}({\rm supp} \,\phi) \v\d   \Big( \int_{\mathbb{R}^2_+} \big(\rho^{\v,\d}|u_r^{\v,\d}|^2+ \rho^{\v,\d}|u^{\v,\d}|^2\big)|\phi_r|\,\dd r\dd t\Big)^{\frac12}\nonumber\\
&\quad\,\,\times \Big( \int_{\mathbb{R}^2_+} (\rho^{\v,\d}+1)\,|\phi_r|\, \dd r\dd t\Big)^{\frac12} \nonumber\\
&\,+C_{\psi}({\rm supp}\,\phi)\sqrt{\d} \Big(\v^2\d\int_{\mathbb{R}^2_+}  (\rho^{\v,\d})^{\alpha-2} |\rho^{\v,\d}_r|^2\,|\phi|\, \dd r\dd t\Big)^{\frac12} \nonumber\\
&\quad\,\, \times
\Big(\int_{\mathbb{R}^2_+}  (\rho^{\v,\d})^{\alpha}|u^{\v,\d}|^2\,|\phi|\, \dd r\dd t\Big)^{\frac12}\nonumber\\
\leqq&\, C_{\psi}({\rm supp}\,\phi,\|\phi\|_{C^1},T,E_0)\sqrt{\d}\nonumber\\
&\quad\,\,\times \Big(\sqrt{\v}+
\big(\int_{\mathbb{R}^2_+}  \rho^{\v,\d}|u^{\v,\d}|^3\,|\phi|\, \dd r\dd t\big)^{\frac13}
\big(\int_{\mathbb{R}^2_+}  |\rho^{\v,\d}|^{3(\alpha-\frac23)}\,|\phi|\, \dd r\dd t\big)^{\frac16}\Big)\nonumber\\
\leqq&\, C_{\psi}({\rm supp}\,\phi,\|\phi\|_{C^1},T,E_0)\big(\sqrt{\v}+1\big) \sqrt{\d}
  \rightarrow0  \qquad  \mbox{as $\d\rightarrow0+$},
\end{align}
where we have used $\alpha=\frac{2N-1}{2N}\in [\frac{3}{4}, 1)$ for $N\geqq 2$.

\medskip
5. Taking $\d\rightarrow0+$ (up to a subsequence) on both sides of \eqref{7.7}, then it follows
from \eqref{7.16}, \eqref{7.20}, \eqref{7.23}, and \eqref{7.27}--\eqref{7.28}  that
\begin{align}\label{7.30}
\partial_t \eta^{\v}+\partial_r q^{\v}
=- \frac{N-1}{r} m^{\v} \big(\eta_\rho^{\v}+u^{\v} \eta_m^{\v}\big) +\mu_1^{\v}+\mu_2^{\v}+\mu_3^{\v}-\sqrt{\v} f^{\v}_r
\end{align}
in the sense of distributions. From \eqref{7.20-1}--\eqref{7.21-2}, we see that
\begin{align}\label{7.31}
- \frac{N-1}{r} m^{\v}\big(\eta_\rho^{\v}+u^{\v} \eta_m^{\v}\big) +\mu_1^{\v}+\mu_2^{\v}+\mu_3^{\v}
\end{align}
are bounded uniformly in $\v>0$ as Radon measures.
From \eqref{7.26}, we have
\begin{align}\label{7.32}
\sqrt{\v} f^{\v}_r \rightarrow 0  \qquad   \mbox{in $W^{-1,\frac43}_{\rm loc}(\mathbb{R}^2_+)$ as $\v\rightarrow0+$}.
\end{align}
Thus, it follows from \eqref{7.31}--\eqref{7.32} that
\begin{align}\label{7.33}
\partial_t \eta^{\v}+\partial_r q^{\v} \,\,  \mbox{is confined in a compact subset of $W^{-1, p_2}_{\rm loc}(\mathbb{R}^2_+)$ for some $p_2\in(1,2)$}.
\end{align}

The interpolation compactness theorem ({\it cf.} \cite{Chen2,Chen5}) indicates that, for $p_2>1$, $p_1\in(p_2,\infty]$, and $p_0\in[p_2, p_1)$,
\begin{align*}
&(\mbox{compact set of} \  W^{-1,p_2}_{\rm loc}(\mathbb{R}_+^2)) \cap (\mbox{bounded  set of} \  W^{-1,p_1}_{\rm loc}(\mathbb{R}_+^2))\\
&\subset (\mbox{compact set of} \  W^{-1,p_0}_{\rm loc}(\mathbb{R}_+^2)),
\end{align*}
which is a generalization of Murat's lemma in \cite{F. Murat,L. Tartar}.
Combining this interpolation compactness theorem for $1<p_2<2$, $p_1>2$, and $p_0=2$ with the facts in \eqref{7.17-1} and \eqref{7.33},
we conclude \eqref{7.1}.
$\hfill\Box$

\medskip
Combining Proposition \ref{prop5.2} with  Lemmas  \ref{lem5.8}--\ref{lem5.9} and \ref{lem7.1}, we have

\begin{theorem}\label{thm5.10}
Let $(\rho_0^\v, m_0^\v)$ be the initial data satisfying \eqref{1.53o}--\eqref{1.53b}.
For each $\v>0$, there exists a spherical symmetry weak solution
$$
(\rho^\v, \M^\v)(t,\mathbf{x}):=(\rho^\v(t,r), m^\v(t,r)\frac{\mathbf{x}}{r})
$$
of the compressible Navier-Stokes equations \eqref{1.1} in the sense of Definition {\rm \ref{definition-NS}}.
Moreover, $(\rho^\v, m^\v)(t,r)=(\rho^\v(t,r), \rho^\v(t,r) u^\v(t,r))$, with
$u^\v(t,r):=\frac{m^\v(t,r)}{\rho^\v(t,r)}$ a.e. on  $\{(t,r)\,:\,\rho^\v(t,r)\ne 0\}$ and $u^\v(t,r):=0$
a.e. on $\{(t,r)\,:\, \rho^\v(t,r)=0\,\, \mbox{or $\,r=0$}\}$,
satisfies the following bounds{\rm :}
\begin{align}
&\rho^\v(t,r)\geqq 0 \quad a.e., \nonumber\\[1mm]
&\big(\frac{m^\v}{\sqrt{\rho^\v}}\big)(t,r)=\sqrt{\rho^\v(t,r)}u^\v(t,r)=0\quad\mbox{a.e. on $\{(t,r)\,:\, \rho^\v(t,r)=0\}$},
\nonumber\\[1.5mm]
&\int_0^\infty  \Big(\frac12 \Big|\frac{m^{\v}}{\sqrt{\rho^\v}}\Big|^2 +e(\rho^{\v},\bar{\r}) \Big)(t,r)\,r^{N-1} \dd r
 +\v\int_{\mathbb{R}_+^2} \Big|\frac{m^{\v}}{\sqrt{\rho^\v}}\Big|^2(s,r)\, r^{N-3}\dd r \dd s\nonumber\\
 &\leqq E_0^\v\leqq E_0+1\,\,\,\mbox{for $t>0$},\label{5.55}\\[2mm]
&\v^2 \int_0^\infty\big|\big(\sqrt{\rho^{\v}(t,r)}\big)_r\big|^2\, r^{N-1}\dd r
  +\v\int_{\mathbb{R}_+^2} \big|\big((\rho^{\v}(s,r))^{\frac{\gamma}{2}}\big)_r\big|^2\,r^{N-1}\dd r\dd s\nonumber\\
&\leqq C (E_0+1)\qquad\mbox{for $t>0$},\nonumber
\\[2mm]
&\int_0^T\int_{d}^D (\rho^{\v})^{\g+1}(t,r)\,\dd r\dd t\leqq C(d,D,T,E_0),\label{5.57}\\[1mm]
&\int_0^T\int_0^D   \big(\rho^\v |u^\v|^3+(\rho^{\v})^{\gamma+\theta}\big)(t,r)\,r^{N-1}\dd r\dd t
\leqq C(D,T, E_0) \label{5.58}
\end{align}
for any fixed $T>0$ and any compact subset $[d,D]\Subset (0,\infty)$.

Let $(\eta, q)$ be an entropy pair defined in \eqref{weakentropy}  for a smooth compact supported function $\psi(s)$ on $\mathbb{R}$.
Then, for $\v\in (0,1]$,
\begin{align*}
\partial_t\eta(\rho^\v,m^\v)+\partial_rq(\rho^\v,m^\v) \qquad  \mbox{is compact in $H^{-1}_{\rm loc}(\mathbb{R}^2_+)$}.\nonumber
\end{align*}
\end{theorem}

\medskip
\section{Proof of the Main Theorems}
In this section, we give a complete proof of Main Theorem II: Theorem \ref{thm1.1},
which leads to Main Theorem I: Theorem \ref{thm:2.1}, as indicated in Remark \ref{remark:2.5a}.

The uniform estimates and compactness properties obtained in Theorem \ref{thm5.10}
imply  that the weak solutions $(\rho^{\v}, m^{\v})$ of the Navier-Stokes equations \eqref{1.3}
satisfy the compensated compactness framework in Chen-Perepelitsa \cite{Chen6}.
Then the compactness theorem established in \cite{Chen6} for the case $\gamma>1$ (also see LeFloch-Westdickenberg \cite{Ph. LeFloch} for $\gamma\in (1,5/3]$)
implies that there exist functions $(\rho,m)(t,r)$ such that
\begin{align*}
(\rho^\v,m^{\v})\rightarrow (\rho,m) \qquad  \mbox{{\it a.e.}\, $(t,r)\in \mathbb{R}^2_+$\, as $\v\rightarrow0+$ (up to a subsequence)}.
\end{align*}

By similar arguments as in the proof of Lemma 5.6, we find that
$m(t,r)=0$ {\it a.e.} on $\{(t,r) \,:\, \rho(t,r)=0\}$.
We can define the  limit velocity $u(t,r)$ by setting $u(t,r):=\frac{m(t,r)}{\rho(t,r)}$ {\it a.e.} on $\{(t,r)\,:\,\rho(t,r)\neq0\}$
and $u(t,r):=0$ {\it a.e.} on $\{(t,r)\,: \rho(t,r)=0\,\,\, \mbox{or $r=0$}\}$.
Then we have
\begin{align}
m(t,r)=\rho(t,r) u(t,r).\nonumber
\end{align}
We can also define $(\frac{m}{\sqrt{\rho}})(t,r):=\sqrt{\rho(t,r)} u(t,r)$, which is $0$ a.e on the vacuum states $\{(t,r) \ :\ \rho(t,r)=0\}$.
Moreover, we obtain that, as $\v\rightarrow0+$,
\begin{equation}\label{7.0-9}
\frac{m^{\v}}{\sqrt{\rho^{\v}}}\equiv\sqrt{\rho^{\v}} u^{\v}\rightarrow \frac{m}{\sqrt{\rho}}\equiv\sqrt{\rho} u \quad
\mbox{strongly in } \ L^2([0,T]\times[0,D], r^{N-1}\dd r \dd t).
\end{equation}

Notice that
$|m|^{\frac{3(\gamma+1)}{\gamma+3}}\leqq C\big(\frac{|m|^3}{\rho^2}+\rho^{\gamma+1} \big)$,
which,  along with \eqref{5.57}--\eqref{5.58}, implies
\begin{align}\label{7.51}
(\rho^\v,m^{\v})\rightarrow (\rho,m) \qquad \mbox{in $L^p_{\rm loc}(\mathbb{R}^2_+)\times L^q_{\rm loc}(\mathbb{R}^2_+)$ as $\v\rightarrow 0+$},
\end{align}
for $p\in[1,\gamma+1)$ and $q\in[1,\frac{3(\gamma+1)}{\gamma+3})$,
where $L^q_{\rm loc}(\mathbb{R}^2_+)$ represents $L^q([0,T]\times K)$ for any $T>0$ and $K\Subset (0,\infty)$.

From the same estimates, we also obtain the convergence of the relative mechanical energy as $\v\rightarrow0+$:
\begin{align}\nonumber
\bar{\eta}^{\ast}(\rho^\v,m^\v)\rightarrow \bar{\eta}^{\ast}(\rho,m)\qquad \mbox{in $L^1_{\rm loc}(\mathbb{R}_+^2)$}.
\end{align}
Since $\bar{\eta}^{\ast}(\rho,m)$ is a convex function, by passing the limit in \eqref{5.55},
we have
\begin{align}\nonumber
\int_{t_1}^{t_2} \int_0^\infty\bar{\eta}^{\ast}(\rho,m)(t,r)\,  r^{N-1}\dd r\dd t
 \leqq (t_2-t_1) \int_0^\infty\bar{\eta}^{\ast}(\rho_0,m_0)(r)\,  r^{N-1}\dd r,
\end{align}
which implies
\begin{align}\label{7.83}
\int_0^\infty\bar{\eta}^{\ast}(\rho,m)(t,r)\, r^{N-1}\dd r\leqq  \int_0^\infty\bar{\eta}^{\ast}(\rho_0,m_0)(r)\, r^{N-1}\dd r
\quad\, \mbox{for {\it a.e.} $t\geqq0$}.
\end{align}
This indicates that there is no concentration formed in the density $\rho(t,r)$ at origin $r=0$.

\medskip
Define
\begin{align}\label{3.113-1}
(\rho,\M)(t,\mathbf{x}):=(\rho(t,r), m(t,r) \frac{\mathbf{x}}{r})=(\rho(t,r), \rho(t,r) u(t,r)\frac{\mathbf{x}}{r}).
\end{align}
From \eqref{7.83}, we know that $\frac{\M}{\sqrt{\rho}}=\sqrt{\rho}u\,\frac{\mathbf{x}}{r}$ is well-defined and in $L^2$ for {\it a.e.} $t>0$.
We now prove that $(\rho,\M)(t,\mathbf{x})$ is a weak solution of the Cauchy problem
for the Euler equations \eqref{1.1-1} in $\mathbb{R}^N$.

\smallskip
Let $\zeta(t,\mathbf{x})\in C_0^1([0,\infty)\times \mathbb{R}^N)$ be a smooth, compactly supported function.
Then it follows from \eqref{5.39} that
\begin{align}\label{7.60}
\int_{\mathbb{R}^{N+1}_+}\big(\rho^{\v} \zeta_t + \M^{\v}\cdot\nabla \zeta\big)\, \dd\mathbf{x}\dd t
 +\int_{\mathbb{R}^N}\rho_0^\v(\mathbf{x}) \zeta(0,\mathbf{x})\, \dd\mathbf{x}=0.
\end{align}

Let $\phi(t,r)$ be the corresponding function defined in \eqref{3.39-1}.
Using \eqref{7.51} and similar arguments as in the proof of Lemma \ref{lem5.8},
we obtain that, for any fixed $\sigma>0$,
\begin{align}\label{3.116}
&\lim_{\v\rightarrow0+}\int_0^\infty \int_{\mathbb{R}^N\backslash B_\s(\mathbf{0})}
  \big(\rho^{\v} \zeta_t + \M^{\v} \cdot\nabla\zeta\big)\, \dd\mathbf{x}\dd t\nonumber\\
&=\omega_N\lim_{\v\rightarrow0+}\int_0^\infty\int_\sigma^\infty \big(\rho^{\v} \phi_t + m^{\v} \phi_r\big)\, r^{N-1}\dd r\dd t\nonumber\\
&=\omega_N\int_0^\infty\int_\sigma^\infty \big(\rho \phi_t + m \phi_r\big)\, r^{N-1}\dd r\dd t\nonumber\\
&=\int_0^\infty \int_{\mathbb{R}^N\backslash B_\s(\mathbf{0})} \big(\rho\zeta_t + \M\cdot\nabla\zeta\big)\, \dd\mathbf{x}\dd t.
\end{align}
Using \eqref{7.83} and by similar arguments as in \eqref{5.42-1}, we have
\begin{align}
&\Big|\int_{0}^{\infty} \int_{B_\s(\mathbf{0})}  (\rho^{\v}-\rho) \zeta_t \,\dd\mathbf{x}\dd t\Big|\nonumber\\
&\leqq C(\|\zeta\|_{C^1}, {\rm supp}\,{\zeta})\Big\{\int_0^\infty\int_0^\sigma \big((\rho^{\v})^{\gamma}+\rho^{\gamma}\big)
 \,|\phi_t|\, r^{N-1}\dd r\dd t \Big\}^{\frac1\gamma}
 \sigma^{N(1-\frac1\gamma)}
 \nonumber\\
&\leqq C(\|\zeta\|_{C^1},{\rm supp}\,{\zeta}, E_0)\, \sigma^{N(1-\frac1\gamma)}\rightarrow 0\qquad  \mbox{as $\sigma\rightarrow0$},\label{3.117a}\\[2mm]
&\Big|\int_{0}^{\infty} \int_{B_\s(\mathbf{0})}  (\M^{\v}-\M) \cdot\nabla \zeta\, \dd\mathbf{x}\dd t\Big|\nonumber\\
&\leqq C\Big\{\int_{0}^{\infty}\int_0^\sigma \Big(\frac{|m^\v|^2}{\rho^{\v}}+\frac{m^2}{\rho}\Big)(t,r)\, |\phi_r|\, r^{N-1}\dd r\dd t\Big\}^{\frac12}\nonumber\\
&\qquad \times
\Big\{\int_{0}^{\infty}\int_0^\sigma \big(\rho^{\v}+\rho\big)(t,r)\, |\phi_r|\, r^{N-1} \dd r\dd t\Big\}^{\frac12}
\nonumber\\
&\,\,\,\leqq C(\|\zeta\|_{C^1},{\rm supp}\,{\zeta},E_0)\sigma^{\frac{N}{2}(1-\frac1\gamma)}\rightarrow 0\qquad  \mbox{as $\sigma\rightarrow0$}, \label{3.117}
\end{align}
which, with \eqref{3.116}--\eqref{3.117}, implies
\begin{align}\label{3.118}
\lim_{\delta\rightarrow0+}\int_{\mathbb{R}^{N+1}_+}
 \big(\rho^{\v} \zeta_t + \M^{\v} \cdot\nabla\zeta\big)\, \dd\mathbf{x}\dd t
 =\int_{\mathbb{R}^{N+1}_+}
 \big(\rho \zeta_t + \M\cdot\nabla\zeta\big)\, \dd\mathbf{x}\dd t.
\end{align}
Letting $\v\rightarrow0+$ in \eqref{7.60} and using \eqref{3.118}, we conclude that
$(\rho, \M)$ satisfies \eqref{1.18}.

Next we consider the momentum equations.
Let $\psi=(\psi_1,\cdots,\psi_N) \in \big(C_0^2(\mathbb{R}\times\mathbb{R}^N)\big)^N$ be a smooth function
with compact support, and let $\chi_\sigma(r)\in C^\infty(\mathbb{R})$ be a cut-off function satisfying \eqref{3.49-1}.
Without loss of generality, we assume that ${\rm supp}\, \psi \subset [-T,T]\times B_D(\mathbf{0})$.
Denote $\Psi_{\sigma}=\psi \chi_\sigma$. Then we have
\begin{align}\label{7.62}
&\bigg|\v\int_{\mathbb{R}^{N+1}_+}
\Big\{\frac{1}{2}\M^{\v}\cdot \big(\Delta \Psi_{\sigma}+\nabla\mbox{div}\Psi_{\sigma} \big)
+\frac{\M^{\v}}{\sqrt{\rho^{\v}}} \cdot \big(\nabla\sqrt{\rho^{\v}}\cdot \nabla\big)\Psi_{\sigma}\nonumber\\
&\qquad\qquad\,\,\,
 + (\nabla\sqrt{\rho^{\v}})\cdot \big(\frac{\M^{\v}}{\sqrt{\rho^{\v}}}\cdot \nabla\big)\Psi_{\sigma}\Big\}\, \dd\mathbf{x} \dd t\bigg|\nonumber\\
&=\left|\sqrt{\v}\int_{\mathbb{R}^{N+1}_+}
  \sqrt{\rho^{\v}}
  \Big\{V^{\v}  \frac{\mathbf{x}\otimes\mathbf{x}}{r^2}+\frac{\sqrt{\v}}{r}\frac{m^\v}{\sqrt{\rho^\v}}
    \big(I_{N\times N}-\frac{\mathbf{x}\otimes\mathbf{x}}{r^2}\big)\Big\}: \nabla\Psi_{\sigma}\, \dd\mathbf{x} \dd t\right|\nonumber\\
&\leqq C\Big\{\int_{\mathbb{R}^{N+1}_+}
|V^\v|^2\,\dd\mathbf{x}
  +\v\int_{\mathbb{R}_+^2} \frac{|m^\v|^2}{\rho^\v}\,r^{N-3}\dd r\dd t \Big\}^{\frac12}
\,\Big\{\v\int_{\mathbb{R}^{N+1}_+}
 \rho^\v |\nabla\Psi_{\sigma}|^2\, \dd\mathbf{x}\dd t\Big\}^{\frac12}\nonumber\\
&\leqq C(\sigma, D, T, E_0)\sqrt{\v}\rightarrow0\qquad\,\,  \mbox{as $\v\rightarrow0$}.
\end{align}
Using \eqref{7.0-9} and \eqref{7.62}, and passing the limit: $\v\rightarrow0+$ (up to a subsequence) in \eqref{3.59-1}, we obtain
\begin{align}\label{7.63}
&\int_{\mathbb{R}^{N+1}_+}
\Big\{\M\cdot\partial_t\Psi_\sigma
+ \frac{\M}{\sqrt{\rho}}\cdot \big(\frac{\M}{\sqrt{\rho}}\cdot \nabla\big)\Psi_{\sigma}
+p(\rho)\,\mbox{div} \Psi_{\sigma}\Big\}\,\dd\mathbf{x} \dd t\nonumber\\
&+\int_{\mathbb{R}^N} \M_0(\mathbf{x})\cdot \Psi_\s(0,\mathbf{x})\,\dd\mathbf{x}  =0.
\end{align}
Notice that, for any $T>0$ and $D>0$,
\begin{equation}\label{3.122}
\int_0^T\int_0^D \big(\frac{m^2}{\rho}+\rho^{\gamma}\big)(t,r)\, r^{N-1} \dd r\dd t\leqq C(D,T,E_0),
\end{equation}
which, with similar arguments as in \eqref{3.61}, leads to
\begin{align}\label{3.123}
&\lim_{\sigma\rightarrow0}\Big\{\int_{\mathbb{R}^{N+1}_+}
\M \cdot\partial_t\Psi_\sigma\,  \dd\mathbf{x}\dd t
+\int_{\mathbb{R}^N} \M_0(\mathbf{x})\cdot \Psi_\s(0,\mathbf{x})\,\dd\textbf{x}\Big\} \nonumber\\
&=\int_{\mathbb{R}^{N+1}_+}
\M \cdot\partial_t\psi\,  \dd\mathbf{x}\dd t
+\int_{\mathbb{R}^N} \M_0(\mathbf{x})\cdot \psi(0,\mathbf{x})\, \dd\mathbf{x}.
\end{align}
Using \eqref{3.63-2}--\eqref{3.63-1} and \eqref{3.122}, we have
\begin{align}
&\left|\int_{\mathbb{R}^{N+1}_+}
\big(\frac{m^2}{\rho}+p(\rho)\big) \big(\psi\cdot \frac{\mathbf{x}}{r}\big) \chi'_\sigma(r)
  \, \dd\mathbf{x}\dd t\right|\nonumber\\
&\leqq C\int_0^\infty\int_{\sigma}^{2\sigma} \big(\frac{m^2}{\rho}+p(\rho)\big) \varphi(t,r) |\chi'_\sigma(r)|\,r^{N-1}\dd r\dd t\nonumber\\
&\leqq C\int_0^T\int_{\sigma}^{2\sigma}  \big(\frac{m^2}{\rho}+p(\rho) \big)\, r^{N-1} \dd r\dd t\rightarrow 0\qquad\,\,  \mbox{as $\sigma\rightarrow0$},\nonumber
\end{align}
which,  with \eqref{3.122} and the Lebesgue dominated convergence theorem, leads to
\begin{align}\label{3.125}
&\lim_{\sigma\rightarrow0}\int_{\mathbb{R}^{N+1}_+}
\Big\{\frac{\M}{\sqrt{\rho}}\cdot \big(\frac{\M}{\sqrt{\rho}}\cdot \nabla\big)\Psi_{\sigma}+ p(\rho)\,\mbox{div} \Psi_{\sigma}\Big\}\,\dd\mathbf{x}\dd t\nonumber\\
&=\int_{\mathbb{R}^{N+1}_+}
\Big\{\frac{\M}{\sqrt{\rho}}\cdot \big(\frac{\M}{\sqrt{\rho}}\cdot \nabla\big)\psi+ p(\rho)\,\mbox{div} \psi \Big\}\,\dd\mathbf{x}\dd t.
\end{align}
Substituting \eqref{3.123}--\eqref{3.125} into \eqref{7.63}, we conclude
that
$(\rho, \M)$ satisfies \eqref{1.19}.

By the Lebesgue theorem, we can weaken the assumption: $\psi\in C_0^2$ as $\psi\in C_0^1$.
This completes the proof of Theorem \ref{thm1.1}.  $\hfill\Box$

\appendix
\section{Construction and Estimates of Approximate Initial Data}

In this appendix, we construct the approximate initial data functions with desired estimates and regularity.
From \eqref{initial}, we know that there exists a constant $R\gg1$ so that
\begin{align}\label{8.1}
0<\frac12 \bar{\rho}\leqq \rho_0(r)\leqq \frac32 \bar{\rho}\qquad\,\, \mbox{for $r\geqq R$}.
\end{align}
We first cut-off the density function $\rho_0(r)$ as
\begin{align}\label{8.2}
\tilde{\rho}_0^\v(r)=
\begin{cases}
(\beta \v)^{\frac14} \qquad &\mbox{if $\rho_0(r)\leqq (\beta\v)^{\frac14} $},\\
\rho_0(r)  &\mbox{if $(\beta \v)^{\frac14} \leqq \rho_0(r)\leqq (\beta \v)^{-\frac12} $},\\
(\beta \v)^{-\frac12}  &\mbox{if $\rho_0(r)\geqq (\beta \v)^{-\frac12} $},
\end{cases}
\end{align}
where $\v\in(0,1]$, and $0<\beta\ll1$ is a given small positive constant, which
is used to ensure $(\beta \v)^{\frac14} \ll (\beta \v)^{-\frac12} $ for all $\v\in(0,1]$. It is easy to check that
\begin{align}\label{8.3-1}
\tilde{\rho}_0^\v(r)\leqq \rho_0(r)+1,
\qquad\tilde{\rho}_0^\v(r)\rightarrow \rho_0(r) \,\,\,\mbox{as $\v\to 0\,\,$ {\it a.e.} $r\in \mathbb{R}_+$}.
\end{align}

To keep the $L^p$--properties of mollification, it is more convenient to smooth out the initial data
in the original coordinate $\mathbb{R}^N$;
so we do not distinguish between functions $(\rho_0,m_0)(r)$ and $(\rho_0,m_0)(\mathbf{x})=(\rho_0,m_0)(|\mathbf{x}|)$
when no confusion arises.

\medskip
It follows from \eqref{1.20a}--\eqref{6.10} that $\rho_0(\mathbf{x})\in L^\gamma_{\rm loc}(\mathbb{R}^N)$.
Using the convexity of $e(\rho,\bar{\rho})$, we have
\begin{align}\label{8.4}
0\leqq e(\tilde{\rho}_0^\v(\mathbf{x}),\bar{\rho})\leqq e(\rho_0(\mathbf{x}),\bar{\rho}).
\end{align}

Combining \eqref{1.20a}
with \eqref{8.3-1}--\eqref{8.4} and the Lebesgue dominated convergence theorem, we obtain
\begin{align}
\lim_{\v\rightarrow0+}  \int_{K}\Big(\big|\tilde{\rho}_0^\v(\mathbf{x})-\rho_0(\mathbf{x})\big|^{\gamma}
+\big|\sqrt{\tilde{\rho}_0^\v(\mathbf{x})}-\sqrt{\rho_0(\mathbf{x})}\big|^{2\gamma}\Big)\,\dd\mathbf{x}=0
\label{8.4-1}
\end{align}
for any $K\Subset\mathbb{R}^N$.

Since we need a better decay property for approximate initial data,
we further cut-off the function $\tilde{\rho}_0^\v(\mathbf{x})$ at the far-field:
\begin{align}\label{8.2-1}
\hat{\rho}_0^\v(\mathbf{x})=
\begin{cases}
\tilde{\rho}_0^\v(\mathbf{x}) \quad &\mbox{if $|\mathbf{x}|\leqq (\beta \v)^{-\frac1{2N}} $},\\[1mm]
\bar{\rho} \quad &\mbox{if $|\mathbf{x}|> (\beta \v)^{-\frac1{2N}} $},
\end{cases}
\end{align}
where we further choose $\beta$ small enough so that $|\mathbf{x}|\geqq (\beta \v)^{-\frac1{2N}} \geqq R+2$ for all $\v\in(0,1]$.
It is clear that  $\hat{\rho}_0^{\v}(\mathbf{x})$ is not a smooth function so that
we need to  mollify $\hat{\rho}_0^{\v}(\mathbf{x})$.
Let  $J(\mathbf{x})$ be the standard mollification function and
$J_{\sigma}(\textbf{x}):=\frac{1}{\sigma^N}J(\frac{\mathbf{x}}{\sigma})$ for $\sigma\in(0,1)$.
For later use, we take $\sigma=\v^{\frac14}$ and define  $\rho_{0}^{\v}(\mathbf{x})$  as
\begin{align}\label{8.6}
\rho_{0}^{\v}(\mathbf{x})
:=\Big(\int_{\mathbb{R}^N} \sqrt{\hat{\rho}_{0}^\v(\mathbf{x-y})} J_\sigma(\mathbf{y})\,\dd\mathbf{y}\Big)^2.
\end{align}
Then $\rho_{0}^{\v}(\mathbf{x})$ is still a spherically symmetric function,
{\it i.e.}, $\rho_{0}^{\v}(\mathbf{x})=\rho_{0}^{\v}(|\mathbf{x}|)$.

\begin{lemma}\label{lem8.2}
For any given $\v\in(0,1]$, $\rho_0^\v(\mathbf{x})$ defined in \eqref{8.6}
is in $C^\infty(\mathbb{R}^N)$
with $(\beta \v)^{\frac14}\leqq \rho^\v_0(\mathbf{x})\leqq (\beta \v)^{-\frac12}$
and satisfies
\begin{align}
&\lim_{\v\rightarrow0+}\Big(\big\| \rho^\v_0-\rho_0\big\|_{L^\gamma_{\rm loc}(\mathbb{R}^N)}
  +\Big|\int_{\mathbb{R}^N} \big(e(\rho_0^{\v}(\mathbf{x}),\bar{\rho})- e(\rho_{0}(\mathbf{x}),\bar{\rho})\big)\,\dd\mathbf{x}\Big|\Big)=0,\label{8.8}\\[1mm]
&\v^2\int_{\mathbb{R}^N} \big|\nabla_{\mathbf{x}} \sqrt{\rho_0^{\v}(\mathbf{x})}\big|^2\,\dd\mathbf{x}
\leqq C\sqrt{\v},\label{8.7}\\[1mm]
&\int_{\mathbb{R}^N} e(\rho_0^{\v}(\mathbf{x}),\bar{\rho}) (1+|\mathbf{x}|)^{N-1+\vartheta}\,\dd\mathbf{x}
 \leqq CE_0\v^{-\frac{N-1+\vartheta}{2N}},\label{8.9}
\end{align}
where $E_0$ is  defined in \eqref{1.53}, and $\vartheta\in (0,1)$.
\end{lemma}

\noindent{\bf Proof.}  We divide the proof into four steps.

\smallskip
1. We first consider the first part of \eqref{8.8}. A direct calculation shows
\begin{align}\label{8.8-2}
&\big|\sqrt{\rho_0^\v(\mathbf{x})}-\sqrt{\rho_0(\mathbf{x})}\big|
\nonumber\\
&\leqq \Big|\int_{\mathbb{R}^N}\big(\sqrt{\hat{\rho}_0^\v(\mathbf{x-y})}
     -\sqrt{\rho_0(\mathbf{x-y})}\big) J_{\sigma}(\mathbf{y})\,\dd\mathbf{y}\Big|\nonumber\\
&\quad+\Big|\int_{\mathbb{R}^N}\big(\sqrt{\rho_0(\mathbf{x-y})}-\sqrt{\rho_0(\mathbf{x})}\big) J_{\sigma}(\mathbf{y})\,\dd\mathbf{y}\Big|.
\end{align}
For any given $M\gg1$, it follows from \eqref{8.8-2} and the H\"{o}lder inequality that
\begin{align}\label{8.8-3}
&\int_{|\mathbf{x}|\leqq M+1} \big|\sqrt{\rho_0^\v(\mathbf{x})}-\sqrt{\rho_0(\mathbf{x})}\big|^{2\gamma}\dd\mathbf{x}\nonumber\\
&\leqq C\int_{\mathbb{R}^N} J_{\sigma}(\mathbf{y}) \int_{|\mathbf{x}|\leqq M+1}
  \Big(\big|\sqrt{\hat{\rho}_0^\v(\mathbf{x-y})}-\sqrt{\rho_0(\mathbf{x-y})}\big|^{2\gamma}\nonumber\\
&\qquad \qquad\qquad\qquad\qquad\quad\,\,\, +\big|\sqrt{\rho_0(\mathbf{x-y})}-\sqrt{\rho_0(\mathbf{x})}\big|^{2\gamma}\Big)\,\dd\mathbf{x}\dd\mathbf{y}\nonumber\\
&\leqq C\sup_{|\mathbf{y}|\leqq \v^{\frac14}}\big\|\sqrt{\rho_0(\cdot+\mathbf{y})}-\sqrt{\rho_0(\cdot)}\big\|_{L^{2\gamma}(\{|\mathbf{x}|\leqq M+1\})} \nonumber\\[1mm]
&\quad + C\big\|\sqrt{\tilde{\rho}_0^\v}-\sqrt{\rho_0}\big\|_{L^{2\gamma}(\{|\mathbf{x}|\leqq M+2\})}
\rightarrow 0
\end{align}
as $\v\rightarrow0+$,
where we have used \eqref{8.4-1}, $\sigma=\v^{\frac14}$, and $\hat{\rho}_0^\v(\mathbf{x})=\tilde{\rho}_0^\v(\mathbf{x})$ for
$|\mathbf{x}|\leqq (\beta\v)^{-\frac1{2N}}$.
Using \eqref{8.8-3}, it is direct to obtain
\begin{equation}\label{8.8-4}
\int_{|\mathbf{x}|\leqq M+1} |\rho_0^\v(\mathbf{x})-\rho_0(\mathbf{x})|^{\gamma}\,\dd\mathbf{x}\rightarrow 0\qquad\,  \mbox{as $\v\rightarrow0$}.
\end{equation}

\smallskip
2. We now consider the second part of  \eqref{8.8}. For any given $M\gg 1$, it follows from \eqref{8.8-4} that
\begin{align}\label{A.16}
\lim_{\v\rightarrow0+}\int_{|\mathbf{x}|\leqq M+1} \big(e(\rho_0^{\v}(\mathbf{x}),\bar{\rho})- e(\rho_0(\mathbf{x}),\bar{\rho})\big)\,\dd\mathbf{x}=0.
\end{align}
For $|\mathbf{x}|>M+1$ with $M\geqq R+1$, noting \eqref{8.1}--\eqref{8.2} and \eqref{8.2-1}--\eqref{8.6}, we have
\begin{align}\label{A.17}
0<\frac12\bar{\rho}\leqq \rho_{0}^\v(\mathbf{x})\leqq \frac{3}{2} \bar{\rho}.
\end{align}
It follows from \eqref{8.2} and \eqref{8.2-1} that
$\big|\sqrt{\hat{\rho}_0^\v(\mathbf{x})}-\sqrt{\bar{\rho}}\big|\leqq \big|\sqrt{\tilde{\rho}_0^\v(\mathbf{x})}-\sqrt{\bar{\rho}}\big|$
for $\mathbf{x}\in\mathbb{R}^N$,
which,  with \eqref{A.17}, yields
\begin{align}\label{A.18}
&\int_{|\mathbf{x}|>M+1} e(\rho^\v_0(\mathbf{x}), \bar{\rho})\,\dd\mathbf{x}\nonumber\\
&\leqq C(\bar{\rho})\int_{|\mathbf{x}|>M+1} \big|\sqrt{\rho_0^\v(\mathbf{x})}-\sqrt{\bar{\rho}}\big|^2\dd\mathbf{x}
\nonumber\\
&  \leqq C(\bar{\rho})\int_{|\mathbf{x}|>M+1}
 \Big|\int_{\mathbb{R}^N}\Big(\sqrt{\hat{\rho}_0^\v(\mathbf{x-y})}-\sqrt{\bar{\rho}}\Big) J_\sigma(\mathbf{y})\dd\mathbf{y}\Big|^2\dd\mathbf{x} \nonumber\\
&\leqq C(\bar{\rho})\int_{|\mathbf{x}|>M} \big|\sqrt{\hat{\rho}_0^\v(\mathbf{x})}-\sqrt{\bar{\rho}}\big|^2\,\dd\mathbf{x}\nonumber\\
&\leqq C(\bar{\rho})\int_{|\mathbf{x}|>M} \big|\sqrt{\tilde{\rho}_0^\v(\mathbf{x})}-\sqrt{\bar{\rho}}\big|^2\,\dd\mathbf{x}
\nonumber\\
&=C(\bar{\rho})\int_{|\mathbf{x}|>M} \big|\sqrt{\rho_0(\mathbf{x})}-\sqrt{\bar{\rho}}\big|^2\,\dd\mathbf{x}\nonumber\\
&\leqq C(\bar{\rho})\int_{|\mathbf{x}|>M} e(\rho_0(\mathbf{x}), \bar{\rho})\,\dd\mathbf{x}.
\end{align}
For any given small $\varrho>0$, there exists $M(\varrho)\gg 1$  such that
\begin{align}\label{A.19}
\int_{|\mathbf{x}|>M(\varrho)} e(\rho_0(\mathbf{x}), \bar{\rho})\,\dd\mathbf{x}\leqq \varrho.
\end{align}
Using \eqref{A.16} and \eqref{A.18}--\eqref{A.19}, we have
\begin{align*}
&\left|\int_{\mathbb{R}^N} \big(e(\rho_0^{\v}(\mathbf{x}),\bar{\rho})-e(\rho_{0}(\mathbf{x}),\bar{\rho})\big)\,\dd\mathbf{x}\right|\\
&\leqq \left|\int_{|\mathbf{x}|\leqq M(\varrho)+1} \big(e(\rho_0^{\v}(\mathbf{x}),\bar{\rho})- e(\rho_{0}(\mathbf{x}),\bar{\rho})\big)\,\dd\mathbf{x}\right|\nonumber\\
&\quad +C(\bar{\rho})\int_{|\mathbf{x}|>M(\varrho)} e(\rho_0(\mathbf{x}), \bar{\rho})\, \dd\mathbf{x}\\
&\leqq C(\bar{\rho})\varrho,
\end{align*}
provided that $\v\ll 1$. Then \eqref{8.8} is proved.

\smallskip
\smallskip
3. Noting \eqref{8.2-1}, we have
\begin{align}
\partial_{\mathbf{x}_i} \sqrt{\rho_0^\v(\mathbf{x})}=
\begin{cases}
\displaystyle\int_{\mathbb{R}^N} \sqrt{\hat{\rho}_0^\v(\mathbf{x-y})} \partial_{\mathbf{y}_i} J_{\sigma}(y)\,\dd\mathbf{y}
  &\quad\mbox{for $|\mathbf{x}|\leqq 1+(\beta\v)^{-\frac1{2N}}$},\\[2.5mm]
\displaystyle 0  &\quad\mbox{for $|\mathbf{x}|\geqq 1+(\beta\v)^{-\frac1{2N}}$},\nonumber
\end{cases}
\end{align}
which, with \eqref{8.2} and \eqref{8.2-1},  leads to
\begin{align}
\v^2\int_{\mathbb{R}^N} \Big|\nabla_{\mathbf{x}} \sqrt{\rho_0^\v(\mathbf{x})}\Big|^2 \dd\mathbf{x}
&=\frac{C\v^2}{\sigma^{2}} \int_{|\mathbf{x}|\leqq 1+(\beta\v)^{-\frac1{2N}}} \sup_{\mathbf{y}\in\mathbb{R}^N} \hat{\rho}_0^\v(\mathbf{y})\,\dd\mathbf{x}\nonumber\\
&\leqq \frac{C\v^2}{\sigma^{2}}(\beta\v)^{-1}\leqq C \v^{\frac12},\nonumber
\end{align}
where we have used $\sigma=\v^{\frac14}$. Thus, \eqref{8.7} is proved.

\smallskip
4. We finally consider \eqref{8.9}.  Noting \eqref{8.2-1},
we see that $\rho^\v_0(\mathbf{x})=\bar{\rho}$ for all $|\mathbf{x}|\geqq 1+(\beta\v)^{-\frac1{2N}}$,
which, with \eqref{8.8}, implies
\begin{align*}
&\int_{\mathbb{R}^N} e(\rho_{0}^{\v}(\mathbf{x}),\bar{\rho}) (1+|\mathbf{x}|)^{N-1+\vartheta}\dd\mathbf{x}\nonumber\\
&=\int_{|\mathbf{x}|\leqq 1+(\beta\v)^{-\frac1{2N}}}e( \rho_{0}^{\v}(\mathbf{x}),\bar{\rho}) (1+|\mathbf{x}|)^{N-1+\vartheta}\dd\mathbf{x} \nonumber\\
&\leqq C\v^{-\frac{N-1+\vartheta}{2N}} \int_{|\mathbf{x}|\leqq 1+(\beta\v)^{-\frac1{2N}}}  e( \rho_{0}^{\v}(\mathbf{x}),\bar{\rho})\,\dd\mathbf{x}
\nonumber\\
&\leqq C(E_0+1) \v^{-\frac{N-1+\vartheta}{2N}}. \nonumber
\end{align*}
Therefore, we have proved \eqref{8.9}.
$\hfill\Box$

Denote $\mathbf{I}_{[4\delta,\delta^{-1}]}(\mathbf{x})$ to be the characteristic function
 $\{\mathbf{x} \in \mathbb{R}^N\, : \, 4\d\leqq |\mathbf{x}|\leqq \delta^{-1}\}$ with $0<\delta\ll1$.
Now, for the approximation of the velocity, we define $u^\v_0(\mathbf{x})$ and $u^{\v,\delta}_0(\mathbf{x})$:
\begin{align}
& u_0^\v(\mathbf{x}):=\frac{1}{\sqrt{\rho^\v_0(\mathbf{x})}}\big(\frac{m_0}{\sqrt{\rho_0}}\big)(\mathbf{x}),\label{A.21}\\
&u_0^{\v,\delta}(\mathbf{x}):=\frac{1}{\sqrt{\rho^\v_0(\mathbf{x})}}
\int_{\mathbb{R}^N} \big(\frac{m_0}{\sqrt{\rho_0}}\mathbf{I}_{[4\delta,\delta^{-1}]}\big)(\mathbf{x-y})J_\delta(\mathbf{y})\,\dd\mathbf{y},\label{A.22}
\end{align}
where $\rho^\v_0(\mathbf{x})$ is the approximate density function defined in Lemma \ref{lem8.2}.

\begin{lemma}\label{lem8.3}
The function $u_0^{\v}(\mathbf{x})$ defined in \eqref{A.21} satisfies
\begin{align}
&\int_{\mathbb{R}^N}\rho_0^\v(\mathbf{x}) |u_0^\v(\mathbf{x})|^2\,\dd\mathbf{x}
\equiv\int_{\mathbb{R}^N}\frac{|m_0(\mathbf{x})|^2}{\rho_0(\mathbf{x})}\,\dd\mathbf{x} \qquad\,\,\,\mbox{for any $\v\in (0,1]$}, \label{8.19}\\[2mm]
&\lim_{\v\rightarrow0+} \|\rho^\v_0 u^\v_0-m_0\|_{L^1_{\rm loc}(\mathbb{R}^N)}=0.\label{8.19-1}
\end{align}
The function $u_0^{\v,\d}(\mathbf{x})$ defined in \eqref{A.22} is in $C_0^\infty(\mathbb{R}^{N})$
and satisfies
\begin{align}
&{\rm supp}\,u_0^{\v,\d}\subset \{\mathbf{x} \in \mathbb{R}^N~ : ~ 2\d\leqq |\mathbf{x}|\leqq 1+\delta^{-1}\},\label{8.20o}\\[1mm]
&\lim_{\d\rightarrow0+} \int_{\mathbb{R}^N}\rho_0^{\v}(\mathbf{x}) |u_0^{\v,\d}(\mathbf{x})|^2\,\dd\mathbf{x}
= \int_{\mathbb{R}^N}\rho_0^\v(\mathbf{x}) |u_0^\v(\mathbf{x})|^2\,\dd\mathbf{x},\label{8.20}\\[1.5mm]
&\lim_{\d\rightarrow0+} \big\|\rho^\v_0 u^{\v,\d}_0-\rho^\v_0 u^{\v}_0\big\|_{L^1_{\rm loc}(\mathbb{R}^N)}=0,\label{8.20-1}\\[1mm]
&\int_{\mathbb{R}^N}\rho_0^{\v}(\mathbf{x}) |u_0^{\v,\d}(\mathbf{x})|^2(|\mathbf{x}|+1)^{N-1+\vartheta}\, \dd\mathbf{x}\leqq C E_0 \delta^{-N+1-\vartheta},\label{8.20-2}
\end{align}
where $E_0$ is defined in \eqref{1.53}.
\end{lemma}

\noindent{\bf Proof.}
\eqref{8.19} follows directly from \eqref{A.21}. Using \eqref{8.8-3} and \eqref{A.21}, we have
\begin{align}
&\int_{|\mathbf{x}|\leqq M}
\big|(\rho^\v_0 u^\v_0-m_0)(\mathbf{x})\big|\,\dd\mathbf{x}\nonumber\\
&=\int_{|\mathbf{x}|\leqq M}
\big|\big(\sqrt{\rho^\v_0}-\sqrt{\rho_0}\big)(\mathbf{x})(\frac{m_0}{\sqrt{\rho_0}})(\mathbf{x})\big|\,\dd\mathbf{x}\nonumber\\
&\leqq \Big(\int_{\mathbb{R}^N}\frac{|m_0(\mathbf{x})|^2}{\rho_0(\mathbf{x})}\,\dd\mathbf{x}\Big)^{\frac12}
\Big( \int_{|\mathbf{x}|\leqq M}
\big|\big(\sqrt{\rho^\v_0}-\sqrt{\rho_0}\big)(\mathbf{x})\big|^2\,\dd\mathbf{x}\Big)^{\frac12}\nonumber\\
&\rightarrow0 \qquad \mbox{as $\v\rightarrow0$}
\end{align}
for any $M\gg 1$, which leads to \eqref{8.19-1}.

From \eqref{A.22}, it is clear that $u_0^{\v,\d}(\mathbf{x})\in C_0^\infty(\mathbb{R}^{N})$
and  ${\rm supp}\,u_0^{\v,\d}\subset \{\mathbf{x} \in \mathbb{R}^N~ : ~ 2\d\leqq |\mathbf{x}|\leqq 1+\delta^{-1}\}$.
For any given small constant $\varrho>0$, there exist small $\epsilon=\epsilon(\varrho)>0$ and large $M=M(\varrho)\gg1$ such that
\begin{align}\label{8.23}
\int_{B_{2\epsilon}(\mathbf{0})\cup\{|\mathbf{x}|\geqq M(\varrho)\}} \frac{|m_0(\mathbf{x})|^2}{\rho_0(\mathbf{x})}\,\dd\mathbf{x}
\leqq \varrho.
\end{align}
Taking $\delta>0$ small enough so that $\epsilon\geqq 6\delta$, then it follows from  \eqref{A.22} that
\begin{align}\label{8.23-1}
\int_{\epsilon\leqq |\mathbf{x}|\leqq M+2} \Big|\Big(\sqrt{\rho_0^{\v}} u_0^{\v,\d}
- \frac{m_0}{\sqrt{\rho}_0}\Big)(\mathbf{x})\Big|^2\dd\mathbf{x}\rightarrow0\qquad\,\, \mbox{as $\d\rightarrow0+$}.
\end{align}
Since $\epsilon\geqq 6\d$, we have
\begin{align}\label{8.24}
&\int_{B_{\epsilon}(\mathbf{0})\cup\{|\mathbf{x}|\ge M+1\}} \big|\sqrt{\rho_0^{\v}(\mathbf{x})} u_0^{\v,\d}(\mathbf{x})\big|^2\,\dd\mathbf{x}\nonumber\\[1mm]
&\leqq \int_{B_{\epsilon}(\mathbf{0})\cup\{|\mathbf{x}|\ge M+1\}}
    \Big|\int_{\mathbb{R}^{N}} \big(\frac{m_0}{\sqrt{\rho_0}}\mathbf{I}_{[4\delta,\delta^{-1}]}\big)(\mathbf{x-y})
     J_\d(\mathbf{y})\,\dd\mathbf{y}\Big|^2\,\dd\mathbf{x}\nonumber\\[1mm]
 &  \leqq \int_{B_{2\epsilon}(\mathbf{0})\cup\{|\mathbf{x}|\ge M\}} \frac{|m_0(\mathbf{x})|^2}{\rho_0(\mathbf{x})}\,\dd\mathbf{x}\leqq \varrho.
\end{align}
It follows from \eqref{A.21} and \eqref{8.23}--\eqref{8.24} that
\begin{align}\label{8.26}
&\int_{\mathbb{R}^N} \big|\big(\sqrt{\rho_0^{\v}} u_0^{\v,\d}
- \sqrt{\rho_0^{\v}} u_0^{\v}\big)(\mathbf{x})\big|^2\,\dd\mathbf{x}\nonumber\\
&=\int_{\mathbb{R}^N} \big|\big(\sqrt{\rho_0^{\v}} u_0^{\v,\d}
- \frac{m_0}{\sqrt{\rho}_0}\big)(\mathbf{x})\big|^2\,\dd\mathbf{x}\nonumber\\
&\leqq \int_{\epsilon\leqq |\mathbf{x}|\leqq M+2} \big|\big(\sqrt{\rho_0^{\v}} u_0^{\v,\d}
- \frac{m_0}{\sqrt{\rho}_0}\big)(\mathbf{x})\big|^2\,\dd\mathbf{x}\nonumber\\
&\quad+C\int_{B_{2\epsilon}(\mathbf{0})\cup\{|\mathbf{x}|\geqq M\}}\frac{|m_0(\mathbf{x})|^2}{\rho_0(\mathbf{x})}\,\dd\mathbf{x}\nonumber\\
&\rightarrow 0\qquad
\mbox{as $\delta\rightarrow0+$},
\end{align}
which leads to \eqref{8.20}.

\smallskip
Using \eqref{8.26}, we have
\begin{align}
&\int_{|\mathbf{x}|\leqq M} \big|(\rho^\v_0 u^{\v,\d}_0-\rho^\v_0 u^{\v}_0)(\mathbf{x})\big|\,\dd\mathbf{x}\nonumber\\
&\leqq \Big( \int_{|\mathbf{x}|\leqq M} \rho^\v_0(\mathbf{x})\,\dd\mathbf{x}\Big)^{\frac12}
\Big(\int_{\mathbb{R}^N} \big|\big(\sqrt{\rho_0^{\v}} u_0^{\v,\d}
- \sqrt{\rho_0^{\v}} u_0^{\v}\big)(\mathbf{x})\big|^2\dd\mathbf{x}\Big)^{\frac12}\nonumber\\
\rightarrow 0 \qquad \mbox{as $\delta\rightarrow 0+$},
\nonumber
\end{align}
which implies  \eqref{8.20-1}.

\smallskip
Finally, noting \eqref{8.20} and $u_0^{\v,\d}(\mathbf{x})=0$ for $|\mathbf{x}|\geqq 1+\delta^{-1}$, we obtain
\begin{align}
&\int_{\mathbb{R}^N}\rho_0^{\v}(\mathbf{x}) |u_0^{\v,\d}(\mathbf{x})|^2(|\mathbf{x}|+1)^{N-1+\vartheta}\,\dd\mathbf{x}
\nonumber\\
&\leqq \int_{|\mathbf{x}|\leqq 1+\delta^{-1}}\rho_0^{\v}(\mathbf{x}) |u_0^{\v,\d}(\mathbf{x})|^2(|\mathbf{x}|+1)^{N-1+\vartheta}\,\dd\mathbf{x}\nonumber\\
&\leqq C\delta^{-N+1-\vartheta}\int_{\mathbb{R}^N}\rho_0^{\v}(\mathbf{x}) |u_0^{\v,\d}(\mathbf{x})|^2\,\dd\mathbf{x}\nonumber\\
&\leqq C\delta^{-N+1-\vartheta} \int_{\mathbb{R}^N}\frac{|m_0(\mathbf{x})|^2}{\rho_0(\mathbf{x})}\,\dd\mathbf{x}\nonumber\\
&\leqq C E_0 \delta^{-N+1-\vartheta},\nonumber
\end{align}
which yields \eqref{8.20-2}.
$\hfill\Box$

\medskip
With $\rho^{\v}_0(\mathbf{x})$, $u_0^{\v}(\mathbf{x})$,
and $u_0^{\v,\d}(\mathbf{x})$ defined above,
we can construct the approximate initial data
$(\rho_0^{\v,\d,b}, m_0^{\v,\d,b})(r)=(\rho_0^{\v,\d,b}, \rho_0^{\v,\d,b} u_0^{\v,\d,b})(r)$ for \eqref{6.1} and \eqref{6.3},
and $(\rho_0^{\v,\d}, m_0^{\v,\d})(r)=(\rho_0^{\v,\d}, \rho_0^{\v,\d}u_0^{\v,\d})(r)$ for \eqref{5.3}:
For $b\geqq 1+\delta^{-1}$,
define
\begin{align}\label{data-b}
(\rho_0^{\v,\d,b}, u_0^{\v,\d,b})(r):=(\rho_0^{\v}(\mathbf{x}), u_0^{\v,\d}(\mathbf{x}))\mathbf{I}_{[\d,b]}(\mathbf{x})
\qquad\,\, \mbox{for $r=|\mathbf{x}|\in[\d,b]$}
\end{align}
to be the initial data for IBVP \eqref{6.1} and \eqref{6.3}.
On the other hand, for IBVP \eqref{5.3}, we define
\begin{align}\label{data-d}
(\rho_0^{\v,\d}, u_0^{\v,\d}(r):=(\rho_0^{\v}(\mathbf{x}), u_0^{\v,\d}(\mathbf{x}))\mathbf{I}_{[\d,\infty)}(\mathbf{x})
\qquad\,\, \mbox{for $r=|\mathbf{x}| \in [\d,\infty)$}.
\end{align}
Then, combining Lemma \ref{lem8.2} with Lemma \ref{lem8.3}, we obtain

\begin{lemma}\label{lem8.4}
The following three results hold{\rm :}
\begin{enumerate}
\item[\rm (i)]  	As $\v\rightarrow0$,
\begin{align}\label{A.39a}
\begin{split}
&(E^\v_0, E_1^\v)\rightarrow (E_0,0), \\
&(\rho_0^{\v}, m_0^{\v})(r)\rightarrow (\rho_0, m_0)(r) \,\,\,\, \mbox{in $L_{\rm loc}^1([0,\infty); r^{N-1}\dd r)$},
\end{split}
\end{align}
where $E^\v_0, E_1^\v$, and $E_0$ are defined in \eqref{1.53}, \eqref{1.53a}, and \eqref{1.20a}, respectively.

\smallskip
\item[\rm (ii)]  For any fixed $\v\in(0,1]$, as $\d\rightarrow0$,
\begin{align}\label{A.39b}
\begin{split}
&(E_0^{\v,\d}, E_1^{\v,\d})\rightarrow (E_0^\v, E_1^\v),\\
&(\rho_0^{\v,\d}, m_0^{\v,\d})(r)\rightarrow (\rho_0^\v, m^\v_0)(r) \,\,\,\, \mbox{in $L_{\rm loc}^1([0, \infty); r^{N-1}\dd r)$},
\end{split}
\end{align}
where $E_0^{\v,\d}$ and $E_1^{\v,\d}$ are defined in \eqref{1.53-1a}--\eqref{1.53-1}.

\smallskip
\item[\rm (iii)] For any fixed $(\v,\d)$, as $b\rightarrow\infty$,
\begin{align}
&(E_0^{\v,\d,b}, E_1^{\v,\d,b})\rightarrow (E_0^{\v,\d}, E_1^{\v,\d}),\label{A.39c}\\
&(\rho_0^{\v,\d,b}, m_0^{\v,\d,b})(r)\rightarrow (\rho_0^{\v,\d}, m_0^{\v,\d})(r) \,\,\,
 \mbox{in $L^1_{\rm loc}((\delta, \infty); r^{N-1}{\rm d}r)$},
 \label{A39.d}
\end{align}
where $E_0^{\v,\d,b}, E_1^{\v,\d,b}, E_2^{\v,\d,b}$, and $\tilde{E}_0^{\v,\d,b}$
are defined
in Lemmas {\rm \ref{lem6.1}}--{\rm \ref{lem6.2}}
and \eqref{4.1d}.
In addition, the upper bounds of $E_0^{\v,\d,b}, E_1^{\v,\d,b}, E_2^{\v,\d,b}$, and
$\tilde{E}_0^{\v,\d,b}$
are independent of $b$ {\rm (}but may depend on $\v, \d${\rm )}, and
\begin{align}
& E_0^{\v,\d,b}+ E_1^{\v,\d,b}\leqq C(E_0+1), \label{A.39f}\\
&\tilde{E}_0^{\v,\d,b}\leqq \int_\delta^b\bar{\eta}^{\ast}(\rho_0^{\v,\delta,b},m_0^{\v,\d,b}) r^{N-1}(1+r)^{N-1+\vartheta}\, \dd r\nonumber\\
&\qquad\,\,\leqq C E_0 \big(\delta^{-N+1-\vartheta}+\v^{-\frac{N-1+\vartheta}{2N}}\big),\label{A.39e}
\end{align}
for some $C>0$ independent of $(\v,\d,b)$, where $\vartheta\in (0,1)$ is any fixed constant.
\end{enumerate}
\end{lemma}

\bigskip
\medskip
\noindent
{\bf Conflict of interest}: The authors declare that they have no conflict of interest.

\bigskip
\medskip
\noindent{\bf Acknowledgments.} The authors would like to thank Professor Didier Bresch
for helpful suggestions.
The research of Gui-Qiang G. Chen was supported in part by the UK
Engineering and Physical Sciences Research Council Awards
EP/L015811/1 and EP/V008854/1, and the Royal Society--Wolfson Research Merit Award (UK).
The research of Yong Wang was partially supported
by the National Natural Sciences Foundation of China No. 12022114, 11771429, 11671237, and  11688101.


\bigskip
\smallskip
\begin{thebibliography}{99}
	
\bibitem{Bressan-2005}
S. Bianchini and A.  Bressan,
Vanishing viscosity solutions of nonlinear hyperbolic systems.
{\it Ann. Math. (2)},  {\bf 161} (2005), 223--342.
	
	
\bibitem{BD-2003-CRMASP}
D. Bresch and B. Desjardins,
On viscous shallow-water equations (Saint-Venant model) and the quasi-geostrophic limit.
{\it C. R. Math. Acad. Sci. Paris}, {\bf 335} (2002), 1079--1084.

\bibitem{BD-2004-CRMASP}
D. Bresch and B. Desjardins,
Some diffusive capillary models of Korteweg type.
{\it C. R. Math. Acad. Sci. Paris, Section M\'{e}canique}, {\bf 332} (2004), 881--886.

\bibitem{BD-2007-JMPA}
D. Bresch and B. Desjardins,
On the existence of global weak solutions to the Navier-Stokes equations for viscous compressible and heat
conducting fluids.
{\it J. Math. Pures Appl.} {\bf 87} (2007), 57--90.



	
\bibitem{BD-2003-CPDE}
D. Bresch, B. Desjardins, and C.~K. Lin,
On some compressible fluid models: Korteweg, lubrication, and shallow water systems.
{\it Comm. Partial Diff. Eqs.} {\bf 28} (2003), 843--868.
	
\bibitem{BDG-2007}
D. Bresch, B. Desjardins, and D. Gerard-Varet,
On compressible Navier-Stokes equations with density dependent viscosities in bounded domains.
{\it J. Math. Pures Appl.} {\bf 87} (2007), 227--235.
	
\bibitem{BJ-2018}
D. Bresch and P.-E. Jabin,
Global existence of weak solutions for compressible Navier-Stokes equations: thermodynamically unstable
pressure and anisotropic viscous stress tensor.
{\it Ann. of Math.}  {\bf 188} (2018), 577--684.

\bibitem{Bressan-2004}
A. Bressan and  T. Yang,
On the convergence rate of vanishing viscosity approximations.
{\it Comm. Pure Appl. Math.} {\bf 57} (2004), 1075--1109.
	
\bibitem{Wang-2012-SIAM}
A. Bressan, F. Huang, Y. Wang, and T. Yang,
On the convergence rate of vanishing viscosity approximations
for nonlinear hyperbolic systems.
{\it SIAM J. Math. Anal.} {\bf 44} (2012), 3537--3563.
	
\bibitem{Chen1}
G.-Q. Chen, Convergence of the Lax-Friedrichs scheme
for isentropic gas dynamics (III).
{\it Acta Math. Sci.} {\bf 6B} (1986), 75--120 (in English);
{\bf 8A} (1988), 243--276 (in Chinese).
	
\bibitem{Chen3}
G.-Q. Chen,
Remarks on R. J. DiPerna's paper: Convergence of the viscosity method for
isentropic gas dynamics [Commun. Math. Phys. 91 (1983), 1--30].
{\it Proc. Amer. Math. Soc.} {\bf 125} (1997), 2981--2986.
	
\bibitem{Chen-1997}
G.-Q. Chen, Remarks on spherically symmetric solutions of the
compressible Euler equations.
{\it Proc. Roy. Soc. Edinburgh}, {\bf 127} (1997), 243--259.
	
\bibitem{Chen2}
G.-Q. Chen, The compensated compactness method and the
system of isentropic gas dynamics. Lecture Notes, Preprint
MSRI-00527-91, Berkeley, October 1990.
	

\bibitem{Chen-Feldman2018}
G.-Q. Chen and M.~Feldman.
\newblock {\em The Mathematics of Shock Reflection-diffraction and Von
  Neumann's Conjectures}.
Research Monograph, Annals of Mathematics Studies, {\bf 197},
Princeton University Press, Princeton, 2018.




\bibitem{Chen6}
G.-Q.  Chen and M. Perepelista,
Vanishing viscosity limit of the Navier-Stokes equations to the Euler
equations for compressible fluid flow.
{\it Comm. Pure Appl. Math.} {\bf 63} (2010), 1469--1504.

\bibitem{Chen7a}
G.-Q.  Chen and M. Perepelista,
Shallow water equations: viscous solutions and inviscid limit.
{\it Z. Angew. Math. Phys.} {\bf 63} (2012), 1067--1084.

\bibitem{Chen7}
G.-Q. Chen and M. Perepelista,
Vanishing viscosity solutions of the compressible Euler equations with spherical symmetry and large initial data.
{\it Commun. Math. Phys.} {\bf 338} (2015), 771--800.
	
\bibitem{Chen-Schrecker}
G.-Q. Chen and M. Schrecker,
Vanishing viscosity approach to the compressible Euler equations for transonic nozzle
and spherically symmetric flows.
{\it Arch. Ration. Mech. Anal.} {\bf 229} (2018), 1239--1279.
	
\bibitem{Courant-Friedrichs}
R. Courant and K.~O. Friedrichs,
{\it Supersonic Flow and Shock Waves}. Interscience Publishers Inc.:  New York, 1948.
	
\bibitem{Dafermos}
C.~M. Dafermos, {\it Hyperbolic Conservation Laws in Continuum
Physics}. Springer-Verlag: Berlin, 2010.
	
\bibitem{Ding} X. Ding,
On a lemma of DiPerna and Chen.
{\it Acta Math. Sci.} {\bf 26B}	(2006), 188--192.
	
\bibitem{Chen5} X. Ding, G.-Q. Chen, and P. Luo,
Convergence of the Lax-Friedrichs scheme for the isentropic gas dynamics
(I)-(II). {\it Acta Math. Sci.} {\bf 5B} (1985), 483--500, 501--540 (in English);
{\bf 7A} (1987), 467--480; {\bf 8A} (1989), 61--94 (in Chinese);
Convergence of the fractional step Lax-Friedrichs scheme and Godunov scheme for the
isentropic system of gas dynamics. {\it Commun. Math. Phys.} {\bf 121} (1989),
63--84.

\bibitem{R.J.DiPerna1}
R.~J. DiPerna, Convergence of the viscosity method for
isentropic gas dynamics.
{\it Commun. Math. Phys.} {\bf 91} (1983), 1--30.
	
\bibitem{R.J.DiPerna2}
R.~J. DiPerna,
Convergence of approximate solutions to conservation laws.
{\it Arch. Ration. Mech. Anal.} {\bf 82} (1983), 27--70.
	
	
\bibitem{Feireisl}	
E. Feireisl, A.  Novotny, and H.  Petzeltov\'{a},
On the existence of globally defined weak solutions to the Navier-Stokes equations.
{\it J. Math. Fluid Mech.} {\bf 3} (2001), 358--392.
	
\bibitem{D. Gilbarg}
D. Gilbarg, The existence and limit behavior of the
one-dimensional shock layer.
{\it Amer. J. Math.} {\bf 73} (1951), 256--274.
	
\bibitem{GMWZ} C.~M.~I.~O. Gu\`{e}s, G. M\`{e}tivier, M. Williams, and K. Zumbrun,
Navier-Stokes regularization of multidimensional Euler shocks.
{\it Ann. Sci. \'{E}cole Norm. Sup. $(4)$}, {\bf 39} (2006), 75--175.
	
\bibitem{Guderley}
G. Guderley, Starke kugelige und zylindrische Verdichtungsst\"{o}sse in der N\"{a}he des
Kugelmittelpunktes bzw. der Zylinderachse. Luftfahrtforschung, {\bf 19(9)} (1942), 302--311.


\bibitem{Guo-Jiu-Xin-2} Z.~H. Guo, Q.~S. Jiu, and Z.~P. Xin,
Spherically symmetric isentropic compressible flows with density-dependent viscosity coefficients.
{\it SIAM J. Math. Anal.} {\bf 39} (2008), 1402--1427.

\bibitem{D. Hoff1}
D. Hoff, Global solutions of the equations of one-dimensional, compressible flow
with large data and forces, and	with differing end states.
{\it Z. Angew. Math. Phys.} {\bf 49} (1998), 774--785.
	
\bibitem{D. Hoff-Liu}
D. Hoff and T.-P. Liu, The inviscid limit for the Navier-Stokes equations of compressible,
isentropic flow with shock data.
{\it Indiana Univ. Math. J.} {\bf 38} (1989), 861--915.

\bibitem{Huang-Wang}
F.~M. Huang and Z. Wang,
Convergence of viscosity solutions for isothermal gas dynamics.
{\it SIAM J. Math. Anal.} {\bf 34} (2002), 595--610.
	
\bibitem{Hugoniot}
P.~H. Hugoniot,
 M\`{e}moire sur la propagation du mouvement dans les corps et sp\'{e}cialement dans les gaz parfaits. 2e Partie.
 {\it J. \'{E}cole Polytechnique Paris}, {\bf 58} (1889), 1--125.
	
\bibitem{Song Jiang}
S. Jiang and P. Zhang,	
On spherically symmetric solutions of the compressible isentropic Navier-Stokes equations.
{\it Commun. Math. Phys.} {\bf 215} (2001), 559--581.
	

\bibitem{Kanel}
Ya. Kanel, On a model system of equations of one-dimensional gas motion.
{\it  Diff. Urav.} {\bf 4} (1968), 721--734 (in Russian).

\bibitem{P. D. Lax}
P.~D. Lax, Shock wave and entropy. In: {\it Contributions to
Functional Analysis}, ed. E.~A. Zarantonello, pp. 603--634,
Academic Press: New York, 1971.
	
\bibitem{Ph. LeFloch}
P. LeFloch and M. Westdickenberg,
Finite energy solutions to the isentropic Euler equations with
geometric effects. {\it J. Math. Pures Appl.} {\bf 88} (2007), 386--429.
	
\bibitem{Li-Wang}
T. Li and D. Wang, Blowup phenomena of solutions to the Euler equations for compressible
fluid flow. {\it J. Diff. Eqs.} {\bf 221} (2006), 91--101.

\bibitem{Lions-CNS}
P.-L. Lions, {\it Mathematical Topics in Fluid Dynamics, {\rm 2}: Compressible Models}.
Oxford Science Publication: Oxford, 1998.
	
\bibitem{Lions P.-L.1}
P.-L. Lions, B. Perthame, and P.~E. Souganidis,
Existence and stability of entropy solutions for the hyperbolic systems of
isentropic gas dynamics in Eulerian and Lagrangian coordinates.
{\it Comm. Pure Appl. Math.} {\bf 49} (1996), 599--638.
	
\bibitem{Lions P.-L.2}
P.-L. Lions, B. Perthame, and E. Tadmor,
Kinetic formulation of the isentorpic gas dynamics and p-systems.
{\it Commun. Math. Phys.} {\bf 163}	(1994), 415--431.
	

\bibitem{liu-xin-yang}
T.-P. Liu, Z. Xin, and T. Yang,
Vacuum states for compressible flow.
{\it Discrete Contin. Dynam. Systems}, {\bf  4} (1998), 1--32.
	
\bibitem{M-T}
T. Makino and S. Takeno,
Initial boundary value problem for the spherically symmetric motion of isentropic gas.
{\it Japan J. Ind. Appl. Math.} {\bf 11} (1994), 171--183.

\bibitem{MV}
A. Mellet and  A. Vasseur,
On the barotropic compressible Navier-Stokes equation.
{\it Comm. Partial Diff. Eqs.} {\bf 32} (2007), 431--452.

\bibitem{MRRS}
F. Merle, P. Raphael, I. Rodnianski, and J. Szeftel,
On the implosion of a three dimensional compressible fluid.
{\it arXiv Preprint}, arXiv:1912.11009v2, 2020.

	
\bibitem{C. Morawetz}
C. Morawetz,
An alternative proof of DiPerna's theorem.
{\it Comm. Pure Appl. Math.} {\bf 44} (1991), 1081--1090.

\bibitem{F. Murat}
F. Murat, Compacit\'e par compensation.
{\it Ann. Scuola Norm. Sup. Pisa Sci. Fis. Mat.} {\bf  5} (1978), 489--507.
	
\bibitem{B. Perthame}
B. Perthame and A. Tzavaras, Kinetic formulation for systems of two conservation
laws and elastodynamics.
{\it Arch. Ration. Mech. Anal.} {\bf  155} (2000), 1--48.
	
\bibitem{Plotnikov}
P.~I. Plotnikov and W. Weigant,
Isothermal Navier-Stokes equations and Radon transform.
{\it SIAM J. Math. Anal.} {\bf 47} (2015), 626--653.

\bibitem{Rankine}
W. J. M. Rankine, On the thermodynamic theory of waves of
finite longitudinal disturbance.
{\it Phil. Trans. Royal Soc. London}, {\bf 1960} (1870), 277--288.
	
\bibitem{Lord Rayleigh}
Lord Rayleigh (J. W. Strutt),
Aerial plane waves of finite amplitude.
{\it Proc. Royal Soc. London}, {\bf 84A} (1910), 247--284.

\bibitem{Rosseland}	
S. Rosseland, {\it The Pulsation Theory of Variable Stars}.
Dover Publications Inc.: New York, 1964.


\bibitem{Schrecker}
M. Schrecker,
Spherically symmetric solutions of the multidimensional compressible, isentropic Euler equations.
{\it Trans. Amer. Math. Soc.} {\bf 373} (2020), 727--746.

\bibitem{D. Serre1}
D. Serre,  La compacit\'{e} par compensation pour les syst\`{e}mes hyperboliques
non lin\'{e}aires de deux \`{e}quations \`{a} une dimension d'espace.
{\it J. Math. Pures Appl.} (9), {\bf 65} (1986), 423--468.
	
\bibitem{Slemrod}
M. Slemrod, Resolution of the spherical piston problem for compressible isentropic
gas dynamics via a self-similar viscous limit.
{\it Proc. Royal Soc. Edinburgh}, {\bf 126A} (1996), 1309--1340.
	
\bibitem{Stokes}
G.~G. Stokes,
On a difficulty in the theory of sound.
[{\it Philos. Mag.} {\bf 33} (1848), 349--356; {\it Mathematical and Physical Papers}, Vol. II, Cambridge Univ. Press: Cambridge, 1883].
Classic Papers in Shock Compression Science, 71--79.
{\it High-Pressure Shock Compression of Condensed Matter}. Springer: New York, 1998.
	
\bibitem{L. Tartar}
L. Tartar, Compensated compactness and applications to partial
differential equations.
In: {\it Research Notes in Mathematics, Nonlinear
	Analysis and Mechanics}, Herriot-Watt Symposium, Vol. 4, R.~J. Knops
	ed., Pitman Press, 1979.
	
\bibitem{VY}
A. Vasseur and C. Yu,  Existence of global weak solutions for 3D degenerate
compressible Navier-Stokes equations.
{\it Invent. Math.} {\bf  206} (2016), 935--974.
	
\bibitem{Whitham}
G. B. Whitham, {\it Linear and Nonlinear Waves}. Wiley: New York, 1974
	
\bibitem{Xin-1993}
Z.~P. Xin,  Zero dissipation limit to rarefaction waves for
the one-dimentional Navier-Stokes equations of compressible
isentropic gases.
{\it Comm. Pure Appl. Math.} {\bf 46} (1993), 621--665.
	
\bibitem{Yang-1}
T. Yang,  A functional integral approach to shock wave solutions of  Euler equations
with spherical symmetry I.
{\it Commun. Math. Phys.} {\bf  171} (1995), 607--638.

\bibitem{Yang-2}
T. Yang, A functional integral approach to shock wave solutions of Euler equations
with spherical symmetry II. {\it J. Diff. Eqs.} {\bf 130} (1996), 162--178.

\end{thebibliography}
\end{document}